%% file: main.tex
\documentclass[10pt,a4paper,reqno]{amsart}
\usepackage[utf8x]{inputenc}
\usepackage{ucs}
\usepackage{tikz}
\usepackage{amsmath}
\usepackage{amsthm}
\usepackage{amsfonts}
\usepackage{amssymb}
\usepackage{hyperref}
\usepackage{dsfont}
\usepackage{mathrsfs}
\usepackage{eqnarray}
\write18{}
\usepackage{xcolor}
\usepackage{graphicx}
\usepackage{pdfpages} 
\usepackage[nothm]{thmbox}
\usepackage{esint}

\hypersetup{
    colorlinks=true,
    linkcolor=red,
    filecolor=magenta,      
    urlcolor=cyan,
}

\usepackage[french,english]{babel}

\numberwithin{equation}{section}

\newcommand{\EE}{\mathbb{E}}

\newcommand{\PP}{\mathbb{P}}

\newcommand{\ZZ}{\mathbb{Z}}
\newcommand{\RR}{\mathbb{R}}
\newcommand{\sub}{\mathrm{sub}}

\newcommand{\NN}{\mathbb{N}}

\newcommand{\supp}[1]{\mathrm{supp}(#1)}
\newcommand{\limisup}[1]{\underset{#1}{\limsup}}
\newcommand{\limiinf}[1]{\underset{#1}{\liminf}}

\newcommand{\infi}[1]{\underset{#1}{\inf}}
\newcommand{\supr}[1]{\underset{#1}{\sup}}
\newcommand{\somme}[1]{\underset{#1}{\sum}}
\newcommand{\limi}[1]{\underset{#1}{\lim}}
\newcommand{\tends}[1]{\underset{#1}{\longrightarrow}}

\newcommand{\fonctionbis}[4]{\begin{array}{r c l}
					   #1 & \to & #2 \\
					 #3 & \mapsto & #4 \\
			   \end{array}}

\DeclareMathOperator{\CAT}{CAT}

\DeclareMathOperator{\Isom}{Isom}

\newcommand{\SL}{\operatorname{SL}}

\newcommand\inte{\operatorname{int}}

\newcommand{\efface}[1]{}

\newtheorem{theoreme}{Theorem}[section]
\newtheorem{proposition}[theoreme]{Proposition}
\newtheorem{corollary}[theoreme]{Corollary}

\newtheorem{lemma}[theoreme]{Lemma}

\theoremstyle{definition}

\newtheorem{definition}[theoreme]{Definition}
\newtheorem{remark}[theoreme]{Remark}
\newtheorem{exxample}[theoreme]{Example}

\numberwithin{equation}{section}

\title[Large deviations on Gromov-hyperbolic spaces]{Large deviations for  random walks on Gromov-hyperbolic spaces}
\author{Adrien Boulanger}
\address{Universit\`a di Bologna, Dipartimento di Matematica, Bologna, Italia}
\thanks{The first author was partially founded by the ERC n°647133 'IChaos'.}
\email{adrien.boulanger@unibo.it}
\author{Pierre Mathieu}
\address{Aix-Marseille Universit\'{e}, CNRS, Centrale Marseille, I2M, UMR 7373, 13453 Marseille, France}
\email{pierre.mathieu@univ-amu.fr}
\author{Cagri Sert}
\address{Institut f\"{u}r Mathematik, Universit\"{a}t Z\"{u}rich, 190, Winterthurerstrasse, 8057 Z\"{u}rich, Switzerland}
\email{cagri.sert@math.uzh.ch}
\thanks{The third author was supported by SNF grants 182089 and 193481.}
\author{Alessandro Sisto}
\address{Departement Mathematik, ETH Z\"{u}rich, 101, R\"{a}mistrasse, 8092 Z\"{u}rich, Switzerland}
\email{sisto@math.ethz.ch}


\begin{document}

\maketitle

\begin{abstract}
	Let $\Gamma$ be a countable group acting on a geodesic Gromov-hyperbolic metric space $X$ and $\mu$ a probability measure on $\Gamma$ whose support generates a non-elementary subsemigroup. Under the assumption that $\mu$ has a finite exponential moment, we establish large deviations results for the distance and the translation length of a random walk with driving measure $\mu$. From our results, we deduce a special case of a conjecture regarding large deviations of spectral radii of random matrix products.
\end{abstract}

\selectlanguage{french} 
\begin{abstract}
Soient $\Gamma$ un groupe dénombrable agissant sur un espace métrique géodesique hyperbolique $X$ et $\mu$ une mesure de probabilité sur $\Gamma$ dont le support engendre un semi-groupe non élémentaire. Sous l'hypothèse de moment exponentiel sur $\mu$, on établit des résultats de grandes déviations pour le déplacement et la longueur de translation d'une marche aléatoire suivant la loi $\mu$. Nous d\'{e}duisons de nos résultats un cas particulier d'une conjecture concernant les grandes déviations des rayons spectraux de produits de matrices al\'{e}atoires.
\end{abstract}

\selectlanguage{english}

\setcounter{tocdepth}{1}
\tableofcontents

\section{Introduction}
\label{secintro}

Let $\Gamma$ be an infinite, countable group acting by isometries on a metric space $(X,d)$, $\mu$ a probability measure on $\Gamma$ and $z_0 \in X$ a base point. A $(\mu,z_0)$-random walk on $X$, or random walk on $X$ for short, is the image under the orbital map $\gamma \mapsto \gamma \cdot z_0$ of the random walk on $\Gamma$ driven by the measure $\mu$. We denote with $(\gamma_n)_{n \in \NN} \in \Gamma^{\NN}$ (resp. $(z_n)_{n \in \NN} \in X^{\NN}$) the sequence of the successive positions of the walk on $\Gamma$ (resp. the sequence of the successive positions of the image random walk on $X$). We refer to Section \ref{subsectionrandomwalk} for basics on random walks. \\

We will say that `$\mu$ has a finite exponential moment' (resp.~finite first moment), if the random variable $d(z_0,z_1)$ has a finite exponential moment (resp.~finite first moment). In the sequel, $(\Omega,\mathbb{P})$ denotes the probability space on which the random walk is defined and $\mathbb{E}$ denotes the corresponding expectation. \\

For a probability measure $\mu$ with finite first moment, the rate of escape of the random walk is defined as the  limit
	\begin{equation}
		\label{eq-def rate of escape}
			 l :=  \underset{n \to \infty}{\lim} \ \frac{\EE(d(z_n,z_0))}{n} \ .
	\end{equation}
(The existence of the limit follows from sub-additivity.) It follows from Kingman's sub-additive ergodic theorem that $l$ is also the $\PP$ almost sure limit of the ratio $d(z_n,z_0)/n$. \\

This article addresses the question of large deviations with respect to this last convergence: we are looking for estimates of the probability that the distance $d(z_n,z_0)/n$ deviates from $l$ by an error of order $1$, either from below or from above, and similarly for the translation length $\tau(\gamma_n)/n$ (see below for definitions). More precisely, we investigate the case where the space $X$ is geodesic and Gromov-hyperbolic and the measure $\mu$ is \textbf{non-elementary}. A probability measure $\mu$ on $\Gamma$ is said to be non-elementary when its support generates a subsemigroup which contains two independent loxodromic elements; see Subsection \ref{subsectionyperbolicity}. Note we do not assume that $X$ is proper. \\

This setting has recently attracted a lot of attention as it encompasses several natural actions such as Gromov-hyperbolic groups acting on their Cayley graphs, rank-one semisimple groups acting on their symmetric spaces or Bruhat--Tits buildings, mapping class groups of surfaces acting on their curve complexes, relatively hyperbolic groups acting on their coned-off spaces, the Cremona group acting on the Picard-Manin hyperbolic space... We refer to the introduction of \cite[Section 1.2]{artmahertiozzo} for more details and references on the topic.\\

In \cite{artmahertiozzo}, \cite{artsunderl2017linear} and \cite{artmahertiozzocremona}, the authors investigate the escape rate of random walks driven by non-elementary measures. They show in particular that it is positive in this setting. Their approach focus on the boundary theory; they also manage to identify the Poisson boundary of the random walk with the Gromov boundary on the underlying Gromov-hyperbolic space under the assumption that the action is WPD. In \cite{artmathieusisto} a different approach was proposed based on deviation inequalities (and thus without any reference to boundary theory). Under the assumption that the action is acylindrical, the authors manage to prove a central limit theorem for the rate of escape on the group itself.

\subsection{Main results}
To formulate our results on large deviations of random walks on $X$, recall that a sequence $(Z_{n})_{n \in \NN}$ of real-valued random  variables is said to satisfy a \textit{large deviation principle}, abbreviated LDP from now on, if there exists a lower-semicontinuous function, called the rate function, $I: \RR \longrightarrow [0, \infty]$ such that for every measurable subset $R$ of $\RR$, we have 
\begin{equation}\label{eq.def.LDP}
\underset{\alpha \in \inte(R)}{-\inf I(\alpha)} \leq \underset{n \rightarrow \infty}{\liminf} \frac{1}{n}\ln \mathbb{P}(Z_{n} \in R) \leq \underset{n \rightarrow \infty}{\limsup} \frac{1}{n}\ln \mathbb{P}(Z_{n} \in R) \leq \underset{\alpha \in \overline{R}}{-\inf I(\alpha)}
\end{equation}
where $\inte(R)$ denotes the interior and $\overline{R}$ the closure of $R$. Our first main theorem is the following. 

\begin{theoreme}
\label{maintheo}
Let $\Gamma$ be a countable group acting by isometries on a geodesic Gromov-hyperbolic space $X$, $\mu$ a non-elementary probability measure on $\Gamma$ with finite exponential moment,  and $z_0 \in X$. Then the sequence of random variables $ ( \frac{1}{n} d(z_0, z_n))_{n \in \NN} $ satisfies a LDP with a proper convex rate function $I:[0,\infty) \to [0,\infty]$ which vanishes only at $l$. 
\end{theoreme}

Note first that the rate function $I$ does not depend on $z_0$ since the group acts by isometries. Indeed, for two different starting points $z_0$ and  $z_0'$, the difference $|d(\gamma_n \cdot z_0', z_0') - d(\gamma_n \cdot z_0, z_0)|$ is bounded by $2 \ d(z_0, z_0')$. Below, we list some more remarks on this result:
\begin{remark}
1. See Theorem \ref{thm.ldp.general} for a version of this result without any moment assumption on the probability measure $\mu$ and any hyperbolicity assumption on the metric space $X$.\\[3pt]
2. By convexity and lower-semicontinuity of $I$, the effective support of $I$, namely the set $D_I=\{\alpha \in [0,\infty) \, | \, I(\alpha)<\infty \}$ is an interval and $I$ is continuous on $D_I$. By Theorem \ref{maintheo}, this in turn implies that for every subset $J$ of $D_I$ satisfying $\overline{\inte(J)}=\overline{J}$ (e.g.\ any interval with non-empty interior), the limit $\lim_{n \to \infty} \frac{1}{n} \ln \mathbb{P}(\frac{1}{n}d(\gamma_n z_0,z_0) \in J)$ exists and is equal to $- \min_{\alpha \in \overline{J}}I(\alpha)$ (see Theorem \ref{thm.support} for more on $D_I$).\\[3pt]
3. The assumption that $\mu$ has a finite exponential moment is sharp regarding the conclusion that the rate function $I$ has unique zero (see Remark \ref{rk.upper.optimal} and also Remark \ref{rk.gouezel}). 
\end{remark}

To the best of our knowledge, exponential decay in large deviations and LDP's had not been studied in the context of Theorem \ref{maintheo} so far. Even in the special case where $\Gamma$ is Gromov-hyperbolic, Theorem \ref{maintheo} seems new. The most similar setting for which such a large deviation principle holds is for Lyapunov exponents associated to random products of matrices. We refer to the introduction of the third author's PhD thesis \cite{thesecagrisert} and the references therein for more details. In particular, in that setting, the proof of exponential decay in large deviations (corresponding to uniqueness of the zero of $I$) goes back to Le Page \cite{lepage}.\\



When $\Gamma$ is Gromov-hyperbolic and $\mu$ has a finite support, a possible alternative approach to prove that the rate function $I$ has unique zero, would be to exploit the spectral gap property of the image of the random walk on the boundary of the group. We refer to \cite[end of page 4]{artsebastiendriftanalitique}. For a surface group with the standard presentation and a driving measure with a finite exponential moment, large deviation estimates follow from the regeneration structure introduced in \cite{artblachmathieuhaissregeneration}.\\

Another important geometric notion of size associated to an isometry $\gamma$ acting on a Gromov-hyperbolic space $(X,d)$ is its translation length defined as
	$$\tau(\gamma) :=\infi{x \in X} \ d(x, \gamma \cdot x) \ .$$

This quantity has the advantage not to depend on a base point and is a conjugacy invariant. On the other hand, it is perhaps harder to study than $d(x, g\cdot x)$ since it is not sub-additive. For example, the lack of sub-additivity prevents one to readily get a convergence as in \eqref{eq-def rate of escape}. 
On the other hand, it is known that for a non-elementary probability measure with bounded support, the averages $\frac{1}{n}\tau(\gamma_n)$ and  $\frac{1}{n} d(z_n,z_0)$ behave similarly from the perspective of law of large numbers. Namely, they converge almost surely to the same constant $l$ (see e.g.\ \cite[Theorem 4.1]{artmahertiozzocremona}). \\

Let us now come to our second main theorem. We say that a set $\mathcal{B} \subset \Isom(X)$ is  \textbf{bounded} if 
	$$ \supr{g \in \mathcal{B}} \ d(x, g \cdot x) \ < \infty \ ,$$
is bounded for some $x \in  X$ (equivalently any). A probability measure $\mu$ on $\Isom(X)$ is said to have \textbf{bounded support} if its support is a bounded set. Our second main result reads


\begin{theoreme}
\label{theo.tau}
Let $\Gamma$ be a countable group acting by isometries on a geodesic Gromov-hyperbolic space $X$ and $\mu$ a non-elementary probability measure on $\Gamma$ of bounded support. Then the sequence of random variables 
	$ ( \frac{1}{n} \tau(\gamma_n))_{n \in \NN} $
satisfies a large deviation principle with the same rate function as the one given by Theorem \ref{maintheo}. 
\end{theoreme}

This theorem refines several previous results on the probabilistic behaviour of translation distance, e.g.\ \cite[Theorem 1.4]{artmahertiozzo}. For example, it  implies both the almost sure and the $L^1$ convergence $$ \frac{\tau(\gamma_n)}{n} \tends{ n \to \infty} l \ $$ (this was shown in \cite[Theorem 4.1]{artmahertiozzocremona}). In particular, specializing to the setting of the Cremona group, it also yields \cite[Theorem 1.2]{artmahertiozzocremona}. 
Another important consequence is expressed in Corollary \ref{corol.conjecture}. Namely, it confirms a special case of a conjecture about large deviations of spectral radii of random matrix products.\\ 

A common and sometimes more convenient way to express a notion of translation length is given by that of asymptotic translation length or stable length defined as 
\begin{equation}\label{eq.def.ell.gamma}
\ell(\gamma)= \limi{n \to \infty}\frac{d(x, \gamma^n \cdot x) }{n} \ .
\end{equation}
The limit exists by sub-additivity and does not depend on $x$. For a geodesic Gromov-hyperbolic space $X$, the difference $|\ell(.)-\tau(.)|$ is uniformly bounded on $\Isom(X)$ (see \cite[Ch.10, Prop. 6.4]{CDP}). Consequently, the previous theorem applies equally to the random variables $\frac{1}{n} \ell(\gamma_n)$ with the same conclusion. \\

The following subsections detail some direct consequences of the two above theorems and discuss some further properties of the rate function $I$. A complete description of the results of this article as well as its structure will be carried out in Section \ref{secoverview}.

\subsection{Properties of the rate function}\label{subsec.intro.properties}
A natural question motivated by the previous results concerns the understanding of the effective support $D_I=\{\alpha \in [0,\infty) \, | \, I(\alpha)<\infty \}$ of the rate function $I$. Note first that by convexity of the rate function $I$, the effective support $D_I$ is an interval in $[0,\infty)$. We denote by $l_{\min} := \inf D_I$ and by $l_{\max} := \sup D_I \in [0,\infty]$. For an equivalent definition of $l_{\min}$ and $l_{\max}$ without reference to a rate function, see \eqref{eq.lmax} and \eqref{eq.lmin}. \\

The function $I$ may be very degenerate. For example, let $\Gamma := \mathbb{F}_2 := \left<A, B \right>$ be the free group with two generators seen as acting on itself. We make it a metric tree $X$ by considering the word distance associated to the generating system $\{ A, A^{-1}, B, B^{-1} \}$ and we mark $z_0$ as the identity of $\Gamma$. Let then $\mu$ be the measure $\mu(A) = \mu(B) = \frac{1}{2}$. In this example the space $X$ is Gromov-hyperbolic and geodesic. The probability measure $\mu$ is supported by the set $\{A,B\}$ and, as such, has a finite exponential moment and generates a non-elementary subsemigroup ($A,B$ themselves are independent and loxodromic). In this case one has for all $n \in \NN$
	$$ d(z_0, z_n) = n \ , $$
so that the function $I$ has value $0$ at $1$ and $\infty$ otherwise. \\

For a boundedly supported probability measure, the function $I$ will be infinite on a neighbourhood of $\infty$ as well. However, it is easy to see that $l_{\min} = 0$ and $ l_{\mathrm{max}} > l$ whenever the subsemigroup generated by $\mu$ contains the identity. Indeed, we may accelerate or decelerate the random walk (with an exponential cost) by adjusting the frequency of 'identity elements' in the trajectories using an argument similar to the one used in the proof of  \cite[Theorem 4.12]{artmathieusisto}. Under more assumptions, one can even be more precise. \\

The following result gives a geometric characterization of $D_I$ only in terms of the support of the probability measure $\mu$. It also relates the effective support with the recently introduced notion of asymptotic joint displacement of a bounded set of isometries of a metric space. To state this result, we need some terminology. A set $\mathcal{B}$ of isometries of a metric space $(X,d)$ is said to be \textbf{non-arithmetic} if there exist $n \in \mathbb{N}$ and $g_1,g_2\in \mathcal{B}^n$ such that $\ell(g_1)\neq \ell(g_2)$. As in \cite{benoist-quint.hyperbolic}, we shall also call a probability measure non-arithmetic if its support is. \\

Let $\mathcal{B}$ be a subset of $\Isom(X)$. 
We call the following two quantities, respectively, asymptotic joint displacement (see \cite{breuillard-fujiwara,oregon.reyes:properties}) and lower asymptotic joint displacement:
\begin{equation}\label{eq.def.l.lsub}
\ell(\mathcal{B})=\lim_{n \to \infty}\sup_{g \in \mathcal{B}^n} \frac{1}{n} d(g\cdot x,x) \qquad \text{and} \qquad \ell_{\sub}(\mathcal{B}) = \lim_{n \to \infty} \ \inf_{ g \in \mathcal{B}^n} \ \frac{1}{n}d(g \cdot x,x).
\end{equation}
Both limits exist by subadditivity and they do not depend on $x$.

\begin{theoreme}[Effective support]\label{thm.support}
Let $\Gamma$ be a countable group acting by isometries on a geodesic Gromov-hyperbolic space $X$ and $\mu$ a non-elementary probability measure on $\Gamma$. Let $I$ be the rate function given by Theorem \ref{thm.ldp.general} (equivalently, by Theorem \ref{maintheo} if $\mu$ has a finite exponential moment). Then, 
\begin{equation*}
l_{\min} = \ell_{\sub}(\supp{\mu}) \qquad \text{and} \qquad l_{\max} = \ell(\supp{\mu}) \ ,
\end{equation*}
and the effective support $D_I$ of $I$ is an interval with non-empty interior (e.g.\ $l_{\min} \neq l_{\max}$) if and only if the probability measure $\mu$ is non-arithmetic. Moreover, if $\supp{\mu}$ is finite, then $D_I=[l_{\min},l_{\max}]$.
\end{theoreme}

\begin{remark}\label{rk.examples.intro}
In Subsection \ref{subsec.support}, we provide examples of probability measures $\mu$ of bounded (infinite) support for which the rate function $I$ explodes at $l_{\min}$ and $l_{\max}$.
\end{remark}

The notion of asymptotic joint displacement is analogous to the classical notion of joint spectral radius from linear algebra. In this geometric setting, it was recently studied by Oreg\'{o}n-Reyes \cite{oregon.reyes:properties} and Breuillard--Fujiwara \cite{breuillard-fujiwara} who proved the geometric analogues of some of the main results on joint spectral radius. The previous result parallels \cite[Theorem 1.7]{sert.LDP} where the effective support of the rate function of the norms of random matrix products was related to joint spectral radii. 

\subsection{Consequences for rank-one linear groups.}  Let us explain a consequence of our main theorem that partially answers a question raised in \cite{sert.LDP}. \\

A simple linear algebraic  group $H$ of rank one over a local field $k$ (e.g.\ $\SL_2(\mathbb{R})$ or $\SL_2(\mathbb{Q}_p)$), has a natural, up to finite index, faithful action by isometries on its symmetric space or the associated Bruhat--Tits tree $(X,d)$. The metric space $(X,d)$ is a Gromov-hyperbolic space. \\

One can find a finite-dimensional representation of $H$ such that for any $x \in X$ and $h \in H$, the displacement functional $d(x, h \cdot x)$ is given by the logarithm of the associated operator norm $\|.\|$ (see e.g.\ \cite[Chapter 6,8]{BQ.book} and \cite[\S 6]{quint.cones}). Moreover, the asymptotic translation length $\ell(h)$ corresponds to the logarithm of the spectral radius $\rho(h)$ of $h$, defined by the spectral radius formula $\rho(h)=\lim_{n \to \infty} \|h^n \|^{\frac{1}{n}}$. In this case, under the assumptions of Theorem \ref{maintheo}, the existence of a convex rate function for $\frac{1}{n}(d(z_n,z_0))$ follows from the main result of \cite{sert.LDP} (as well as, from Theorem \ref{maintheo}). \\

It was conjectured \cite[Conjecture 6.2]{sert.LDP} (see also \cite[\S 5.15]{breuillard-sert}) that if the support of the probability measure $\mu$ on $H$ generates a Zariski-dense subsemigroup (equivalently, if $\mu$ is non-elementary), then the sequence $\frac{1}{n}\ln \rho(\gamma_n)$ satisfies a LDP and  the rate function coincides with the rate function of the sequence $\frac{1}{n}\ln \|\gamma_n\|$. Under the assumption that the probability measure has finite support, this conjecture follows from Theorem \ref{theo.tau} for simple rank one groups.

\begin{corollary}\label{corol.conjecture}
Let $H$ be a simple linear algebraic  group of rank one over a local field $k$ endowed with an absolute value $|.|$. Let $\mu$ be a finitely generated probability measure on $H$ whose support generates a Zariski dense subsemigroup in $H$. Let $\|.\|$ be an operator norm on a finite-dimensional representation $V$ of $H$ as above and $I:[0,\infty) \to [0,\infty]$ be the rate function of the LDP of $\frac{1}{n}\ln \|\gamma_n\|$. Then, the sequence $\frac{1}{n}\ln \rho(\gamma_n)$ of random variables satisfies a LDP with rate function $I$.
\end{corollary}

The assumption that the support is finite may be replaced by the one that the measure has compact support. The authors decided not to write the article in this generality in order not to burden the proofs. 

\section{Detailed presentation of the article}
\label{secoverview} The article is mostly self-contained and  proofs only use a combination of elementary geometric and probabilistic arguments.  In particular, unlike in \cite{artmahertiozzo}, \cite{artmahertiozzocremona} or in \cite{artsunderl2017linear}, we make no use of any boundary whatsoever.


\subsection{Deviations from above and below}
In Section \ref{secrappel} we recall some basics on random walks, large deviation principles and hyperbolic geometry. As we shall see there, the proof of Theorem \ref{maintheo} (and Theorem \ref{thm.ldp.general} below) boils down to studying, the exponential decay and the limiting behaviour of the probabilities 
\begin{equation}\label{eq.lmax}
\frac{- \ln \left( \PP(  d(z_n, z_0)   \ge a n ) \right) }{n}
\end{equation} for every $a \in (l,l_{\max})$, and 
\begin{equation}\label{eq.lmin}
\frac{- \ln \left( \PP(  d(z_n, z_0)  \le  a n ) \right) }{n}
\end{equation}
for every $a \in (l_{\min},l)$; where $l_{\max}$ is defined as the infimum of $a$'s such that the limsup in \eqref{eq.lmax} is finite and similarly for $l_{\min}$. We refer to \eqref{eq.lmax} and \eqref{eq.lmin} as deviations from above and below, respectively. A thorough investigation of these is the overall objective of Sections \ref{secupper}, \ref{seclower}, \ref{seclinearprogress}, \ref{secintertwine}, \ref{secdevlinearprog} and \ref{sechittingmeasure}. \\

A very first observation is that under the finite exponential moment assumption, a general sub-additivity argument due to Hamana \cite{arthamadaupperdev}, that we recall in Appendix \ref{appendixA}, gives an upper bound on the probability of deviations from above: for any $\varepsilon > 0$, one has
\begin{equation}
	\label{eq:upp} 
	\limiinf{n \to \infty} \ \frac{ - \ln \left( \PP( d(z_n, z_0) - l n  \ge \varepsilon n ) \right) }{n} > 0 \ . 
\end{equation}

Inequality (\ref{eq:upp}) is very general; it holds for any group acting by isometries on any metric space. 

\begin{remark}\label{remarque1}
We observe that the exponential decay of the probability of a deviation from below cannot hold in the same generality as (\ref{eq:upp}). In the examples below, we equip a group $\Gamma$ with any left-invariant metric.  We choose $z_0=\mathrm{id}$ to be the identity element in $\Gamma$. We assume the rate of escape does not vanish for otherwise it makes no sense to compute deviations from below. \\[3pt]
1. Let $\Gamma$ be an amenable group and $\mu$ a symmetric probability measure with positive drift and whose finite support generates $\Gamma$ (see e.g.\ \cite{artkaimanovich1983}
). Then, Kesten's theorem implies that the probability $\PP(z_n=z_0)$ does not decay exponentially fast: 
$$-\frac 1 n \ln \PP(z_n=z_0) \tends{n \to \infty} 0 .$$ Therefore deviations from below have a sub-exponential decay. \\[3pt]
2. It is also possible to give examples of random walks on non-amenable groups for which deviations from below have a sub-exponential decay. Indeed start with an amenable group $\tilde\Gamma$ and a finitely supported symmetric driving measure $\tilde\mu$ as in 1. Then let $\Gamma$ be the direct product of $\tilde\Gamma$ with the free group on two generators $\mathbb F_2:=\left<A,B\right>$. Then $\Gamma$ is non-amenable. We endow $\Gamma$ with the metric $d=\tilde d + d^{(2)}$ given by a chosen metric $\tilde d$ on $\tilde\Gamma$ and the usual word metric $d^{(2)}$ on $\mathbb F_2$. Let $\mu$ be the product measure of $\tilde\mu$ on $\tilde\Gamma$ with the lazy simple random walk driving measure $\frac 12\delta_{\mathrm{id}}+\frac 18(\delta_A+\delta_{A^{-1}}+\delta_B+\delta_{B^{-1}})$ on $\mathbb F_2$. The two components of the random walk driven by $\mu$, say $(z_n)$, are then a random walk on $\Gamma$ driven by $\tilde\mu$ for the first component, say $(\tilde z_n)$ and a lazy simple symmetric random walk on $\mathbb F_2$ for the second component, say $(z^{(2)}_n)$. The two random walks $(\tilde z_n)$ and $(z^{(2)}_n)$ are independent. The rate of escape $l$ of the random walk $(z_n)$ is therefore the sum of the rate of escape of the random walk $(\tilde z_n)$ with respect to $\tilde d$, say $\tilde l$, and the rate of escape of the lazy simple random walk $(z^{(2)}_n)$, say $l^{(2)}$.  For any real $a$ such that $l^{(2)}<a<l=l^{(2)}+\tilde l$, we have that 
$$\PP(d(\mathrm{id},z_n)\leq an)\geq 
\PP(\tilde z_n=\mathrm{id})\PP(d^{(2)}(\mathrm{id},z^{(2)}_n)\leq an)\ .$$
As in example 1., the term $\PP(\tilde z_n=\mathrm{id})$ has a sub-exponential decay. Since $l^{(2)}<a$, the second term $\PP(d^{(2)}(\mathrm{id},z^{(2)}_n)\leq an)$ tends to $1$. Therefore $\PP(d(\mathrm{id},z_n)\leq an)$ has a sub-exponential decay. 
\end{remark}

Let us come back to the setting of Theorem \ref{maintheo}. We denote with $(y,z)_x$ the \textbf{Gromov product} of $y,z \in X$ with respect to $x$: \begin{equation}\label{eq.defn.gromov.product}  (y,z)_x := \frac 12 (d(y,x) + d(z,x) - d(y,z)).
\end{equation}
Our main geometric tool is the existence of a Schottky set as defined in the next 

\begin{definition}[Schottky set]
\label{defschottkyset}
Let $X$ be a metric space, $z_0 \in X$ and $S$ a non-empty finite subset of $\Isom(X)$. We say that $S$ is a \textbf{Schottky set} if there is a constant $C > 0$ such that for any pair $y,z \in X$ we have 
	$$ \frac{\# \ \{ s \in S \ , \ (y, s \cdot z)_{z_0} \le C \}}{\# \ S} \ge \frac{2}{3} \ . $$
\end{definition}

In Appendix \ref{appendixB}, we use a variation of the ping-pong lemma to prove that, when $X$ is Gromov-hyperbolic and geodesic and if the probability measure $\mu$ is non-elementary then there exists $p\in\NN$ such that the support of $\mu^{*p}$ contains a Schottky set. \\ 

We then deal separately with large deviations from above and from below. \\ 

As far as deviations from above are concerned, we already mentioned that 
the fact that a deviation from above has an exponentially small probability follows from 
Hamana's argument. 
In Section \ref{secupper}, we explain how the existence of the limit $\lim -\frac 1 n \ln \PP( d(z_n, z_0)\ge a n ) $ for all $a>l$
follows from a  sub-additivity argument.  In that argument, in order to compare 
$\PP( d(z_{n+m}, z_0)\ge a (n+m) )$ with the product $\PP( d(z_n, z_0)\ge a n )\PP( d(z_m, z_0)\ge a m )$, following  \cite{artdalbofeteke}, we use a Schottky set. We implement this approach using an insertion trick as in \cite{arthamanakesten1}. \\

Let us now discuss deviations from below. 
It is immediate, again by sub-additivity, that the limit $\lim -\frac 1 n \ln \PP( d(z_n, z_0)\le a n ) $ exists for all $a<l$ and defines a convex function; see Section \ref{seclower}. These already establishes the existence of LDP with a convex rate function for the sequence of random variables $\frac{1}{n}d(z_n,z_0)$ (see \S \ref{subsec.basics.LDP}). Regarding their large deviations, the hardest (and hopefully most interesting) part is to show that the limit is positive.  \\


Our starting point is a clever way to decompose a trajectory of a random walk that was introduced 
by A. Asselah and B. Schapira \cite{artshapiraasselah} to study large deviations for the range 
of random walks on $\mathbb Z^d$. Adapted to our context, it yields the following quite general criterion for deviations from below to be exponentially small.

\begin{proposition}\label{pasdelabel?}
Let $\Gamma$ be a countable group acting on a metric space $X$ and $\mu$ a probability measure on $\Gamma$. 
Then there is a convex function $\Psi : \ [ 0 , \infty ) \to [0,\infty]$ such that for all $a \neq l_{\min}$
\begin{equation*}
\frac{-\ln \PP \big( d(z_n,z_0) \le a n \big)  }{n}  \tends{n \to \infty} \Psi(a) \ .
\end{equation*}
Furthermore, if $\mu$ has a finite exponential moment and 
satisfies 
\begin{equation}\label{equationliminf}
\limiinf{p \to \infty} \ \supr{x \in X} \ \frac{\EE  \big( ( x, z_p)_{z_0} \big)}{p} = 0,
\end{equation}
then $\Psi$ vanishes only on $[l, \infty]$.
\end{proposition}
Proposition \ref{pasdelabel?} is proved in Section \ref{seclower}. Note that, in Proposition \ref{pasdelabel?}, we do not need assume $X$ is Gromov-hyperbolic or geodesic. 

\begin{remark}
	The above proposition can be more generally stated for defective adapted cocycles as defined in \cite{artmathieusisto}. However we restrain from doing so in order not to burden this article with many definitions.
\end{remark}

As a corollary of the previous proposition, we have the following 

\begin{corollary}\label{corollary1}
Let $\Gamma$ be a finitely generated amenable group and $\mu$ a symmetric finitely supported probability measure on $\Gamma$ whose support generates $\Gamma$. Equip $\Gamma$ with any left-invariant metric $d$. Assume the rate of escape does not vanish. Then 
	\begin{equation}\label{no-good} 	 \inf_{p} \ \supr{x \in X} \ \frac{\EE  \big( ( x, z_p)_{z_0} \big)}{p} \not= 0 \ .
	\end{equation}
\end{corollary}

Corollary \ref{corollary1} follows from Proposition \ref{pasdelabel?} and Kesten's theorem. As in Remark \ref{remarque1}, one also shows that there exist examples of random walks on non-amenable groups for which (\ref{equationliminf}) fails. \\


It now remains to show that  \eqref{equationliminf} holds in the setting of Theorem \ref{maintheo}. 
This will be a consequence of more precise exponential bounds on the tail of the law of the Gromov product $(z_n,x)_{z_0}$ stated in Proposition \ref{theoedpi} below. 

\subsection{LDP and walking-away uniformly on general metric spaces}

We start quantifying the rough idea that, given any point $x \in X$, with high probability, the random walk tends to walk away from $x$. The next Theorem \ref{theoremwau}  plays the central role in the proof of Proposition \ref{theoedpi}. It is proved in Sections \ref{seclinearprogress} and \ref{secintertwine}.

\begin{theoreme}[Walking-away uniformly] \label{theoremwau}
Let $\Gamma$ be a countable group acting by isometries on a metric space $X$, $\mu$ a probability measure on $\Gamma$ with a finite exponential moment and $z_0 \in X$. If the subsemigroup generated by $\mu$ contains a Schottky set and has unbounded orbits, then there is $\varepsilon, c_1, c_2 > 0$ such that for any $x \in X$ and all $n \in \NN$ we have
	$$ \PP( d(z_n,x) - d(z_0, x) \le \varepsilon n) \le c_1 \ e^{ -c_2 n} \ . $$
\end{theoreme} 

Note that we do not require $X$ to be Gromov-hyperbolic nor geodesic.

\begin{remark}
In the setting of Gromov-hyperbolic spaces, Theorem \ref{theoremwau} can be extracted from \cite{artsunderl2017linear} which builds on \cite{artmahertiozzo} and on ideas of \cite{artmathieusisto}. We however decided to give a short alternative proof to keep the article self-contained and use-of-boundary free. Moreover, the proof proposed here also adapts to the setting of a finite first moment to give an alternative proof of \cite[Theorem 1.1, Theorem 1.2]{artmahertiozzo}, see Subsection \ref{subsec.finitefirstmoment}.
\end{remark}

In fact, as we shall see, the analysis carried out so far allows us to get the following intermediary and general result which is weaker in conclusion but more general in assumptions (e.g.~no moment assumption on $\mu$ or Gromov-hyperbolicity assumption on $X$) in comparison to Theorem \ref{maintheo}. To state it, we introduce the following weakening of LDP which is relevant when the probability measures driving the random walk do not have a finite exponential moment: in \eqref{eq.def.LDP}, we say that the sequence $Z_n$ satisfies a \textit{weak LDP} if the lower bound holds for every measurable set $R$ and the upper bound holds for bounded measurable sets $R$. We have

\begin{theoreme}\label{thm.ldp.general}
Let $\Gamma$ be a countable group acting by isometries on a metric space $X$, $\mu$ a probability measure on $\Gamma$, and $z_0 \in X$. Suppose that the subsemigroup generated by the support of $\mu$ contains a Schottky set. Then,\\[2pt]
\hspace*{15pt} 1. the sequence $\frac{1}{n}d(z_n,z_0)$ satisfies a weak LDP with convex rate function $I:[0,\infty) \to [0,\infty]$.\\
\hspace*{15pt} 2. If, moreover, $\mu$ has a finite exponential moment and the subsemigroup generated by the support of $\mu$ has unbounded orbits, then the sequence $\frac{1}{n}d(z_n,z_0)$ satisfies an LDP, the rate function is proper and there exists $\varepsilon>0$ such that $I(x) >0$ for every $x \in [0,\varepsilon) \cup (l,\infty)$. 
\end{theoreme}
We note that the existence of $\varepsilon>0$ with the property that $I(x)>0$ for every $x \in [0,\varepsilon)$ directly follows from Theorem \ref{theoremwau} (see \S \ref{seclinearprogress}). \\

In view of the existence of Schottky sets in non-elementary semigroups (proved in Appendix \ref{appendixA}), taking the previous theorem for granted, to show Theorem \ref{maintheo}, what remains to be proven is that \textit{when $X$ is Gromov-hyperbolic} and $\mu$ is non-elementary, $\varepsilon$ can be taken to be the drift $l$ (which, as explained, we aim to achieve using Proposition \ref{pasdelabel?} by verifying \eqref{equationliminf}).

\begin{remark}\label{rk.gouezel}
By the discussion in Remark \ref{remarque1}, one cannot expect to get $\varepsilon=l$ in the generality of the previous theorem. However, for Gromov-hyperbolic spaces, after the appearance of a first version of this article, Gou\"{e}zel improved the moment aspect, by showing in the setting of Theorem \ref{maintheo} that $I(x)>0$ for every $x<l$ (i.e.\ $\varepsilon=l$) without the finite exponential moment assumption (see \cite{gouezel.positive}).
\end{remark}

Theorem \ref{theoremwau} in particular implies that the rate of escape does not vanish. More precisely, it implies the following linear progress with exponential tail property.

\begin{definition}[Linear progress]
\label{defslinearprogress}
Let $X$ be a metric space. We say that a random path $(z_n)$, with values in $X$, has \textbf{linear progress with exponential tail} if there is a constant $\varepsilon > 0$ such that 
	$$  \limiinf{n \to \infty} \ \frac{ - \ln \left( \PP( d(z_n, z_0)  \le \varepsilon n ) \right) }{n} > 0 \ . $$
\end{definition}

Note that for Gromov-hyperbolic spaces, the linear progress with exponential tail property was proved in \cite{artmahertiozzo} under the extra assumption that $\mu$ has finite support. 

\subsection{Exponential-tail and punctual deviations}

Sections \ref{secdevlinearprog} and \ref{sechittingmeasure} are devoted to deducing Proposition \ref{theoedpi} from the walking-away uniformly theorem. This proposition readily implies \eqref{equationliminf} and completes the proof of Theorem \ref{maintheo}. To prove Proposition \ref{theoedpi}, we shall rely on deviation inequalities. We start with the next result which is a variant of  \cite[Theorem 11.1]{artmathieusisto}. It is proved in Section \ref{secdevlinearprog}.

\begin{proposition}[exponential-tail deviation inequalities]\label{theoremdevineq}
Let $\Gamma$ be a countable group acting by isometries on a geodesic Gromov-hyperbolic space $X$, $\mu$ a non-elementary probability measure on $\Gamma$ with a finite exponential moment and $z_0 \in X$. If the random walk has linear progress with exponential tail, there are $ c_1, c_2 > 0$ such that for all $ 0 \le i \le n$ and all $R > 0$ one has 
	$$ \PP((z_n,z_0)_{z_i} \ge R) \le c_1 \ e^{ -c_2 R} \ .  $$
\end{proposition} 

In Section \ref{sechittingmeasure}, combining Proposition \ref{theoremdevineq} and the walking-away property from Theorem \ref{theoremwau}, we finally derive exponential bounds on the Gromov product $(z_n,x)_{z_0} $ as announced. 

\begin{proposition}[uniform punctual deviations]
\label{theoedpi}
Let $\Gamma$ be a countable group acting by isometries on a geodesic Gromov-hyperbolic space $X$ and $\mu$ a non-elementary probability measure on $\Gamma$.  Then, there are constants $C, \alpha >0$ such that for any $p \in \NN$ and any $x \in X$, $R>0$ we have 
	$$ \PP( (z_p,x)_{z_0} \ge R ) \le C e^{-\alpha R} \ .$$
\end{proposition} 

Integrating with respect to $R$ the bound in Proposition \ref{theoedpi}, one easily checks condition \eqref{equationliminf}. The proof of Theorem \ref{maintheo} is now complete. \\ 

We observe that, taking $n$ to $\infty$ in Proposition \ref{theoedpi}, we immediately derive bounds on the harmonic measure. 
We refer to Section \ref{secrappel} for all definitions regarding the next statement. 

\begin{corollary}[harmonic measure]
	\label{theoremharmonicmeasure}
Let $\Gamma$ be a countable group acting by isometries on a geodesic Gromov-hyperbolic space $X$, $\mu$ a non-elementary probability measure on $\Gamma$ with a finite exponential moment and $z_0 \in X$. There exists $D, C > 0$ such that for any $\zeta \in \partial X$ and any $r > 0$ the harmonic measure $\nu$ on $\partial X$ satisfies 
	$$ \nu(B(\zeta, r)) \le C \ r^D \ , $$
where $B(\zeta,r)$ stands for the ball (with respect to the Gromov metric) on $\partial X$ centred at $\zeta$ of radius $r$.
\end{corollary} 

Harmonic measures were studied in great detail for  proper Gromov-hyperbolic spaces (see for example \cite{artkiferharm,artkiferledrappharm,artblacherehaissinskymathieuharmonicmeasure,artbenoisthlinharm}). In particular the Hausdorff dimension of $\nu$ can then be computed and its multi-fractal spectrum described as in \cite{arttanakaharmo}. If $\Gamma$ is Gromov-hyperbolic and $\mu$ has a finite support, the inequality in Corollary \ref{theoremharmonicmeasure} holds when $D$ is replaced by the Hausdorff dimension \cite{artblacherehaissinskymathieuharmonicmeasure}. In our context of a more general action, an upper bound on the harmonic measure of a ball as in Corollary \ref{theoremharmonicmeasure} is proved in \cite{artmaherharm} but only when $\mu$ has a finite support.  

\subsection{LDP for translation length and support of the rate function}

Section \ref{sec.tau} is dedicated to the proof of Theorem \ref{theo.tau}. The proof uses Theorem \ref{maintheo} and can be split in two steps. \\

In a first part, using the existence of a Schottky set and an insertion trick in a similar way as in Section \ref{secupper}, we show that, given a prescribed speed $\alpha \geq l$ the event $\tau(\gamma_n) \ge \alpha n$ is, at the exponential scale, as likely as the event $d(z_0, z_n) \ge \alpha n $. \\

In the second part, for all prescribed speeds $ 0 \le \alpha < l$, we show that the event $\tau(\gamma_n) \le \alpha n$ is, at the exponential scale, as likely as $d(z_0, z_n) \ge \alpha n $. This step relies on Proposition \ref{prop.moving.tau} that uses an argument that finds, among the cyclic permutations of a given trajectory, a word whose displacement is uniformly close to the translation distance, which itself is invariant by cyclic permutation. \\ 

Section \ref{sec.support} is devoted to the proof of Theorem \ref{thm.support}. There, we also record some deterministic consequences of our results and the ingredients that we develop. For example, the following is a deterministic consequence of the combination of  Theorems \ref{maintheo}, \ref{theo.tau} and \ref{thm.support}. 

\begin{proposition}\label{prop.joint.spectrum}
Given a countable, bounded and non-elementary subset $\mathcal{B}$ of $\Isom(X)$, the  sequences of subsets $\frac{1}{n}d(\mathcal{B}^n \cdot z_0,z_0)$ and $\frac{1}{n}\tau(\mathcal{B}^n)$ of $\mathbb{R}$ converge to $[\ell_{\sub}(\mathcal{B}),\ell(\mathcal{B})]$ with respect to the Hausdorff metric.
\end{proposition}

In other words, the sequences $\frac{1}{n}d(\mathcal{B}^n \cdot z_0,z_0)$ and $\frac{1}{n}\tau(\mathcal{B}^n)$ become more and more dense in the interval $[\ell_{\sub}(\mathcal{B}),\ell(\mathcal{B})]$ as $n$ grows. In fact, Theorems \ref{maintheo} and \ref{theo.tau} can be seen as  quantitative refinements of this convergence.  \\

The previous proposition parallels the convergence result proven in \cite[Theorem 1.3]{breuillard-sert} for the vectors of singular values and moduli of eigenvalues of powers of a set of matrices. The interval $[\ell_{\sub}(\mathcal{B}),\ell(\mathcal{B})]$ corresponds to what is called the joint spectrum of $\mathcal{B}$ in that article. 

\subsection*{Acknowledgements} The authors are very grateful to  A. Asselah et B. Schapira who explained to them the strategy developed in \cite{artshapiraasselah} from which our Section \ref{seclower} is inspired. They would also like to thank Mathieu Dussaule and Peter Haissinsky for helpful conversations and Nguyen-Bac Dang for his explanations on the Picard-Manin space. Finally, the authors are also thankful to the anonymous referees for a number of corrections and suggestions that clarified the exposition of this article.


\section{First definitions and preliminary remarks}
\label{secrappel}

\subsection{Basics on random walks} 
\label{subsectionrandomwalk}
As a general reference on the topic, we recommend \cite{livrewoess,pete}. Let $\Gamma$ be an infinite, countable group and $\mu$ be a probability measure on $\Gamma$. 
Let $(\Omega, \PP)$ be a probability space and $(\omega_i)_{i \in \NN} : \Omega \to \Gamma$ a sequence of I.I.D.\ random variables following the law $\mu$. We call such a sequence the \textbf{increments} of the random walks. We then form the sequence of random variables 
 	$$ \gamma_n := \omega_1 \cdot \omega_2 \cdot ... \cdot \omega_n \ .$$
 	
Let $\Gamma$ act on a metric space $(X,d)$ with a marked point $z_0 \in X$. The push-forward of the random walk with respect to the orbital map is defined by 
	$$ \fonctionbis{\Gamma}{X}{\gamma}{\gamma \cdot z_0 \ .}  $$
We denote with
 	$$ z_n := \omega_1 \cdot \omega_2 \cdot ... \cdot \omega_n \cdot z_0 \ $$ the image of the sequence $(\gamma_n)$. 
We call $(z_n)$ \textbf{the positions} of the image random walk under the orbital map. 
We will often use the notation $d_n := d(z_0, z_{n})$ for short. 

\begin{remark} 
Note that the sequence of random variables $(z_n)_{n \in \NN}$ may not have the Markov property, even though the random walk 
$(\gamma_n)_{n \in \NN}$  is a Markov process. 
\end{remark}


Using $d(z_0, z_{m+n}) \le d(z_0, z_{n}) + d(z_n, z_{n+m})$ given by the triangle inequality and the fact that $d(z_n, z_{n +m})$ and $d(z_0, z_{m})$ have the same law, one deduces that the sequence $(\EE(d_n))_{n \in \NN}$ is sub-additive. Therefore, Fekete's lemma implies that the following limit exists
		$$ l := \limi{n \to \infty} \frac{\EE  ( d_n)}{n} = \infi{n \in \NN} \ \frac{\EE  (d_n)}{n}  \ . $$ 
We call $l$ the \textbf{rate of escape} of the image random walk. Note that Kingman's sub-additive ergodic theorem \cite{artkingmansubadd} (see also \cite{artsteelekingmans}) implies that the sequence $ \left( \frac{d_n}{n} \right)_{n \in \NN}$ also $\PP$-almost surely converges towards $l$.

\begin{remark}\label{rk.upper.optimal}
We observe that if one has a large deviations estimates as in Theorem \ref{maintheo}, then the measure $\mu$ has a finite exponential moment. 
	
Indeed the triangle inequality implies that, for any $a, \alpha \in \RR$, we have 
	$$ \PP(d_1 \ge \alpha n) \ \PP(d_{n-1} \le a n) \le \PP(d_n \ge (\alpha - a) n) \ . $$ 
In particular, for $a > l$ and $\alpha > 2 a$, the definition of $l$ imposes
	$$\PP(d_{n-1} \le a n) \to 1$$
 whereas the large deviations estimates from above imply that
	$$n \mapsto \PP(d_n \ge (\alpha - a) n)$$
 has an exponential decrease. Therefore,  the sequence $\left(\PP(d_1 \ge \alpha n)\right)_{n \in \NN}$ must also decrease exponentially fast.
\end{remark}

\subsection{Some preliminaries on large deviations theory}\label{subsec.basics.LDP} Here, we briefly justify that to prove the existence of limits in deviations from below and above is equivalent to the existence of the rate function in the language of large deviations theory. To keep the reading smooth, we postpone to Appendix \ref{appendix.rate} some further basic arguments in large deviations such as the explanation of how to identify the rate function using the limit Laplace generating function of the sequence $\frac{1}{n}d_n$.\\



Recall that for a sequence $(Z_n)_{n \in \mathbb{N}}$ of real-valued random variables, the definition of the large deviation principle (LDP) with a rate function $I:[0,\infty) \to [0,\infty]$ is given in \eqref{eq.def.LDP} and weak LDP is defined before Theorem \ref{thm.ldp.general}. The rate function is uniquely defined \cite[Lemma 4.1.4]{dembo-zeitouni}. We also introduce the notion of exponential tightness which, in our case, is an easy consequence of the finite exponential moment assumption (see Lemma \ref{lemma.proper.rate}).

\begin{definition}\label{def.exp.tight} A sequence $Z_n$ of real-valued random variables is said to be exponentially tight if for every $R>0$, there exists a compact set $K\subset \mathbb{R}$ such that $\liminf_{n \to \infty}-\frac{1}{n}\ln \mathbb{P}(Z_n \in K^c) \geq R$. 
\end{definition}

We have the following useful criterion for the existence of a LDP with a proper rate function. For its proof, see  \cite[Theorem   4.1.11 \& Lemma 1.2.18]{dembo-zeitouni}.

\begin{theoreme}\label{thm.LDP.criterion}
Let $Z_{n}$ be a sequence of real-valued random variables. Denote by $\mu_{n}$ the distribution of $Z_{n}$. For 
each $\alpha \in \mathbb{R}$, define:
$$
I_{li}(\alpha):= \sup_{\varepsilon > 0} - \underset{n \rightarrow \infty}{\liminf} \frac{1}{n}\ln  \mu_{n}((\alpha-\varepsilon,\alpha+\varepsilon)) \quad \text{and} \quad I_{ls}(\alpha):=\sup_{\varepsilon > 0} - \underset{n \rightarrow \infty}{\limsup} \frac{1}{n}\ln  \mu_{n}((\alpha-\varepsilon,\alpha+\varepsilon))
$$
Suppose that for all $\alpha \in \mathbb{R}$, we have $I_{li}(\alpha)=I_{ls}(\alpha)$. Then, the sequence $Z_{n}$ satisfies a weak LDP with the rate 
function $I$ given by $I(\alpha):=I_{li}(\alpha)=I_{ls}(\alpha)$. If, moreover, the sequence $Z_n$ is exponentially tight, then $I$ is proper and $Z_n$ satisfy LDP with the rate function $I$.
\end{theoreme}

Note that in the setting of Theorem \ref{thm.ldp.general}, the random variables $Z_n$ in the previous result are given by $\frac{1}{n}d(z_n,z_0)$. 

It is not hard to see that the hypotheses of the LDP criterion provided by the previous theorem boil down to the existence as limits of decay rates of one sided intervals. The following lemma spells out the precise conditions.

\begin{lemma}\label{lemma.Psi.and.I}
In the setting of Theorem \ref{thm.ldp.general}, suppose that\\[3pt] 
1)(deviations from above) for every $a \in (l,l_{\max})$, we have 
\begin{equation}\label{eq.from.above}
\lim_{n \to \infty} \frac{- \ln \left( \PP(  d(z_n, z_0)  \geq a n ) \right) }{n}=:\Psi(a) \in (0,\infty)
\end{equation} is a convex function of $a$ on $(l,l_{\max})$ tending to $0$ towards $l$,\\[3pt]
2)(deviations from below) for every $a \in (l_{\min},l)$, we have 
\begin{equation}\label{eq.from.below}
\lim_{n \to \infty} \frac{- \ln \left( \PP(  d(z_n, z_0)  \leq a n ) \right) }{n}=:\Psi(a) \in (0,\infty)
\end{equation} is a convex function of $a$ on $(l_{\min},l)$, tending to $0$ towards $l$.
Then, the sequence $\frac{1}{n}d(z_n,z_0)$ satisfies a weak LDP with the convex rate function $I:[0,\infty) \to [0,\infty]$ given by the extension of $\Psi$ by continuity to $[l_{\min},l_{\max}]$ and $I(a)=\infty$ for every $a \notin [l_{\min},l_{\max}]$. If, moreover, $\mu$ has finite exponential moment, then $\frac{1}{n}d(z_n,z_0)$ satisfy a LDP with the rate function $I$ which is proper.
\end{lemma}
\textbf{Proof.} It follows from Lemma \ref{lemma.proper.rate} that if $\mu$ has finite exponential moment, then the sequence $\frac{1}{n}d(z_n,z_0)$ is exponentially tight. The rest of the proof consists of a tedious verification that the hypotheses of this lemma imply those of Theorem \ref{thm.LDP.criterion} (namely that $I_{li}(\alpha)=I_{ls}(\alpha)$ for every $\alpha \geq 0$) and the extension of $\Psi$ have the common values of $I_{ls}$ and $I_{li}$. The details are straightforward and hence omitted for brevity.  \hfill $\blacksquare$\\

Finally, we mention the recent work of Corso \cite{corso} where the author proves the existence of LDP with a convex rate function for random walks on free products using, as in our work, Theorem \ref{thm.LDP.criterion} as a starting point.




\subsection{Basics on hyperbolicity}
\label{subsectionyperbolicity}
As general references on the topic one can recommend \cite{livregromovhyperbolic}, \cite{surveykapobenakli} and \cite{artvasailagromovhyp} for the non-proper setting. 

 \begin{definition}
	A metric space $(X,d)$ is said to be \textbf{Gromov-hyperbolic} if there is a constant $\delta > 0$ such that for any four points $\{x_i\}_{ 0 \le i \le 3}$ we have
	$$ (x_1, x_2)_{x_0} \ge \min \{ (x_3,x_1)_{x_0}, (x_3,x_2)_{x_0} \} - \delta, $$ where for $x,y,z \in X$, $(y, z)_{x}$ denotes the Gromov product  as defined in \eqref{eq.defn.gromov.product}. 
 \end{definition} 
 
In this article, we will mostly deal with geodesic spaces. Recall that a metric space $(X,d)$ is geodesic if the distance between any two points $x,y$ is given by the length of a rectifiable path whose endpoints are $x$ and $y$.\\

The following definition is to explain the terminology involved in the statement of Proposition \ref{theoremharmonicmeasure}. Let $X$ be a Gromov-hyperbolic metric space and $x_0 \in X$ a base point.

\begin{definition}
The \textbf{Gromov boundary}, denoted by $\partial X$, is defined to be the set of all sequences $(x_n)_{n \in \NN} \in X^{\NN}$ such that $\lim_{n,m \to \infty} \ (x_n, x_m)_{x_0} = \infty$ 
modulo the equivalence relation $(x_n) \sim (y_n)$ if $(x_n,y_n)_{z_0} \tends{n \to \infty} \infty $. We denote by $[(x_n)]$ the class of such a sequence.
\end{definition}

One can easily verify that the construction of $\partial X$ does not depend on the base point $x_0$. \\

Choose $\zeta := [(x^{\zeta}_n)] \in \partial X$ and $r > 0$ and set 
$$ B(\zeta, r) := \{ \zeta_2 := [(y_n)] \in \partial X \ , \ \limiinf{n \to \infty} \ e^{-(x^{\zeta}_n, y_n)_{x_0}} \le r \} \ . $$
We  define a topology on $\partial X$ by choosing the above sets as a neighborhood basis at $\zeta$. 	
The resulting topological space $\partial X$ is metrizable. The sets $B(\zeta, r)$  are `almost' balls of radius $r$. We refer to \cite[Section 5]{artvasailagromovhyp} for more details. 

\begin{definition}
  Let $\lambda, C > 0$ and $I$ a sub interval of $\NN$. A $(\lambda,C)$\textbf{-quasi-geodesic indexed by $I$} (simply called quasi-geodesic when not ambiguous) is a sequence $(x_n)_{n \in I}$ such that for any $n, m \in I$ 
 $$ \lambda^{-1} |n-m| -C	\le d(x_n,x_m) \le \lambda |n-m| +  C \ . $$
 In other words, a quasi-geodesic is a quasi-isometric embedding of $I$ into $X$. 
\end{definition}

One can easily verify that quasi-geodesics indexed by $\NN$ define a unique point in $\partial X$. Recall the statement of the fundamental Morse lemma. 

\begin{lemma}[Morse lemma] For any $\lambda, C > 0$ there is a constant $L = L(\lambda, C, \delta)$ such that any $(\lambda,C)$-quasi-geodesic having the same endpoints are $L$-close to one another.
\end{lemma}

The following definitions are to explain the terminology `non-elementary'.

\begin{definition}
An isometry $\gamma$ of a Gromov-hyperbolic space $X$ is called \textbf{loxodromic} if for a point $x \in X$ (equivalently any) the sequence	$ (\gamma^n \cdot x)_{n \in \ZZ} $ is a quasi-geodesic.
\end{definition}

In particular, a loxodromic element defines two points in the Gromov boundary $\gamma_+$ and $\gamma_-$ corresponding to the classes of the two quasi-geodesics defined by the future and the past. We say that two loxodromic elements $\gamma_1, \gamma_2$ are \textbf{independent} if the four points $\gamma_1^{\pm}, \gamma_2^{\pm}$ are distinct.

\begin{definition}\label{def.nonel}
A semigroup acting on $X$ by isometries is called \textbf{non-elementary} if it contains two independent loxodromic elements. 
\end{definition} 
We note that some authors use the term ``general type'' for subgroups containing two independent loxodromic elements (see \cite[\S 3]{CCMT} for a detailed discussion).  
\begin{remark}
1. For a $\CAT(0)$ Gromov-hyperbolic space $X$, the condition in  Definition \ref{def.nonel} is equivalent to requiring that the group generated by the semigroup contains two independent loxodromic elements (see e.g.~\cite[\S 6.2]{Das-Simmons-Urbanski}).\\[2pt]
2. For a general Gromov-hyperbolic space $X$, it follows from \cite[Theorem 6.2.3 and Proposition 6.2.14]{Das-Simmons-Urbanski} that a semigroup $S$ of isometries of $X$ contains two independent loxodromic elements if and only if $S$ is unbounded and the group generated by $S$ has two independent loxodromic elements.\footnote{We thank an anonymous referee for this remark.}
\end{remark}

Finally, we say that a probability measure $\mu$ on a group $\Gamma$ acting by isometries on a Gromov-hyperbolic space $X$ is \textbf{non-elementary} when its support generates a non-elementary subsemigroup.\\

Non-elementary groups have a lot of elements spreading apart points of $X$. The proof of the following lemma is a variation around the proof of the well known ping-pong lemma. As we could not find any ready-to-use reference in this generality, we inserted  a proof in Appendix \ref{appendixB}.

\begin{proposition}[Existence of Schottky sets]
\label{lemmapingpong} Let $\Gamma$ be  group acting by isometries on a geodesic Gromov-hyperbolic space $X$, $z_0 \in X$ and $\mu$ a non-elementary probability measure on $\Gamma$. Then there is $p \in \NN$ such that $\supp{\mu^{*p}}$ contains a Schottky set.
\end{proposition}

We note in passing that a slight modification of Abels--Margulis--Soifer's proof of \cite[Theorem 4.1]{AMS} yields a Schottky set in the sense of Definition \ref{defschottkyset} for symmetric spaces of non-compact type. See also \cite{DGLM} for related considerations.


\section{Deviations from above} 
\label{secupper}

The goal of this section is to prove that \eqref{eq.from.above} holds in the setting of Theorem \ref{thm.ldp.general}. 

\begin{proposition}\label{propupper}
Let $\Gamma$ be a countable group acting by isometries on a metric space $X$, $\mu$ a probability measure on $\Gamma$ and $z_0 \in X$. Suppose that the semigroup generated by the support of $\mu$ contains a Schottky set. Then, there is a non-negative convex function $\Psi : \ [ 0 , \infty ) \to [0,\infty]$ such that for any $a \neq l_{\max}$
\begin{equation*}
\frac{- \ln \PP \big( d(z_n,z_0) \ge a n \big)  }{n}  \tends{n \to \infty} \Psi(a) \ .
\end{equation*}
Moreveor, if $\mu$ has a finite exponential moment, then $\Psi$ vanishes only on $[0,l]$.
\end{proposition}
For the definition of the constant $l_{\max}$ in the above, see \eqref{eq.lmax}. 
The part of the previous result concerning $\Psi > 0$ on $(l,\infty)$ follows from Hamana's argument taken from \cite{arthamadaupperdev}. Namely, we will show in Appendix \ref{appendixA} that 

\begin{proposition}\label{propositionfromhamana}
Let $X$ be a metric space and $\mu$ a probability measure on $\Isom(X)$ with a finite exponential moment. Then for any $a > l$ we have 
\begin{equation}
	\label{eqpsotidevabove}
	 \limiinf{n \to \infty}  \ \frac{ - \ln \PP \big( d(z_n,z_0) \ge a n \big)  }{n} > 0 \ ,
\end{equation}
\end{proposition}

The proof of Proposition \ref{propositionfromhamana} only requires sub-additivity, which, for random walks, comes from the triangle inequality and the independence of the increments as shown in Section \ref{secrappel}. The rest of this section is devoted to answer the second part of the question: show that the limit defining $\Psi$ exists and that $\Psi$ is convex. \\

The next proposition gives an almost sub-additivity relation.

\begin{proposition}
	\label{propupperexist}
There is a constant $c > 1$ and an integer $p \in \NN$ such that for any $x,y$ in $\mathbb{R}_+$ and $n,m \in \mathbb{N},$ we have
	\begin{equation}	
		\label{equpperlimitexist2}
	 		 \PP( d_{m+n+p} \ge x + y - c) \ge c^{-1} \cdot \PP(d_m \ge x) \ \PP(d_n \ge y) \ . 
	\end{equation}
\end{proposition}

Before proving the above proposition, let us see how to use it to show that the limit defining $\Psi$ exists and is convex. \\

\textbf{Proof of (Proposition \ref{propupperexist} $\Rightarrow$ $\Psi$ exists and is convex).} Throughout the proof, $p$ is fixed as in Proposition \ref{propupperexist}. \\

To apply Fekete's lemma, we substitute in \eqref{equpperlimitexist2}, $m-p$ for $m$ and $n-p$ for $n$ to get that for any $x,y > 0$ and $m,n \ge p$:
\begin{equation} \label{equpperlimit3} 
\PP( d_{m  +n - p} \ge x + y - c) \ge c^{-1} \cdot \PP(d_{m-p} \ge x) \ \PP(d_{n -p} \ge y) \ . \end{equation} 

We now replace $x$ with $am + c$ and $y$ with $an + c$ in order to get that for all $m,n > p$
\begin{equation*}
	 \PP( d_{m+n -p} \ge a(m+n) + c) \ge  c^{-1} \cdot \PP(d_{m-p} \ge am +c)  \PP(d_{n-p} \ge an + c).
\end{equation*}
	
Thus we see that the sequence 	$(-\ln \left( c^{-1}\PP( d_{n-p} \ge an + c)\right))_{ n\geq p}$ is sub-additive. Let us define 
	$$ \psi_n(a):= \frac{ -\ln \left( c^{-1}\PP( d_{n-p} \ge an + c)\right) }{ n}  \ .$$
Fekete's lemma implies that, for all $a$, $(\psi_n(a))_{n \geq p}$ converges; we denote with $\Psi(a)$ the limit. \\

We now show that $\Psi$ is convex. Indeed, using Inequality \eqref{equpperlimit3} one gets that, for any $a,b > 0$ and for any 
$n \geq p$, we have
	$$ \psi_{2n} \left( \frac{a+b}{2} \right) \le  \frac{1}{2} \left( \psi_n ( a ) + \psi_n (b) \right) \ ,  $$

which shows, letting $n \to \infty$, that $\Psi$ is convex. \\ 

We now show that as $n \to \infty$ the sequence $-\frac 1 n \ln \PP(d_n\geq an)$ converges to $\Psi(a)$ for $a \neq l_{\max}$. We start with the observation that for any $\varepsilon > 0$ we have for $n$ large enough 
 $$ \PP( d_{n-p} \ge (a - \varepsilon) n + c) \ge \PP( d_{n-p} \ge a(n-p)) \ge \PP( d_{n-p} \ge a n + c). $$
Therefore 
\begin{equation}
\label{eq-inegalite dev par en haut}
	\begin{split}
			 \Psi(a-\varepsilon)\leq  \liminf -\frac 1 n \ln\PP(d_n\geq an) 
		\leq \limsup & -\frac 1 n \ln\PP(d_n\geq an)\leq \Psi(a)\ . 
	\end{split}
\end{equation}

The above inequality implies that if $a > l_{\max}$ then $\Psi(a) = \infty$. In particular, using again the above inequality, if $a > l_{\max}$ we get
$$ \liminf - \frac 1 n \ln\PP(d_n\geq an) \ge \Psi \left(\frac{l_{\max} + a }{2} \right) = \infty = \Psi(a)  \ .$$ 

We conclude showing that $- \frac 1 n \ln\PP(d_n\geq an)$ converges to $\Psi(a)$ for $a \in (0, l_{\max})$. Since $\Psi$ is convex and finite on $(0, l_{\max})$ it is in particular continuous. Letting $\varepsilon \to 0$ in \eqref{eq-inegalite dev par en haut} we get that the sequence $-\frac 1 n \ln \PP(d_n\geq an) $ converges to $\Psi(a)$ on $(0, l_{\max})$.  \hfill $\blacksquare$ \\

In the presence of a finite first moment (in particular, finite exponential moment), the almost sure convergence of $n^{-1} d_n$ to $l$ shows directly that $\Psi = 0$ on $[0,l]$ (if $l \neq l_{\max}$) using that $\Psi$ is convex (in particular continuous). Therefore in view of Proposition \ref{propositionfromhamana}, one is left to show that Proposition \ref{propupperexist} holds. Our strategy is inspired by the replacement trick proposed in \cite{arthamanakesten2} and by the use of a  Schottky set, inspired from \cite{artdalbofeteke}. \\

\textbf{Proof of Proposition \ref{propupperexist}.}
To ease the notation, we denote by $S_\mu$ the support of $\mu$. For an element $g$ in the semigroup generated by $S_\mu$, we write $|g|$ to denote the least number of factors needed to write $g$ as a product of elements of $S_\mu$. By hypothesis, there exists a Schottky set $S$ in the semigroup generated by $S_\mu$. Let $C>0$ be the associated constant as in Definition \ref{defschottkyset}. Let $p \in \mathbb{N}$ be such that any element of $S$ can be written as a product of at most $p$ elements of $S_\mu$. For $i=0,\ldots,p$, we fix some elements $h_i \in S_{\mu^{\ast i}}$ with $h_0=\mathrm{id}$. We let
\begin{equation*}
\begin{aligned}    
c_1 &=2C+\max\{d(h_i \cdot z_0,z_0,)  \; | \; s \in S \; , \; i=0,\ldots, p\} \qquad \text{and}
\\ c_2 &=\min\{\mu^{\ast i}(h_i) \cdot \mu^{\ast |s|}(s)  \; | \; s \in S \; , \; i=0,\ldots, p\}.
\end{aligned}
\end{equation*}

Let $m,n \in \mathbb{N}$ and $x,y \in \mathbb{R}_+$ be as in the statement. Using the defining property of the Schottky set $S$, for every $g_m \in S_{\mu^{\ast m}}$ and $g_n \in S_{\mu^{\ast n}}$, we fix an element 
$s=s(g_m,g_n)$ such that $(g_m^{-1} \cdot z_0, s g_n \cdot z_0)_{z_0} \leq C$, equivalently,
\begin{equation}\label{eq.fix.ss}
d(g_msg_n \cdot z_0, z_0) \geq d(g_m \cdot z_0, z_0) + d(g_n \cdot z_0, z_0) -2C.
\end{equation}

Abbreviating $d(g \cdot z_0,z_0) \geq x$ by $g \geq x$, we have 
\begin{equation*}
\begin{aligned}
&\mathbb{P}(d_{m+n+p} \geq x+y-c_1)\\&= \sum_{\substack{g_m \in S_{\mu^{\ast m}}\\g_n \in S_{\mu^{\ast n}}}} \sum_{k=0}^p \sum_{\substack{g_k \in S_{\mu^{\ast k}}\\g_{p-k} \in S_{\mu^{\ast (p-k)}}}} 1_{g_m g_k g_n g_{p-k} \geq x+y-c_1}  \mu^{\ast m}(g_m) \mu^{\ast k}(g_k) \mu^{\ast n}(g_n) \mu^{\ast (p-k)}(g_{p-k}) \\
& \geq c_2 \cdot \sum_{\substack{g_m \in S_{\mu^{\ast m}}\\g_n \in S_{\mu^{\ast n}}}}  1_{g_m s(g_m,g_n) g_n h_{p-|s(g_m,g_n)|} \geq x+y-c_1} \cdot \mu^{\ast m}(g_m) \mu^{\ast n}(g_n)\\
&\geq c_2 \cdot \sum_{\substack{g_m \in S_{\mu^{\ast m}}\\g_n \in S_{\mu^{\ast n}}}} 1_{g_m \geq x} 1_{g_n \geq y} \mu^{\ast m}(g_m) \mu^{\ast n}(g_n)= c_2\mathbb{P}(d_m \geq x) \mathbb{P}(d_n \geq y).
\end{aligned}
\end{equation*}
In the above, to pass from first line to the second, we used the I.I.D. property. To pass from second to the third, for each $g_m$ and $g_n$, we specialized to the $k \leq p$ such that $g_k=s(g_m,g_n)$ and to $g_{p-k}=h_{p-k}$, and used the definition of $c_2$. To pass to the last line we used \eqref{eq.fix.ss} and the definition of $c_1$. Therefore the proposition follows by setting $c= \max\{\frac{1}{c_2},c_1\}$. \hfill $\blacksquare$\\

\section{Deviations from below} 
\label{seclower} 

This section is dedicated to investigating the deviations from below. The strategy of the proof of the following proposition is inspired from \cite{artshapiraasselah}.

\begin{proposition}\label{proplowerdev}
Let $\Gamma$ be a countable group acting on a metric space $X$ and $\mu$ a probability measure on $\Gamma$. 
Then there is a convex function $\Psi : \ [ 0 , \infty ) \to [0,\infty]$ such that for all $a \neq l_{\min}$
\begin{equation}\label{eqproplower}
\frac{-\ln \PP \big( d(z_n,z_0) \le a n \big)  }{n}  \tends{n \to \infty} \Psi(a) \ .
\end{equation}
Furthermore, if $\mu$ has a finite exponential moment and 
satisfies 
\begin{equation}\label{eqassumptionproplower}
\limiinf{p \to \infty} \ \supr{x \in X} \ \frac{\EE  \big( ( x, z_p)_{z_0} \big)}{p} = 0,
\end{equation}
then $\Psi$ vanishes only on $[l, \infty]$.
\end{proposition}
For the definition of the constant $l_{\min}$ in the above, see \eqref{eq.lmin}. To prove the previous result, we shall start by showing that the limit defining the function $\Psi$ exists. This only requires sub-additivity. We will then prove the most difficult part of the proof, namely that $\Psi > 0$ under the assumption \eqref{eqassumptionproplower}. \\

\textbf{Proof that the limit exists.} The proof does not require Assumption \eqref{eqassumptionproplower}. By the triangle inequality and independence, we have
	$$ \PP(d_{n+m} \le a (n + m) ) \ge  \PP(d_{m} \le a m ) \  \PP(d_{n} \le a n ) \ .$$

Therefore the sequence $(- \ln  \PP(d_{n} \le a n ))_{n\in\NN}$ is sub-additive. Let us define 
	$$( \Psi_n(a) )_{n \in \NN} :=  \left( \frac{- \ln \Big( \PP(d_{n} \le a n ) \Big)}{n} \right)_{n \in \NN} \ .$$
Fekete's lemma then gives that the sequence $( \Psi_n(a) )_{n \in \NN}$ converges; we denote the limit with $\Psi(a)$. \\

To show that $\Psi$ is convex, let $a,b \in \mathbb{R}$.  
Using again the triangle inequality, we get
	$$ \PP \left(d_{2n} \le \frac{a+b}{2} \cdot 2n \right) \ge \PP(d_{n} \le an) \ \PP(d_{n} \le bn) \ , $$
and then 
	$$ \Psi_{2n} \left( \frac{a+b}{2} \right) \le \frac{\Psi_n(a) + \Psi_n(b)}{2} \ . $$ 

We conclude by letting $n$ tend to $\infty$. \hfill $\blacksquare$ \\
	 
\textbf{Proof that $\Psi > 0$.} We will now use Assumption \eqref{eqassumptionproplower} and the finite exponential moment hypothesis. \\

Let us start by noticing that Proposition \ref{proplowerdev} is invariant under acceleration: given $k \in \NN$ a measure $\mu$ with a finite exponential moment satisfies the conclusion of Proposition \ref{proplowerdev} if and only if the measure $\mu^{*k}$ satisfies it. \\

Given a trajectory, we chop it into pieces of size $j \in \NN$ and write the distance between the base point $z_0$ and the endpoint $z_n$ (where $n = mj$ for some integer $m$) as a summation of I.I.D.\ random variables and a defect term. \\


By definition of the Gromov product, we have for any $m,j > 0$ 
	$$ 2(z_0, z_{mj})_{z_{(m-1)j}} = d(z_0, z_{(m-1)j}) + d(z_{mj}, z_{(m-1)j}) - d(z_0, z_{mj}) \ . $$

Equivalently,
	$$ d_{mj} = d_{(m-1)j} + d(z_{mj},z_{(m-1)j}) - 2(z_0, z_{mj})_{z_{(m-1)j}} \ . $$

By an immediate induction we get 
	$$ d_{mj} = \somme{1 \le i \le m}  d(z_{ij},z_{(i-1)j}) - 2 \somme{1 \le i \le m} (z_0, z_{ij})_{z_{(i-1)j}} \ . $$	
	
Since the Gromov product is non-negative, one has the following set inclusion 
$$ \left\{ d_{mj} \le a n \right\} \subset  \left\{ \somme{1 \le i \le m}  d(z_{ij},z_{(i-1)j})  \le \frac{a+l}{2} n \right\} \bigcup \left\{ \somme{1 \le i \le m} (z_0, z_{ij})_{z_{(i-1)j}} \ge \frac{l-a}{4} n  \right\}  \ , $$

which implies that 
\begin{equation}	
	\label{eqlargedev1}
	\begin{split}
	 \PP( d_{mj} \le a n ) \le  \PP & \left(  \somme{1 \le i \le m}  d(z_{ij},z_{(i-1)j})  \le \frac{a+l}{2} n \right)  + 
	 \PP \left( \somme{1 \le i \le m} (z_0, z_{ij})_{z_{(i-1)j}} \ge \frac{l-a}{4} n \right)  \ . 
	 \end{split} 
\end{equation} 

We shall see that there exists $j$  such that both the above probabilities decay exponentially fast to $0$. The argument for the first one only uses classical large deviations estimates for I.I.D.\ random variables whereas the control of the second one will be handled using Assumption \eqref{eqassumptionproplower}. \\

We start with the top probability appearing in \eqref{eqlargedev1}. The random variables $ (d(z_{ij},z_{(i-1)j}))_{i \in \NN}$ are I.I.D.\ and follow the law of $d_j$. Therefore large deviations estimates for I.I.D.\ random variables with a finite exponential moment imply that 
$ \PP \left( \somme{1 \le i \le m}  d(z_{ij},z_{(i-1)j})  \le \frac{a+l}{2} n \right)$ has an exponential decay as soon as 
$\frac{\EE(d_j)}{j}> \frac{a+l}2$. \\

On the other hand, we already know  that $\frac{\EE(d_j)}j$ converges to $l$ and $l>\frac{a+l}{2}$. 
Thus we conclude that there exists $j_0$ such that for all $j \geq j_0$, we have 

$$\liminf -\frac 1 n \ln \PP \left( \somme{1 \le i \le m}  d(z_{ij},z_{(i-1)j})  \le \frac{a+l}{2} n \right)>0\ .$$


We now deal with the second probability appearing in \eqref{eqlargedev1} using Assumption \eqref{eqassumptionproplower}. 

Let us set $\varepsilon := \frac{l-a}{4}$ and let $\lambda > 0$. 
We start with the Chernoff bound 
\begin{equation}
	\label{equationmarkov}
		\begin{split}
			  \PP \left(  \somme{1 \le i \le m} (z_0, z_{ij})_{z_{(i-1)j}} \ge \frac{l-a}{4} n \right)  \le  
			  		e^{ - \lambda \varepsilon n}  \cdot  \EE  & \left( \exp \left( \lambda \somme{1 \le i \le m} (z_0 , z_{ij})_{z_{(i-1)j}}  \right) \right)  \ .
		\end{split}
\end{equation}
We introduce the random variables
	$$ \Pi_m(\lambda,j) := \exp \left( \lambda \somme{1 \le i \le m} (z_0, z_{ij})_{z_{(i-1)j}}  \right) \ , $$

and note that 
 $$ \Pi_m(\lambda,j)  = \Pi_{m-1}(\lambda,j) \cdot  \exp \left( \lambda (z_0, z_{mj})_{z_{(m-1)j}} \right) \ . $$
 
Let us denote with $(\mathcal{F}_i)_{i \in \NN}$ the filtration naturally associated to the random walk. We compute
 \begin{align*}
 	  \EE  \left( \Pi_m(\lambda,j) \right) & = \EE  \Big( \EE  \big( \Pi_{m-1}(\lambda,j) \cdot  \exp \left( \lambda (z_0, z_{mj})_{z_{(m-1)j}} \right) \big| \mathcal{F}_{(m-1)j} \big) \Big) \\
 	 & =	\EE  \Big( \EE  \big( \Pi_{m-1}(\lambda,j) \cdot \exp \left( \lambda (\gamma_{(m-1)j}^{-1} z_0, \gamma_{(m-1)j}^{-1} z_{mj})_{z_0} \right) \big| \mathcal{F}_{(m-1)j} \big) \Big) \\
 	 & = \EE  \Big(  \Pi_{m-1}(\lambda,j) \cdot \EE \big( \exp \left( \lambda (\gamma_{(m-1)j}^{-1} z_0, \gamma_{(m-1)j}^{-1} z_{mj})_{z_0} \right) \big| \mathcal{F}_{(m-1)j} \big) \Big) \ .
 \end{align*}
The last equality holds because $ \Pi_{m-1}(\lambda,j) $ is measurable with respect to $\mathcal{F}_{(m-1)j}$. 
Moreover, since $ \gamma_{(m-1)j}^{-1} z_{mj}$ is independent of $\mathcal{F}_{(m-1)j}$ and since $ \gamma_{(m-1)j}^{-1} z_{mj}$ follows the same law as $z_j$, 
we have 
 \begin{align*}
 	  \EE  \left( \Pi_m(\lambda,j) \right) & \le \EE  \Big( \Pi_{m-1}(\lambda,j) \Big) \cdot \supr{x \in X} \ \EE  \Big( \exp \left( \lambda (x, \gamma_{(m-1)j}^{-1} z_{mj})_{z_0} \right) \Big) \\
 	 & \le	 \EE  \Big( \Pi_{m-1}(\lambda,j) \Big)  \cdot \supr{x \in X} \ \EE  \Big( \exp \left( \lambda (x, z_j)_{z_0} \right) \Big) \ .
 \end{align*} 
An immediate induction yields
  $$	 \EE  \left( \Pi_m(\lambda,j) \right) \le \delta(j, \lambda)^m \ , $$
  
where 
  	$$ \delta(j,\lambda) := \supr{x \in X} \ \EE  \Big( \exp \left( \lambda (x, z_j)_{z_0} \right) \Big) \ . $$

Therefore, 
\begin{align*}	
			  \PP \left( 2 \somme{1 \le i \le m} (z_0, z_{ij})_{z_{(i-1)j}} \ge \frac{l-a}{2} n \right)  \le
			  		e^{ - \lambda \varepsilon n} \delta(j, \lambda)^m 
			  		& \le e^{ - \lambda \varepsilon n + m \ln \left( \delta(j, \lambda) \right) } \\
			  		& \le e^{ m [  \ln \left( \delta(j, \lambda)\right) -\lambda j \varepsilon ]} \ . 
\end{align*}	
We shall prove, using Assumption \eqref{eqassumptionproplower}, that for all $\varepsilon' > 0$ there exist $j \geq j_0 \in \NN$ and $\lambda > 0$ such that 
\begin{equation}
	\label{eqpreuvepropdevtolargedev}
	  \frac{ \ln(\delta(j, \lambda))}{ \lambda j} \le \varepsilon' \ . \end{equation}
This is enough to conclude:  we choose $\varepsilon':= \varepsilon / 2 $ with $j$ and $\lambda$ such that  \eqref{eqpreuvepropdevtolargedev} holds. Then  
$$ \PP \left( 2 \somme{1 \le i \le m} (z_0, z_{ij})_{z_{(i-1)j}} \ge \frac{l-a}{2} n \right) \le
			  		e^{ - \lambda \varepsilon n} \delta(j, \lambda)^m  \le e^{- n  \frac{ \lambda \varepsilon}{ 2}} \ ,  $$
does indeed  decrease exponentially fast to $0$ as $n \to \infty$. \\ 

It remains to prove Inequality \eqref{eqpreuvepropdevtolargedev}. Note first that for any $x \in X$ we have 
 	$$ \EE  \Big( \exp \left( \lambda (x, z_j)_{z_0} \right) \Big) \le 1 + \lambda \ \EE  \left( (x, z_j)_{z_0} \right) + \lambda^2 \   \EE  \left( \left( (x, z_j)_{z_0} \right)^2  \exp \left( \lambda (x, z_j)_{z_0} \right) \right) \ ,  $$
since  $e^x \le 1 + x + x^2 e^x$. Using the upper bound $ (x, z_j)_{z_0} \le d(z_0, z_j)$, we get that 
 	$$ \EE  \Big( \exp \left( \lambda (x, z_j)_{z_0} \right) \Big)   \le 1 + \lambda \ \EE  \left( (x, z_j)_{z_0} \right) + \lambda^2 \  \EE_{z_0 } \left(  d_j^2 \  e^{\lambda d_j} \right) \ .  $$
	
Assumption \eqref{eqassumptionproplower} provides us with some $j_1$ such that, for all $j \geq j_1$,  we have  
	$$ \supr{x \in X} \ \EE  \left( (x, z_j)_{z_0} \right) \le \frac{\varepsilon' j }{2} \  . $$ 
We choose $j \geq \max(j_0, j_1)$.  Then, taking the sup over $x \in X$, we get 
 	$$ \delta(j, \lambda)  \le 1 + \frac{\lambda \varepsilon' j}{2}  + \lambda^2 \ \EE  \left( d_j ^2 \ e^{ \lambda d_j} \right) \ . $$
 	
We now choose $\lambda = \lambda(j)$ small enough such that $\lambda^2 \ \EE  \left( d_j ^2 \ e^{ \lambda d_j} \right) \le \frac{\lambda \varepsilon' j}{2}$. Then $\delta(j, \lambda) \le  1 + \lambda \varepsilon' j$, 
and therefore,  since $\ln(1 +x ) \le x$, we have $ \ln (\delta(j, \lambda)) \le \lambda \varepsilon' j$.
\hfill $\blacksquare$


\section{Walking-away uniformly}
\label{seclinearprogress}

\begin{definition}
	A sequence of random variables $(Z_n)_{n \in \NN}$ taking values in a metric space $X$ is said to satisfy the \textbf{walking-away uniformly property} if there are constants $\varepsilon, \alpha, C > 0$ such that for all $x \in X$ and for all $n \in \NN$
	$$ \PP \left(  d(Z_n,x) - d(z_0,x) \le \varepsilon n \right) \le C e^{-\alpha n} \ .$$ 
\end{definition}

Note that the above definition does not actually depend on the random variables $(Z_n)_{n \in \NN}$ but only on their laws. We shall use this fact in the proof of the following theorem by exhibiting a special set of random variables which have the desired law. 

\begin{theoreme}
	\label{theowalkingaway}
Let $\Gamma$ be a countable group acting by isometries on a metric space $X$  and $\mu$ a probability measure on $\Isom(X)$ with finite exponential moment whose support generates a subsemigroup which contains a Schottky set and which has unbounded orbits. Then, $(z_n)_{n \in \NN}$ satisfies the walking-away uniformly property. 
\end{theoreme}

Notice that Lemma \ref{lemma.Psi.and.I}, Propositions \ref{propupper} and \ref{proplowerdev}, and the previous theorem completes the proof of Theorem \ref{thm.ldp.general}.


\subsection{Overview of the argument}
\label{subsec.overviewargument} The proof of the above theorem is quite intricate. Let us start by noticing that Theorem \ref{theowalkingaway} is invariant under acceleration: given $k \in \NN$ a measure $\mu$ with a finite exponential moment satisfies the conclusion of Theorem \ref{theowalkingaway} if and only if the measure $\mu^{*k}$ satisfies it. Moreover, we can assume without loss of generality that the identity element is in the support of $\mu$. For a probability measure $\mu$ with $\mu(\mathrm{id})>0$, it is clear that if the semigroup generated by the support of $\mu$ contains a Schottky set, then the support of some convolution power of $\mu$ contains a Schottky set. Therefore, to start with, we can assume that the support of $\mu$ contains a Schottky set.\\

We start by showing that the above theorem is also invariant under sampling. More precisely, we will sample the positions $(z_n)_{n \in \NN}$ along the times when drawing increments in a given set $S$. We shall then use this sampling with respect to a Schottky set.  \\

To make it precise, we will first exhibit a special family of increments $(\omega_i)_{i \in \NN}$ (following the law $\mu$) using the following random variables. Let $S \subset \supp{\mu}$ be any finite set and
	$$\zeta := \min_{\gamma \in S} \ \mu(\gamma) > 0 \ .$$ 

Let $(\eta_i)_{i \in \NN}$ be independent random variables following the Bernoulli law of parameter $\zeta$. Let also $(V_i)_{i \in \NN}$ I.I.D.\ random variables independent of the $\eta_i$'s taking values in $\Gamma$ with (common) distribution 
	$$ \PP(V_i = \gamma) := 
		\left\{ 
			\begin{array}{lr}
				 (1 - \zeta)^{-1} \big( \mu(\gamma) - \frac{\zeta}{\# S} \big) & \text{ if } \gamma \in S		\\		
				 (1 - \zeta)^{-1} \mu(\gamma) & \text{ if } \gamma \notin S \\
			\end{array}
		\right. $$
This distribution defines a probability measure on $\Gamma$ since $\mu(\gamma) - \frac{\zeta}{\# S} \ge 0$ by definition of $\zeta$ and since its total mass is $1$ by construction. Note also that the random variables $(d(V_i \cdot z_0, z_0))_{i \in \NN}$ have a finite exponential moment since the measure $\mu$ has it (the laws of the $d(V_i \cdot z_0, z_0)$ are proportional to $\mu$ on all but finitely many $\gamma \in \Gamma$).\\

Let us now introduce the last set of random variables that we will need. Let $(S_i)_{i \in \NN}$ be I.I.D.\ random variables uniformly distributed on $S$ independent of all the $\eta_i$'s and of the $V_i$'s:
	$$ \PP(S_i = \gamma) :=
		\left\{ 
			\begin{array}{lr}
				 (\# S)^{-1} & \text{ if } \gamma \in S		\\		
				 0 & \text{ if } \gamma \notin S \ . \\
			\end{array}
		\right. $$

In total, we are left with three sets of random variables that are all independent from one another. Finally, note that the following defined random variables (also taking values in $\Gamma$)
	$$ \omega_i :=
		\left\{ 
			\begin{array}{lr} 
				S_i  & \text{ if } \eta_i = 1 \\
				V_i  & \text{ if } \eta_i = 0 \\
			\end{array}
		\right. $$
follow the law of $\mu$. Indeed, by construction of the $V_i$, the $S_i$ and the $\eta_i$, one has 
	\begin{align*}
		 \PP(\omega_i = \gamma)  & = 
		\left\{ 
			\begin{array}{lr}
				\PP(V_i = \gamma) \PP(\eta_i = 0) + \PP(\eta_i = 1)\PP(S_i = \gamma) & \text{ if } \gamma \in S \\
				\PP(V_i = \gamma) \PP(\eta_i = 0) 	& \text{ if } \gamma \notin S \\
			\end{array}
		\right. \\
		& = \mu(\gamma) \ .
	\end{align*}

We endow our new probability space with the filtration $(\mathcal{F}_i)_{i \in \NN}$ corresponding to events which can be expressed using the random variables defined above only with indices $\le i$. \\
		
Let $p \in \NN$ fixed. We now define the $(S,p)$-sampling that we will use through the following sequence of stopping times, defined inductively as $\tau(0)=0$ and
	$$ \left\{ 
			\begin{array}{l}
				\tau(1) := \inf \ \{ k \ge p \ , \ \eta_k = 1 \} \\
				\tau(i) := \inf \ \{ k \ge p + \tau(i-1)  \ , \ \eta_k = 1 \}	\ \text{ if } i > 1  \ .
			\end{array}  
		\right. $$

The reason why we introduce an extra parameter $p$ will become clear later. Intuitively, we will use this parameter in order to guarantee that the average distance the random walk travels between positions at times $\tau(i)$ and $\tau(i +1)$ is large compared to the constant $C$ appearing in Definition \ref{defschottkyset}. \\

Note that the random variables $(\tau(i+1) - \tau(i))_{i >0}$ are I.I.D.\ following the law of $\tau(1)$ since the $(\eta_k)_{k \in \NN}$ are I.I.D. \\

The sampling on $\Gamma$ is defined according to the previously defined stopping time. Namely, it is the random walk whose successive positions are 
	$$ \gamma_{\tau(n)} = \omega_1 \cdot ...\cdot  \omega_{\tau(n)} \ .$$

By construction, the random variable $\gamma_{\tau(n)}$ follows the law $\mu_{\tau}^{*n}$, where 
	$$\mu_{\tau}(\gamma) := \PP(\omega_1 \cdot ... \cdot \omega_{\tau(1)} = \gamma ) \ .$$

\begin{definition}
\label{defshcottkyrandomwalk}
The corresponding image random walk on $X$, whose positions are $ z_{\tau(n)} = \gamma_{\tau(n)}  \cdot z_0 \ ,$ is called the $(S,p)$\textbf{-sampling of} $(z_n)_{n \in \NN}$.
\end{definition}

The following proposition guarantees that one can prove Theorem \ref{theowalkingaway} for the sampled random walk instead of the original random walk. 

\begin{proposition}
	\label{propsampling}
Let $\mu$ be a probability measure with a finite exponential moment on a group $\Gamma$ which acts on a metric space $X$. Let $S \subset \supp{\mu}$, $p >0$ and $z_0 \in X$. The image random walk driven by $\mu$ satisfies the walking-away uniformly property if and only if its $(S,p)$-sampling satisfies it too.
\end{proposition}

In order to keep this subsection as an overview, we postpone the proof of the above Proposition to Subsection \ref{subsecpropsampling}. The proof makes use of the following simpler lemma whose proof is also postponed. 

\begin{lemma}
	\label{lemfiniteexpmomentsampling}
Let $\mu$ be a probability measure with a finite exponential moment on a group $\Gamma$ which acts on a metric space $X$, $z_0 \in X$ and $\tau(1)$ as above. Then, the random variables $d(z_{\tau(1)}, z_p)$, $d(z_{\tau(1)}, z_0)$ and $d(z_{\tau(1) -1}, z_p)$ have a finite exponential moment.
\end{lemma}	

We will then prove that the $(S,p)$-sampled random walk satisfies the walking-away uniformly property. In order to do so, we shall introduce a last type of random walks. Intuitively, a $(S,p)$-sampling can be thought as a process in two steps. First, we ignore the first $p$ increments and we do not draw 'bad elements' from $S$ (corresponding to $\eta_k = 0$) for a random time which follows a geometric law. Secondly, we draw an element uniformly from the set $S$. We shall make this precise by showing that a $(S,p)$-sampled random walk can be seen as a random walk whose odd increments correspond to the first step described above and the even ones to the second step, as in the following definition. \\

Let $\mu_1$ be a probability measure on $\Gamma$, $(X_i)_{i \in 2 \NN+1}$ I.I.D.\ random variables following the law $\mu_1$ and $(Y_i)_{i \in 2\NN}$ I.I.D.\ random variables uniformly distributed on the set $S$ and independent of the $X_i$'s. 
 
\begin{definition}
Let $(X_i)_{i \in 2 \NN+1}$ and $(Y_i)_{i \in 2\NN}$ as above. We call (the laws of) the following sequence of random variables a $(\mu_1, S)$\textbf{-random walk}  
	$$ z_n^{\mu_1,S} := 
		\left\{ 
			 \begin{array}{lr} 
			 	X_1 \cdot Y_2 \cdot X_3 \cdot ... \cdot X_{n} \cdot z_0 & \text{ if } n \ \text{ is odd } \\ 
			 	X_1 \cdot Y_2 \cdot X_3 \cdot ... \cdot Y_{n} \cdot z_0 & \text{ if } n \ \text{ is even } 
		 	\end{array}
		 \right.	
	$$
\end{definition} 

The following lemma relates the position at time $n$ of a $(S,p)$-sampled random walk to the position at time $2n$ of a $(\mu_1,S)$-random walk.	
	
\begin{lemma}
\label{lemsamplingtoschottky}
Let $\mu$ be a probability measure on a group $\Gamma$ which acts on a space $X$, $z_0 \in X$ and $(\tau(i))_{i \in \NN}$ as above. The sequence of random variables $(z_{\tau(n)})_{n \in \NN}$ follows the law of the sequence $(z^{\mu_1,S}_{2n})_{n \in \NN}$ with
		$$ \mu_1(\gamma) :=  \PP(\omega_1 \cdot ... \cdot \omega_p \cdot ... \cdot \omega_{\tau(1)-1} = \gamma ) \ .$$
In particular a $(S,p)$-sampled random walk satisfies the walking-away uniformly property if and only if its associated $(\mu_1, S)$-random walk satisfies it too.
\end{lemma}

\textbf{Proof.} We first set the random variables $X_i$ and $Y_i$ as
$$ 
	\left\{ \begin{array}{l}
			Y_{2i} := \omega_{\tau(i)} \\
			X_{2i +1} := \omega_{\tau(i)+1} \cdot .... \cdot \omega_{\tau(i + 1)-1} \ .
		\end{array}	 
	\right.
$$	

Then, by definition,  $z_{\tau(n)}=X_1 \cdot Y_2 \cdot X_3 \cdot ... \cdot X_{2 n -1} \cdot Y_{2n} \cdot z_0$. 

It follows from the independence properties of the random variables $X_i$'s, $V_i$'s and $\eta_i$'s that the random variables 
$(X_{2i-1},Y_{2i})_{i\geq 1}$ are I.I.D. Using the fact that, on the set $\tau(1)=k$, 
we have $Y_2=S_k$ and $X_1=V_1\cdot ...\cdot V_{k-1}$, it is also easy to see that $X_1$ and $Y_2$ are independent. 	\hfill $\blacksquare$ \\

The next step is to find a criterion on $\mu_1$ which guarantees that if $S$ is a Schottky set then the associated $(\mu_1,S)$-random walk satisfies the walking-away uniformly property. The following result is the key and its proof will occupy Section \ref{secintertwine}.

\begin{proposition}
	\label{propsamplingschottky}
	For any Schottky set $S$ there is a constant $M > 0$ such that the following holds. For any probability measure $\mu_1$ with a finite exponential moment and 
	$$ \somme{\gamma \in \Gamma} \ \mu_1(\gamma) \ d(z_0, \gamma \cdot z_0)  > M, $$
the $(\mu_1,S)$-random walk satisfies the walking-away uniformly property.
\end{proposition}

Let us see how to deduce Theorem \ref{theowalkingaway} with all the material introduced above. Recall that we fix $\mu$ a probability measure on $\Gamma$ and $S$ a Schottky set contained in the support of $\mu$. Proposition \ref{propsampling} implies that it is sufficient to prove the walking-away uniformly property for the $(S,p)$-sampled random walk. Because of Lemma \ref{lemsamplingtoschottky}, we know that the $(S,p)$-sampled random walk is also a $(\mu_1,S)$-random walk with 
	$$\mu_{1}(\gamma) := \PP(\omega_1 \cdot ... \cdot \omega_p \cdot ... \cdot \omega_{\tau(1)-1} = \gamma ) \ .$$
It remains to show that  the resulting $(\mu_1,S)$-random walk satisfies the conditions of Proposition \ref{propsamplingschottky}. 
Note that Lemma \ref{lemfiniteexpmomentsampling} already asserts that $\mu_1$ has a finite exponential moment. The following lemma ensures that we can choose $p$ such that the mean $\sum_{\gamma \in \Gamma} \ \mu_1(\gamma) \ d(z_0, \gamma \cdot z_0) $ exceeds $M$.

\begin{lemma}
\label{lemassumptionschottkywalk}
Let $\Gamma$ be a countable group acting by isometries on a metric space $X$ and $\mu$ be a probability measure on $\Isom(X)$ whose support generates a subsemigroup with unbounded orbits and assume that $\mu$ has a finite first moment. Then 
$$	\limisup{p \to \infty} \ \EE (d(z_0, z_{\tau(1)-1})) = \infty \ . $$
\end{lemma}

The following subsections are devoted to the proofs of all the above lemmata, except Proposition \ref{propsamplingschottky} which will be proven in Section \ref{secintertwine}.

\subsection{Proof of Lemma \ref{lemfiniteexpmomentsampling}}
\label{subseclemsampling} The differences of any two of the three random variables appearing in Lemma \ref{lemfiniteexpmomentsampling} obviously have  a finite exponential moment. It is therefore sufficient to prove Lemma \ref{lemfiniteexpmomentsampling} for one of them only, say $d(z_p, z_{\tau(1)})$. The proof is a straightforward computation. It only uses  that $\mu$ has a finite exponential moment together with the fact that $\tau(1)-p$ follows a geometric  law of parameter $1-\zeta$. Given $\lambda > 0$ we compute
\begin{align*}
	 \EE \left( e^{\lambda \ d(z_{\tau(1)}, z_p)} \right) & \le \EE \left( \exp \left( \lambda \ \somme{p \le i \le \tau(1)-1} d(z_i, z_{i+1}) \right) \right) \\
	 & \le \somme{ k \in \NN} \ \EE \left( \exp \left( \lambda \ \somme{p \le i \le k-1} d(z_i, z_{i+1}) \right) \Big| \tau(1) = k \right) \PP( \tau(1) = k) \\
 	 & \le \zeta \ \somme{ k \in \NN} \ \EE \left( \exp \left( \lambda \ \somme{p \le i \le k-1} d(z_i, z_{i+1}) \right) \Big| \tau(1) = k \right) (1-\zeta)^{k-p} \ . 
\end{align*}

We shall now see that for all $\varepsilon > 0$ there is $\lambda > 0$ such that for every $k \in \NN$ 
	\begin{equation}\label{uneequation}  E_k := \EE \left( \exp \left( \lambda \ \somme{p \le i \le k-1} d(z_i, z_{i+1}) \right) \Big| \tau(1) = k \right) < (1 + \varepsilon)^{k} \ . \end{equation} 
It concludes the proof since we can choose $\varepsilon$ such that $(1 + \varepsilon)(1 - \zeta)^{-1} < 1$. \\

Let us check (\ref{uneequation}). 
The event $\tau(1) = k$ is defined as $ \eta_p = 0, \ \eta_{p+1} = 0, ...., \eta_{k} = 1$. Therefore, by construction of the $X_i$'s, we have  
\begin{align*}
		E_k & = \EE \left( \exp \left( \lambda \ \somme{p \le i \le k -1} d(V_i \cdot z_0, z_0) \right) e^{\lambda d(z_0, S_k \cdot z_0)} \Big| \tau(1) = k \right) \ .
\end{align*}
But $\tau(1)$ is a function of the $\eta_i$'s only and therefore is independent of $S_k$ and independent of the $(V_i)_{1 \le i \le k -1}$. It yields 
\begin{align*}
		E_k & = \EE \left( \exp \left( \lambda \ \somme{p \le i \le k -1} d(V_i \cdot z_0, z_0) \right) e^{\lambda d(z_0, S_k \cdot z_0)} \right) \\ 
			& = \left( \EE \left( e^{\lambda d(V_1 \cdot z_0, z_0)} \right) \right)^{k-1 -p} \  \EE \left( e^{\lambda d(z_0, S_k \cdot z_0)} \right) \ , 
\end{align*}
since the $V_i$'s are I.I.D.\ and independent of $S_k$. This concludes the proof since we already saw that $d(V_i \cdot z_0, z_0)$ has a finite exponential moment and since $S_k$ has finite support.	 \hfill $\blacksquare$  	

\subsection{Proof of Proposition \ref{propsampling}}
\label{subsecpropsampling}
We will prove that:  if the random walk $(z_{\tau(n)})_{n \in \NN}$ satisfies the walking-away uniformly property then $(z_n)_{n \in \NN}$ satisfies it too. This is the only implication we need in this paper. The proof of the other implication is very similar. \\

Let $\varepsilon, C$ and $\alpha$ such that for any $x \in X$ we have
	$$ \PP \left(  d(z_{\tau(n)},x) - d(z_0,x) \le \varepsilon n \right) \le C e^{-\alpha n} \ .$$ 

We set $\beta := \EE(\tau)$. We will show that $(z_{\beta n})_{n \in \NN}$ satisfies the walking-away uniformly property, which implies the result using again the invariance under acceleration. \\

Rewriting $d(z_0, z_{\beta n})-d(z_0,x)$ as $d(z_0, z_{\beta n}) - d(z_0, z_{\tau(n)}) + d(z_0, z_{\tau(n)})-d(z_0,x)$, we have
	\begin{equation*}
		\begin{split}
	 \Big\{ d(z_{\beta n},x) & - d(z_0,x) \le  \frac{\varepsilon n}{2} \Big\}  \subset \\
	  & \left\{ d(z_{\tau(n)},x) - d(z_0,x) \le \varepsilon n  \right\}  \cup \left\{d(x, z_{\tau(n)}) - d(x, z_{\beta n}) \ge \frac{\varepsilon n}{2} \right\} \ .
	 	\end{split}
	\end{equation*}

And then:
	\begin{equation*}
		\begin{split}
			\PP \left( d(z_{\beta n},x) - d(z_0,x) \le  \frac{\varepsilon n}{2} \right)  \le & \\
	 	 \PP \left( d(z_{\tau(n)},x) - d(z_0,x) \le \varepsilon n  \right) &  + \PP \left( d(x, z_{\tau(n)}) - d(x, z_{\beta n}) \ge \frac{\varepsilon n}{2} \right) \ .
	 	\end{split}
	\end{equation*}
	
Since we assumed that $(z_{\tau(n)})_{n \in \NN}$ satisfies the walking-away uniformly property we already know that 
$\PP \left( d(z_{\tau(n)},x) - d(z_0,x) \le \varepsilon n \right) $ has an exponential decay to $0$, uniformly in $x$. \\ 

It remains then to show that 
	$$ \PP \left( d(x, z_{\tau(n)}) - d(x, z_{\beta n}) \ge \frac{\varepsilon n}{2} \right) $$
decreases exponentially fast in $n$, uniformly in $x$. We will actually show that for all $a > 0$ 
	$$ \PP \left( d(x, z_{\tau(n)}) - d(x, z_{\beta n}) \ge a n \right) $$
decreases exponentially fast, uniformly in $x$. By the triangle inequality we have  
$$ \somme{ i \in [\![\tau(n),\beta n]\!]  } d(z_i,z_{i +1}) \ge  d(z_{\beta n}, z_{\tau(n)}) \ge d(x, z_{\tau(n)}) -d(x, z_{\beta n}) \ , $$  
where, $[\![\tau(n),\beta n]\!]$  denotes the set of  natural numbers in the interval bounded by $\{ \tau(n), \beta n \}$. Therefore 
	$$ \PP \left(   \somme{ i \in [\![\tau(n),\beta n]\!]  } d(z_{i},z_{i +1 }) \ge a n \right) \ge  \PP \big(  d(x, z_{\tau(n)}) - d(x, z_{\beta n}) \ge a n \big) \ . $$ 
	Note that the left hand side does not depend on $x$ anymore. 
Define $Z_i := d(z_{i},z_{i +1}) $. The desired result will follow once we prove that for all $a > 0$
$$ \PP \left(  \somme{ i \in [\![\tau(n),\beta n]\!]  } Z_i \ge a n \right)$$ decreases exponentially. The above summation is a summation of I.I.D.\ random variables over a random time interval. In order to control it, we shall first control the random time with a large deviations estimate for I.I.D.\ random variables and conclude by controlling the summation using again a large deviations estimate for I.I.D.\ random variables. Recall that, by construction of the sampling, one has $\tau(n) = \sum_{ 0 \le i \le n-1} (\tau(i+1) - \tau(i)) $, the $(\tau(i+1) - \tau(i))$'s being I.I.D.\ distributed as $\tau(1)$ (in particular they have a finite exponential moment). \\
 
Let $\alpha > 0$ such that $ \alpha \cdot \EE(Z_1) \le a/4$. We use the large deviations estimate for $\tau(n)$ (which  is a summation of I.I.D.\ random variables with a finite exponential moment): let $c_1, c_2 > 0$ such that 
 	$$ \PP(  | \tau(n) - \beta n | \ge \alpha n ) \le c_1 \ e^{-c_2 n} \ . $$
Recall that $\beta$ is the mean of $\tau$. Therefore, 
	$$ \PP \left(   \somme{ i \in [\![\tau(n),\beta n]\!]  } Z_i \ge a n \right) \le  c_1 \ e^{-c_2 n}  +  \PP \left(  \left\{ \somme{ i \in [\![\tau(n),\beta n]\!]  } Z_i  \ge a n \right\} \cap \{ | \tau - \beta n | \le \alpha n \} \right) \ . $$

Since the $Z_i$'s are non-negative,  one has 
$$ \PP \left(  \left\{ \somme{ i \in [\![\tau(n),\beta n]\!]  } Z_i  \ge a n \right\} \cap \{ | \tau - \beta n | \le \alpha n \} \right) \le  \PP \left(   \somme{ (\beta - \alpha)n \le i \le (\beta + \alpha)n }  Z_i  \ge a n \right) \ . $$
We conclude rewriting the right member of the above inequality as 
 $$ \PP \left(   \somme{ (\beta - \alpha)n \le i \le (\beta + \alpha)n } ( Z_i  - \EE(Z_i)) \ge( a  - 2 \alpha \cdot \EE(Z_1) ) n \right) \ . $$

Recall that we chose $\alpha$ such that $ a  - 2\alpha \cdot \EE(Z_1) \ge \ \frac{a}{2}$. The $Z_i$'s are I.I.D.\ with a finite exponential moment. Hence they satisfy large deviations estimates and  the above probability decreases exponentially fast. \hfill $\blacksquare$ 

\subsection{Proof of Lemma \ref{lemassumptionschottkywalk}} 
By Lemma \ref{lemfiniteexpmomentsampling} we know that $ \EE (d(z_p, z_{\tau(1)})) $ is finite; besides, by construction of $\tau$, it does not depend on $p$. Therefore, by the triangle inequality and linearity of the expectation, Lemma  \ref{lemassumptionschottkywalk} will follow once we have proved that 
		$$ \limisup{p \to \infty} \ \EE(d(z_0,z_p)) = \infty \ .$$

We start noticing that, for any $R > 0$, the following stopping time 
$$ \tau_R := \inf \{ k \in \NN \ , \ d(z_0,z_k) \ge R \} \ $$
is almost surely finite. 
Indeed, there is at least one element $\gamma_0$ in the subsemigroup $\Gamma_{\mu}$ generated by $\supp{\mu}$ such that $\gamma_0 \cdot B(z_0,R) \cap B(z_0,R) = \emptyset \ $: recall we assumed that $\Gamma_{\mu}$ has unbounded orbits. Therefore there exists $k_0$ such that $\PP(\gamma_{k_0}=\gamma_0)>0$. With probability one, there will be infinitely many times $k$ such that 
$\gamma_k^{-1}\gamma_{k_0+k}=\gamma_0$. This last property implies that almost any path eventually leaves the ball of radius $R$ around $z_0$. \\

We conclude the proof of Lemma \ref{lemassumptionschottkywalk} with the following 

\begin{lemma}
	Let $\mu$ be a probability measure on a group $\Gamma$ acting by isometries on a metric space $X$ and $z_0 \in X$. If for any $R > 0$ the time $\tau_R$ is almost surely finite, then 
		$$ \limisup{p \to \infty} \ \EE(d(z_0,z_p)) = \infty \ .$$
\end{lemma}

\textbf{Proof.} We have for any $n \in \NN$ and any $R > 0$ 
\begin{align*}
	\PP(z_n \notin B(z_0, R)) & \ge \PP( \tau_{2R} \le n \ , \ d(z_{\tau_{2R}}, z_n) \le R) \\
	& \ge \somme{ 0 \le k \le n} \PP( \tau_{2R} = k  \ , \ d(z_{k}, z_n) \le R) \\
	& \ge \somme{ 0 \le k \le n} \PP( \tau_{2R} = k ) \cdot  \PP(d(z_{k}, z_n) \le R)  
\end{align*}
	since the event $ \tau_{2R} = k$, that only depends on the first $k$ increments of the walk and $d(z_{k}, z_n)$, that only depends on the later increments of the walk, are independent. 
The random variable $d(z_{k}, z_n)$ follows the same law as  $d(z_{0}, z_{n-k})$. Therefore 
\begin{align*}
	 \PP(z_n \notin B(z_0, R)) & \ge \somme{ 0 \le k \le n} \PP( \tau_{2R} = k ) \cdot  \PP(d(z_{0}, z_{n-k}) \le R) \\
	  & \ge \ \PP( \tau_{2R} \le n)  \infi{ 0 \le k \le n} \PP(d(z_{0}, z_{k}) \le R) \\
	  & \ge \ \PP( \tau_{2R} \le n)  (1 - \supr{ 0 \le k \le n} \PP(z_k \notin B(z_0, R)) ) \ . 
 \end{align*}

Let $a_{n} := \supr{0 \le k \le n} \PP(z_k \notin B(z_0, R)) \ge \PP(z_n \notin B(z_0, R))$. We have shown that  
	$$ a_n \ge  \ \PP( \tau_{2R} \le n)  (1 - a_n ) \ , $$ 

Recall we are assuming that $\tau_{2R}$ is almost surely finite. Therefore there exists   $n$  such that $\PP( \tau_{2R} \le n) \ge 1/2$. For such a $n$,  we get $a_n \ge \frac{1}{5}$. Therefore there exists $ 0 \le k \le n$ such that $\PP(z_k \notin B(z_0, R) ) \ge 1/5$,
which implies in particular for the same $k$ that $ \EE(d(z_0, z_k)) \ge R/5$
thus concluding the proof. 	 \hfill $\blacksquare$ $\blacksquare$
	
\section{Proof of Proposition \ref{propsamplingschottky}}
\label{secintertwine}
We prove the following more precise version of  Proposition \ref{propsamplingschottky}.

\begin{proposition} Let $S$ be a Schottky set and $\mu_1$ be a probability measure with a finite exponential moment such that 
	$$  \somme{\gamma \in \Gamma} \ \mu_1(\gamma) \ d(z_0, \gamma \cdot z_0)
 > 6C + 6 S_{\sup} $$
where $C$ is as in Definition \ref{defschottkyset} of a Schottky set given by Proposition \ref{lemmapingpong} and $S_{\sup} := \supr{s \in S} \ d(z_0, s \cdot z_0)$. Then the corresponding $(\mu_1,S)$-random walk has the walking-away uniformly property.
\end{proposition} 
 
\textbf{Proof.} In order not to burden the notations, we shall denote by $(Z_n)_{n \in \NN}$ (instead of $(z_n^{\mu_1, S})_{n \in \NN}$) the successive positions in $X$ of the $(\mu_1,S)$-random walk. To simplify a bit the exposition, let us first note that one can suppose \textbf{the even increments of the walk to be $\mu_1$-increments and the odd ones to be Schottky increments.} Indeed, since we assumed that $\mu_1$ has a finite exponential moment, the walking-away uniformly property does not depend on the  first increment of the walk. With the notation introduced in Part  \ref{seclinearprogress} to define the $(\mu_1,S)$ random walk, we have $Z_n = \Upsilon_n \cdot z_0$ with
$$ \Upsilon_n := 
		\left\{ 
			 \begin{array}{lr} 
			   Y_1 \cdot X_2 \cdot ... \cdot X_{n}  & \text{ if } n \ \text{ is even } \\ 
			 	Y_1 \cdot X_2 \cdot ... \cdot Y_{n}  & \text{ if } n \ \text{ is odd } \ . 
		 	\end{array}
		 \right.	$$ 

We start with the obvious equality  
	$$ \PP \Big( d(Z_{2n}, x) - d(z_0,x) \le \varepsilon n \Big) = \PP \Big( \somme{0 \le i \le n-1} d(Z_{2i +2}, x) - d(Z_{2i},x) \le \varepsilon n \Big) \ . $$
	
For any $x,y,z \in X$, we let  
	$$ B_x(z,y) := d(z, x) - d(y,x) \ , $$
	
so that 
	$$\somme{0 \le i \le n-1} d(Z_{2i +2}, x) - d(Z_{2i},x) = \somme{0 \le i \le n-1} B_{x}(Z_{2i+2}, Z_{2i}) := S_n(x) \ . $$

Using Markov inequality for a small enough $\lambda>0$, we get that 
	$$ \PP \Big( S_n(x) \le \varepsilon n \Big) \le e^{ \lambda \varepsilon n} \ \EE  \Big( e^{ -\lambda  S_n(x)} \Big) . $$ 
We will be done once we prove that there exist $\lambda > 0$ and $ 0 < \delta < 1$ (which may depend on $\lambda$) such that for all $x$ 
	\begin{equation}\label{uneautreequation} \EE  \Big( e^{ -\lambda  S_n(x)} \Big) \le \delta^n \ . \end{equation} 

Recall that we denoted by $(\mathcal{F}_n)_{n \in \NN}$ the filtration of $\Omega$ with respect to the increments of the walk. Conditioning on $\mathcal{F}_{2n -2}$, we have
\begin{align*}
\EE  \Big( e^{ -\lambda  S_n(x)} \Big)  = \EE  \Big( \EE  \big( e^{-\lambda S_{n-1}(x)}  \cdot  e^{- \lambda \left( B_{x}(Z_{2n}, Z_{2(n-1)}) \right)} \big| \mathcal{F}_{(2n-2)} \big) \Big) \  
 \end{align*}

\begin{align*}  
	=\EE  \Big( \EE  \big( e^{-\lambda S_{n-1}(x)}  \cdot  \exp \left( - \lambda B_{\Upsilon^{-1}_{2n-2} x}( \Upsilon^{-1}_{2(n-1)} Z_{2n},z_0) \right) \big| \mathcal{F}_{(2n-2)} \big) \Big) \  
 \end{align*} 
since $S_{n-1}(x)$ is $\mathcal{F}_{2n -2}$ measurable. 

Because $\Upsilon^{-1}_{2n-2}$ is $\mathcal{F}_{2n -2}$ measurable and $\Upsilon^{-1}_{2(n-1)} Z_{2n}$ is independent of $\mathcal{F}_{2n -2}$, we have 
\begin{align*}
	\EE  \big( \exp & \left( - \lambda B_{\Upsilon^{-1}_{2n-2} x}( \Upsilon^{-1}_{2(n-1)} Z_{2n},z_0) \right) \big| \mathcal{F}_{(2n-2)} \big) \\
	& \le \supr{y \in X} \ \EE  \big( \exp  \left( - \lambda B_{y}( \Upsilon^{-1}_{2(n-1)} Z_{2n},z_0) \right) \big| \mathcal{F}_{(2n-2)} \big) \\
	& = \supr{y \in X} \ \EE  \big( \exp  \left( - \lambda B_{y}( Z_{2},z_0) \right) \big| \mathcal{F}_{(2n-2)} \big) \ .
\end{align*}

We get by an immediate induction that
	$$ 	 \EE  \Big( e^{ -\lambda  S_n(x)} \Big) \le \delta(\lambda)^n , $$
where 
	$$ \delta(\lambda) := \supr{y \in X} \  \EE  \left( \ e^{- \lambda B_{y}(Z_2,z_0)} \right) \ .$$

We end this proof by showing  
	\begin{lemma}
		\label{lemproplinear}
There is $\lambda > 0$  such that
	$$\supr{y \in X} \ \EE  \left( \  e^{- \lambda \left( B_{y}(Z_2, z_0) \right)} \right) < 1 \ . $$
	\end{lemma}	
	
\textbf{Proof.}  We denote by $A^{c}$ the complement of a set $A$. Given $y \in X$, we use the decomposition  
 $$ B_{y}(Z_2, z_0) = B_{y}(Z_2, z_0) \ \mathds{1}_A + B_{y}(Z_2, z_0) \ \mathds{1}_{A^c} \ ,$$

where $A := \{ (Z_2, y)_{z_0} \le C\}$ and $C$ is the constant given by Proposition \ref{lemmapingpong}. \\

Note that on $A$,  since the first increment of the walk is in $S$, we have 
	\begin{align*} 
		B_{y}(Z_2, z_0) & \ge d(Z_2, z_0) - 2C  \\
			& \ge d(Z_2, Z_1) - 2C - d(Z_1, z_0) \\
			& \ge d(Z_2, Z_1) - 2C - S_{\sup} \ .
	\end{align*}

On $A^c$ we use the trivial lower bound
	$$ B_{y}(Z_2, z_0) \ge - d(Z_2, z_0) \ge - d(Z_2, Z_1) - S_{\sup} \ . $$ 

We thus obtain the inequality  
\begin{align*}
	 e^{-\lambda B_{y}(Z_2, z_0)} & \le e^{ -\lambda (d(Z_2, Z_1) - 2C - S_{\sup} )} \ \mathds{1}_A + e^{ - \lambda ( - d(Z_2, Z_1) - S_{\sup})} \ \mathds{1}_{A^c} \\
	 	& = e^{\lambda S_{\sup}} \cdot \Big( e^{ -\lambda ( d(Z_2, Z_1) - 2C) } \ \mathds{1}_A + e^{\lambda  d(Z_2, Z_1) } \ \mathds{1}_{A^c}  \Big) \ .
\end{align*}

Since the distances appearing in the exponentials do not depend on $Y_1$ but only on $X_2$, we have \begin{align*}
	 \EE  \left( e^{-\lambda B_{y}(Z_2, z_0)} \big| X_2 \right) & \le e^{\lambda S_{\sup}} \cdot \Big( e^{ -\lambda ( d(Z_2, Z_1) - 2C) } \ \PP( A \ | \ X_2) + e^{\lambda  d(Z_2, Z_1) } \ \PP( A^c \ | \ X_2) \Big) \\ 
	 & \le e^{\lambda S_{\sup}} \cdot \Big( e^{ -\lambda ( d(Z_2, Z_1) - 2C) } \ \alpha(X_2) + e^{\lambda  d(Z_2, Z_1) } \ (1 - \alpha(X_2)) \Big) \ ,
\end{align*}

where we set 
	\begin{align*}
		 \alpha(X_2)  := \PP( A \ | \ X_2)   =  \PP( (Z_2, y)_{z_0} \le C \ | \ X_2) 
		 	= \frac{ \# \ \{s \in S \ , \  (s \cdot X_2 \cdot z_0,  y)_{z_0} \le C  \ \}}{ \# S} \ ,
 	\end{align*} 
using the fact that the first increment is uniformly distributed on $S$. Because $S$ is a Schottky set, we readily get that 
\begin{equation}
	\label{eqref1000}
	 \alpha(X_2) \ge \frac{2}{3} \ ,
\end{equation}
for all $y$. \\ 

Next, we use the lower bound on $\EE  (d(Z_2, Z_1))$ to argue that, in the upper bound above, 
out of the two competing exponentials, the main contribution comes from the term $e^{ -\lambda ( d(Z_2, Z_1) - C) }$. \\

Recall the general upper bound, $e^x \le 1 + x + x^2e^{\vert x\vert}$  and set 
	$$ R(\lambda, X_2) := \lambda^2 \ (d(Z_2,Z_1) +2C)^2 e^{ \lambda d(Z	_2,Z_1)}\ . $$
We then estimate  \begin{align*}
	  & \EE   \left(  e^{-\lambda B_{y}(Z_2, z_0)} \big| X_2 \right)  \\ 
	  &\le   e^{\lambda S_{\sup}} \cdot \Big( \alpha(X_2)-\lambda  ( d(Z_2, Z_1) - 2C)\ \alpha(X_2) +1-\alpha(X_2) +   \lambda  d(Z_2, Z_1)\ (1-\alpha(X_2)) +e^{2\lambda C}       R(\lambda, X_2)                                 \Big)\\
	  &=   e^{\lambda S_{\sup}} \cdot \Big(1-\lambda   d(Z_2, Z_1) \ (2\alpha(X_2)-1) +   2\lambda  C\alpha(X_2) +e^{2\lambda C}       R(\lambda, X_2)                                                  \Big)\\
	  &\le   e^{\lambda S_{\sup}} \cdot \Big(  1-\frac 13\lambda   d(Z_2, Z_1)  +   2\lambda  C +e^{2\lambda C}       R(\lambda, X_2)                                                                                        \Big)\ .
	  \end{align*}
We used the bound \eqref{eqref1000} and the fact that $\alpha(X_2)\leq 1$. Taking the expectation in this last inequality and using the lower bound on $\EE(d(Z_2, Z_1))$, we get that 
	 $$  \EE   \left(  e^{-\lambda B_{y}(Z_2, z_0)} \right)  \le 
	 e^{\lambda S_{\sup}} \Big(  1-2\lambda   S_{\sup}   +e^{2\lambda C}       \EE(R(\lambda, X_2))                                                                                        \Big)\ .
$$

Choose $ \lambda_0>0$ such that
	$$C_1:= \EE \left( (d(Z_1,Z_2)+C)^2 \ e^{ \lambda_0 d(Z_1,Z_2)} \right) < \infty \ .$$ 
	Then, for all $\lambda\le\lambda_0$, 
	 $$  \EE(R(\lambda, X_2)) \le  C_1 \lambda^2  \  $$
and  we get that 
	 $$  \EE   \left(  e^{-\lambda B_{y}(Z_2, z_0)} \right)  \le 
e^{\lambda S_{\sup}} \Big(  1-2\lambda   S_{\sup}   +C_1\lambda^2\ e^{2\lambda C}      \Big)\ .$$
The right hand side of this last inequality is $<1$ for some positive but small enough $\lambda$ and this completes the proof. \hfill $\blacksquare$ $\blacksquare$

\subsection{The finite first moment case}	
\label{subsec.finitefirstmoment}
As emphasised in the introduction, one can adapt Sections \ref{seclinearprogress} and \ref{secintertwine} to the setting where the measure $\mu$ has only a finite first moment to recover that $l >  0$ in this setting. \\

The general strategy is entirely the same, in particular the exact same sampling is to be performed. The only things to be modified are the statements of the various lemmas appearing in Subsection \ref{subsec.overviewargument}. We will not give all the details since it mainly repeats previously given arguments. Let us however indicate to the reader the changes and the non changes that one should perform to get positivity of the escape rate under a finite first moment. \\

Under the assumption that $\mu$ has a finite first moment, the proof of Proposition \ref{propsampling} gives 
	$$ \limi{n \to \infty} \ \frac{1}{n} \EE ( d(z_{\beta n},z_{\tau(n)}))=0 \ ,  $$ 
which implies in particular that the random walk $(z_{\beta n})$ has positive escape rate if and only if the random walk $(z_{\tau(n)})$ has positive escape rate. \\

Lemma \ref{lemfiniteexpmomentsampling} should be rephrased by replacing every occurrences of 'finite exponential moment' with 'finite first moment'. \\

Lemmas \ref{lemsamplingtoschottky} and \ref{lemassumptionschottkywalk} are identical (they do not require a finite exponential moment). \\

The assumption of Proposition \ref{propsamplingschottky} is to be modified with the assumption that $\mu_1$ has a finite first moment. Its conclusion should be replaced with 'the $(\mu_1,S)$-random walk has positive escape rate'. The proof is even simpler in this case. Indeed, using the notations previously used, we start with the same decomposition but with taking the expectation:
	$$ \EE \Big( d(Z_{2n}, x) - d(z_0,x)\Big) =  \somme{0 \le i \le n-1} \EE \Big( B_x(Z_{2i +2}, Z_{2i})\Big) \ . $$
	
We then skip all the Markov Inequality/conditioning to go directly to the following modified version of Lemma \ref{lemproplinear}, which shows that $l>0$ by taking $z_0 = x$ in the above identity. 

	\begin{lemma}
		\label{lemproplinear2}
There is $c > 0$  such that
	$$\inf_{y \in X} \ \EE \left( B_{y}(Z_2, z_0) \right) > c \ .$$
	\end{lemma}	
	
The proof follows the same lines as in the proof of Lemma \ref{lemproplinear}.

\begin{remark}
One could be even more precise and get the following weak walking-away uniformly property (compare with \cite[Definition 1.4]{EL}). \\

\textit{Let $\Gamma$ be a countable group acting by isometries on a geodesic Gromov-hyperbolic space $X$, $\mu$ an admissible probability measure on $\Gamma$ with a finite first moment and $z_0 \in X$.  Then there is a constant $c > 0$ such that for all $x \in X$ and for all $n \in \NN$}
	$$ \EE \left(  d(z_n,x) - d(z_0,x) \right) \ge c  \ n\ .$$ 
\end{remark}

\section{Deviation inequalities}
\label{secdevlinearprog}
Recall that the walking-away uniformly property, treated in the previous two sections, directly implies linear progress with exponential tail (Definition \ref{defslinearprogress}). The goal of this section is to show that a random walk which satisfies linear progress with exponential tail also satisfies the following property.

\begin{definition}{\cite{artmathieusisto}}
Let $(z_n)_{n \in \NN}$ be a random path in a metric space $X$. We say that $(z_n)_{n\in\NN}$ satisfies the \textbf{exponential-tail deviation inequality} if there are constants $ C_1, C_2 > 0$ such that for all $ 0 \le i \le n$ and all $R > 0$ one has
	$$ \PP((z_n,z_0)_{z_i} \ge R) \le C_1 \ e^{ -C_2 R} \ .  $$
\end{definition}

We adapt the proof of \cite[Theorem 11.1]{artmathieusisto} to prove the following

\begin{proposition}
	\label{proplinearprogresstodev}
Let $\Gamma$ be a countable group acting by isometries on a geodesic Gromov-hyperbolic space $X$ and $\mu$ a probability measure on $\Gamma$ with a finite exponential moment. If the random walk has linear progress with exponential tail then it satisfies the exponential-tail deviation inequalities.
\end{proposition}

\begin{remark}
1) We note that, unlike Theorem \ref{theowalkingaway}, the previous proposition assumes, among others, that $X$ is a Gromov-hyperbolic metric space.\\[2pt] 
2) In the case where $\Gamma$ acts acylindrically on a geodesic Gromov-hyperbolic space, this proposition is already proved in \cite[Theorem 10.7]{artmathieusisto}.
\end{remark}

\textbf{Proof.} Given a geodesic $\upsilon$, we denote by $\pi_{\upsilon}$ a choice of nearest point projection from $X$ to $\upsilon$. Given two points $x,y \in X$ we denote by $[x,y]$ the choice of any geodesic path joining $x$ to $y$. Given any $y \in X$ we define 
	$$ N_{\upsilon}(y) := \{ x \in X \ , \ d(\pi_{\upsilon}(y),\pi_{\upsilon}(x)) \ge d(x, \pi_{\upsilon}(x)) \} \ . $$  
Note that the above set actually depends on $\pi_{\upsilon}(y)$ only and, in particular, not on $d(y, \upsilon)$.  We refer to \cite[Section 2]{artmahertiozzo} and \cite[Section 3]{artmahercomplexcurve} for more details about the nearest point projection. \\

Let $I$ be an interval of $\ZZ$ and $(x_i)_{i \in I}$ be a discrete path whose endpoints lie on the geodesic $\upsilon$. Given $k \in I \setminus \partial I$ we define
\begin{equation*}
	\left\{
		\begin{array}{l}
			k_1 := \sup \ \{ j < k  \ , \ [x_j,x_{j+1}] \cap N_{\upsilon}(x_k) \neq \emptyset \} \   \\
			k_2 := \inf \ \{ j \ge k  \ , \ [x_j,x_{j+1}] \cap N_{\upsilon}(x_k) \neq \emptyset \} \ .
		\end{array} 
	\right.
\end{equation*} 

\begin{figure}[h!]
\begin{center}
	\def\svgwidth{1 \columnwidth}
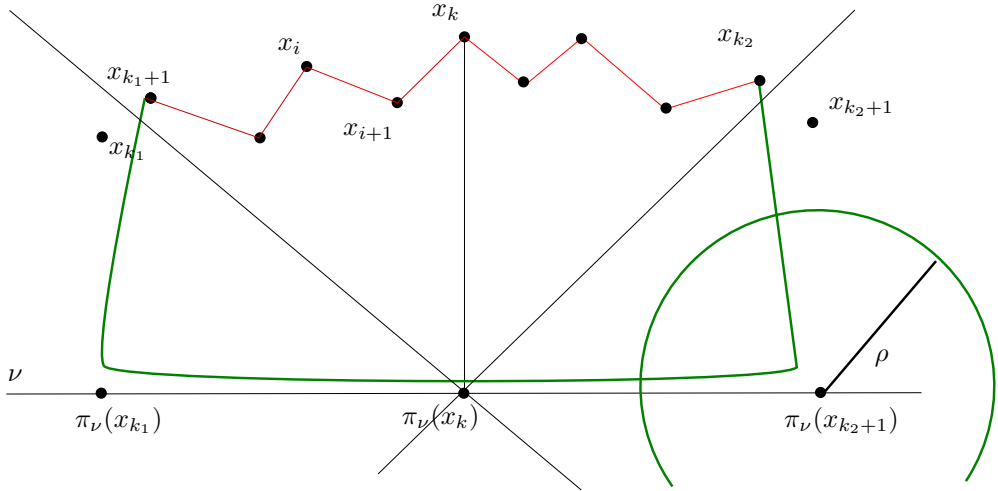
	\end{center}
\caption{The green path represents the geodesic from $x_{k_1+1}$ to $x_{k_2}$. Because of Items (2) and (3) of Lemma \ref{lemmathieussistohyp} their projections on $\upsilon$ must remain close to those of $x_{k_1}$ and $x_{k_2 +1}$. From Item (1) and by construction of $N_{\upsilon}(x_k)$ one of the distances $ d ( \pi_{\upsilon}(x_{k_1}), \pi_{\upsilon}(x_k))$, $ d (\pi_{\upsilon}(x_{k_2 +1}), \pi_{\upsilon}(x_k))$ must be at least of the order $\varepsilon( k_2 - k_1)/ 5$. This prevents the red path to enter at least one of the balls $B(\pi_{\upsilon}(x_{k_1},\rho)$ or $B(\pi_{\upsilon}(x_{k_2 +1},\rho)$, in green in the picture (with $\rho = \varepsilon( k_2 - k_1)/ 10$). This implies that the length of the red path must be exponential in $\rho$ since it avoids a ball through which the geodesic connecting its endpoints passes (in green in the picture).}
\label{figmathieusistogeomlem}
\end{figure}

Note that $d(x_{k_i},x_k) \ge d(x_k, \pi_{\upsilon}(x_k)) - 100 \delta $ for $i= 1,2$ (see \cite[Lemma 11.4, Claim 1]{artmathieusisto} and Figure \ref{figmathieusistogeomlem}). The following lemma is the geometric key of the proof.

\begin{lemma}{\cite[Lemma 11.4]{artmathieusisto}}	
\label{lemmathieussistohyp}
For any $\varepsilon > 0$ there are constants $c_1, c_2 > 0$ such that if 
\begin{enumerate}
	\item		$ d(x_{k_1}, x_{k_2}) \ge \varepsilon \ (k_2 - k_1) \ ; $
 	\item	 	$ d(x_{k_1}, x_{k_1 +1}) \le d(x_{k_1}, x_{k_2})/100 \ ; $
 	 \item	 	$ d(x_{k_2}, x_{k_2 +1}) \le d(x_{k_1}, x_{k_2})/100 \ , $

\end{enumerate}
then 
	$$ \somme{ i \in [k_1+1, \ k_2-1]} d(x_{i}, x_{i +1}) \ge c_2 \ e^{ c_1 (k_2 - k_1)} \ .$$ 
\end{lemma}

The statement above is a simplified version  of \cite[Lemma 11.4]{artmathieusisto}. The proof follows the same line and is illustrated in Figure \ref{figmathieusistogeomlem}. \\

We now use Lemma \ref{lemmathieussistohyp} with the successive positions of the random walk $(x_i = z_i)$. Recall that we want to show that there are constants $C_1, C_2$ such that for any $n > k > 0$ and any $R > 0$, we have
	$$ \PP( (z_n, z_0)_{z_k} > R) \le C_1 e^{-C_2 R} \ .$$ 

We fix $k,n$ and $R > 0$. For a path $(z_j)$ satisfying $(z_n, z_0)_{z_k} > R$, we define the times $k_1, k_2$ as in Lemma \ref{lemmathieussistohyp}. \\

We distinguish two cases, depending on whether or not $k_2 -k_1$ is large with respect to $R$. \\

The next lemma addresses the case of paths with a small value for $k_2 - k_1$. 

\begin{lemma}
	\label{lempropdev1}
There are constants $c_3, c_4, C > 0$ (independent of $k,n$ and $R$) such that
	$$ \PP( (z_n, z_0)_{z_k} > R \ , \ k_2 - k_1 \le c_3 R) \le C e^{- c_4 R} \ .$$
\end{lemma}

\textbf{Proof.} We will look at all the possible values of $k_1, k_2$ and conclude using the union bound. \\

Since we assumed that $(z_n, z_0)_{z_k} > R$ and by construction of $N_{\upsilon}(z_k)$, we have that $d(z_{k_1},z_k) \ge R - 100 \delta$ \cite[Lemma 11.4, Claim 1]{artmathieusisto}. Let $ 0 \le m < c_3 R$ for some $c_3$ that we will fix later on. Choose $\alpha < k$ and $\beta > k$ such that $\beta - \alpha = m$. We have 
\begin{align*}
	 \PP \big( (z_n, z_0)_{z_k} > R \ , \ k_2 = \beta \ , \  k_1 = \alpha  \big) & \le \PP \left( d(z_{\alpha}, z_k) \ge R - 100 \delta \right) \ .
\end{align*}

Using the triangle inequality we get 
\begin{align*}
	 \PP \left( d(z_{\alpha}, z_k) \ge R- 100 \delta \right) 			& \le \PP \left(  \somme{ \alpha \le i \le k-1 } d(z_{i}, z_{i+1}) \ge R - 100 \delta \right) \\ 
	 	& \le \PP \left(  \somme{ \alpha \le i \le \beta - 1 } d(z_{i}, z_{i+1}) \ge R - 100 \delta \right) \\ 
	 	& = \PP \left(  \somme{ 0 \le i \le m-1 } d(z_{i}, z_{i+1}) \ge R - 100 \delta \right) \ . 
\end{align*}

Taking the Laplace transform and using Markov's inequality we get that, for all $\lambda > 0$,  
\begin{align*}
 	\PP \left(  \somme{0 \le i \le m -1 } d(z_{i}, z_{i+1}) \ge R - 100 \delta \right) & \le e^{ -\lambda(R - 100 \delta)} \ \left( \EE  \left( e^{\lambda d(z_0,z_1)} \right) \right)^{m} \\ 
 	& \le C \ e^{ -\lambda R} \ \left( \EE  \left( e^{\lambda d(z_0,z_1)} \right) \right)^{c_3 R} \ .
 	 	\end{align*}

From this last inequality, provided we choose $\lambda$ such that $\EE  \left( e^{\lambda d(z_0,z_1)} \right)<\infty$ and $c_3$ small enough, we deduce that 
$$ \PP \big( (z_n, z_0)_{z_k} > R \ , \ k_2 = \beta \ , \  k_1 = \alpha  \big)\le Ce^{-c R}\ ,$$
for some constants $C$ and $c>0$. 
The lemma now follows by summing over the possible choices of $\beta$ and $\alpha$. \hfill $\blacksquare$ \\

The next lemma deals with the remaining case corresponding to $k_2 - k_1 \ge c_1 R$ and concludes the proof of Proposition \ref{proplinearprogresstodev}.  

\begin{lemma}	For any $c_1>0$, there are constants $c_5, c_6$ (independent of $k,n$ and $R$) such that we have
	 	$$ \PP( (z_n, z_0)_{z_k} > R \  , \ k_2 - k_1 \ge c_1 R ) \le c_5 \ e^{- c_6 R} \ .$$
\end{lemma}

\textbf{Proof.} We shall prove that, for any $m > 0$, then 
	 	$$ \PP( (z_n, z_0)_{z_k} > R \  , \  k_2 - k_1 = m ) \le c_5 \ e^{- c_6 m} \ .$$
The lemma follows by summing over all $ m \ge c_1 R$ (with slightly different values for $c_5$ and $c_6$). \\

Let us then fix $m > 0$. In the same way as for the proof of Lemma \ref{lempropdev1} we  first fix $k_1 = \alpha$ and $k_2 = \beta$ with 
$\beta - \alpha = m$ and then use the union bound. \\

Recall that since we assumed that the walk has linear progress with exponential tail one has constants $\varepsilon, c_7, c_8 > 0$ such that 
$$ \PP \left( d(z_{\alpha},z_{\beta}) \le \varepsilon m \right) \le c_7 e^{-c_8 m} \ . $$

There are also constants $c_9,c_{10}$ such that
\begin{align*}
	 \PP \left( d(z_{\alpha},z_{\alpha +1 }) \ge \frac{d(z_{\alpha},z_{\beta})}{100} \right) & \le 	 \PP \left( d(z_{\alpha},z_{\alpha +1}) \ge \frac{\varepsilon m }{100} \right)+  c_7 e^{-c_8 m} \\
	& \le c_9 e^{-c_{10} m} \ ,
\end{align*}

since $d(z_{\alpha},z_{\alpha +1})$ has a finite exponential moment. A similar bound applies to $d(z_{\beta},z_{\beta +1 })$.\\ 

It remains to estimate the probability of the event, say $A$, when  $d(z_\alpha,z_\beta)\le \varepsilon m$, $d(z_{\alpha},z_{\alpha +1})\le d(z_\alpha,z_\beta)/100$ and $d(z_{\beta},z_{\beta +1 })\le d(z_\alpha,z_\beta)/100$.  

According to Lemma \ref{lemmathieussistohyp}, on $A$, one has 
	$$ \somme{ \alpha \le i \le \beta} d(z_{i}, z_{i+1}) \ge c_2 \ e^{ c_1 m} \ .$$ 
The probability of the above event is (super)-exponentially small in $m = \beta - \alpha$. 
\hfill $\blacksquare$ $\blacksquare$

\section{Hitting measure}
\label{sechittingmeasure}

The purpose of this section is to prove the uniform punctual deviations Proposition \ref{theoedpi} from the introduction. We recall its statement:

\begin{proposition}[uniform punctual deviations]
\label{prophittingmeasuresec8}
Let $\Gamma$ be a countable group acting by isometries on a geodesic Gromov-hyperbolic space $X$ and $\mu$ a non-elementary probability measure on $\Gamma$.  Then, there are constants $C, \alpha >0$ such that for any $p \in \NN$ and any $x \in X$, $R>0$ we have 
\begin{equation}\label{eq.unif.punct.dev}
\PP( (z_p,x)_{z_0} \ge R ) \le C e^{-\alpha R} \ .
\end{equation}
\end{proposition}
The above proposition implies that Assumption \eqref{eqassumptionproplower} holds since it implies that for any $x \in X$ 
	$$ \EE ( (z_p,x)_{z_0}) \le \frac C\alpha \ .$$

\begin{remark}
If we further assume $\mu$ is symmetric, then there is an easy way to deduce Proposition \ref{prophittingmeasuresec8} from Proposition \ref{theoremdevineq}. Indeed, let us rewrite as follows the square of the quantity we want to bound  
\begin{align*}
	 \PP_{z_0} \left( (z_m,x)_{z_0} \ge R \right)^2 & = \PP_{z_0} \left( (z_m,x)_{z_0} \ge R \right) \cdot \PP_{z_0} \left( (\widehat{z}_m,x)_{z_0} \ge R \right) \\
	 & = \PP_{z_0} \left( (z_m,x)_{z_0} \ge R \ , \  (\widehat{z}_m,x)_{z_0} \ge R \right)  
\end{align*}
where $\widehat{z}_m$ is an independent copy of $z_m$. 
The hyperbolicity of $X$ implies that, for any four points $(x_i)_{0 \le i \le 3}$ such that $ (x_1, x_2)_{x_0} \ge R$ and $(x_2, x_3)_{x_0} \ge R$, we have 
	$$ (x_3, x_1)_{x_0} \ge \min(R, R) - \delta = R - \delta \ .$$
Therefore 
$$ \PP_{z_0} \left( (z_m,x)_{z_0} \ge R \ , \  (\widehat{z}_m,x)_{z_0} \ge R \right) \le \PP_{z_0} \left( (z_m, \widehat{z}_m)_{z_0} \ge R - \delta \right)\ .$$	
Because we assumed the measure $\mu$ to be symmetric, the random variable $(z_m, \widehat{z}_m)_{z_0}$ has the same law as $(z_{2m}, z_0)_{z_m}$. Therefore, we have 
$$\PP_{z_0} \left( (z_m, \widehat{z}_m)_{z_0} \ge R - \delta \right)= \PP_{z_0} \left( (z_{2m}, z_0)_{z_m} \ge R - \delta \right) $$
and hence \eqref{eq.unif.punct.dev} follows from the exponential-tail deviation inequality  (Proposition \ref{theoremdevineq}).
\end{remark}

\textbf{Proof of Proposition \ref{prophittingmeasuresec8}} We will use the walking-away property from Theorem \ref{theoremwau},  the linear progress property from Definition \ref{defslinearprogress}  and exponential-tail deviation inequality from Proposition \ref{theoremdevineq}.\\

The geometric key of the proof is 
\begin{lemma}
\label{lemgeomdemhitting}
Let $z,x,q,y \in X$. There is $R_0 = R_0(\delta) > 0$ such that for every $ R \ge R_0$ if
\begin{itemize}
	\item  $(x,q)_z \ge R \ ;  $
	\item $	\frac{4R}{5} \le d(z,y) \le R \ ; $
	\item $	(z,q)_y \le \frac{R}{5} \ ,$ 
then 
	$$ d(z,x) - d(y,x) \ge \frac{R}{10} \ . $$
\end{itemize}
\end{lemma}

\begin{figure}[h!]
\vspace*{-0.6cm}
\begin{center}
	\def\svgwidth{0.64 \columnwidth}
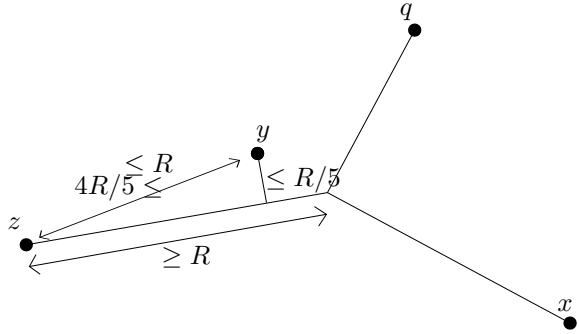
	\end{center}
	\vspace*{-1cm}
\caption{The geodesics from $z$ to $q$ and $x$ fellow-travel for at time at least $R$. 
}
\label{figlemgeosimple}
\end{figure}

The proof of the above lemma is also illustrated in Figure \ref{figlemgeosimple}.\\

\textbf{Proof.}
Since $d(z,y) \geq \frac{4R}{5}$ and $(z,q)_y \leq \frac{R}{5}$, by expanding $(z,q)_y$, we get that $d(q,z)-d(q,y)\geq \frac{2R}{5}$. Therefore, using once more that $d(y,z)\geq \frac{4R}{5}$, we obtain $(q,y)_z \geq \frac{3R}{5}$. Using this, the hypothesis $(x,q)_z \geq R$ and $(x,y)_z \geq \min\{(x,q)_z,(q,y)_z\} -\delta$, we get $(x,y)_z \geq \frac{3R}{5}-\delta$. Expanding $(x,y)_z$ and using $d(z,y) \leq R$, we obtain $d(z,x)-d(y,x) \geq \frac{R}{5}-2\delta$ and hence $d(z,x)-d(y,x) \geq \frac{R}{10}$ for all $R$ large enough. \hfill $\blacksquare$

Let us see how to use Lemma \ref{lemgeomdemhitting} (with $z =z_0, x=x, q =z_p$ and with $y = z_i$ for some $i$) to get Proposition \ref{prophittingmeasuresec8}. \\

For a real number $R>0$, we denote $A_R := \{ d(z_{\beta R},z_0) \le R \} $. 
The linear progress with exponential tail property implies there exists $\beta, c_1 > 0$ such that  $\PP(A_R) \le c_1^{-1} \ e^{- c_1 R}$ for every $R>0$. \\

Since $\mu$ has finite exponential moment, using large deviations estimates for I.I.D.\ random variables we know that there is $\alpha, c_2 > 0$ such that for every $R>0$ and  for $ 0 \le j \le \alpha R$ we have
 $$ \PP \left( \somme{ 1 \le i \le j} d(z_i, z_{i +1}) \ge R/2 \right) \le c_2^{-1} \ e^{-c_2 R} \ . $$
In particular, denoting $B_R := \{ \exists i \in [0, \alpha R] \ , \ d(z_{i},z_0) \ge R/2 \} $ and using the triangle inequality together with the union bound, it gives a constant $c_3 > 0$ such that for every $R>0$, we have $\PP(B_R) \le c_3^{-1} \ e^{- c_3 R}$. Note that we may choose $\alpha$ so that $\alpha\le\beta$. \\

Using that $\mu$ has a finite exponential moment and the union bound, we get a constant $c_4 > 0$ such that, denoting $C_R := \{ \exists i \in [0, \beta R] \ , \ d(z_{i},z_{i+1}) \ge R/5 \} $, we have $\PP(C_R) \le c_4^{-1} \ e^{-c_4 R}$ for every $R>0$. \\

Finally, using the exponential-tail deviation inequality, the union bound and denoting $D_{p,R} := \{ \exists i \in [0, \min(\beta R,p)] \ , \ (z_{p},z_0)_{z_i} \ge R/5 \}$, we get a constant $c_5 > 0$ such that $\PP(D_{p,R}) \le c_5^{-1} e^{-c_5 R}$ for every $p \in \mathbb{N}$ and $R>0$. \\


It now remains to prove that there is a constant $c_6 > 0$ such that for all $p \in \NN$, $x \in X$ and $R>0$, we have
\begin{equation}\label{eq.new.label1}
\PP(\{(z_p,x)_{z_0} \ge R \} \cap A_R^c \cap B_R^c \cap C_R^c \cap D_{p,R}^c) 
 \le c_6^{-1} \ e^{-c_6 R} \ .
\end{equation}
 
Let $p \in \mathbb{N}$, $x \in X$ and $R>0$ be given. Note that a path in $\{(z_p,x)_{z_0} \ge R \} \cap A_R^c \cap B_R^c \cap C_R^c$ is such that all the steps are of length at most $R/5$; the path remains in  $B(z_0, R/2)$ for the first $\alpha R$  steps but is outside the ball $B(z_0, R)$ at time $\min(p,\beta R)$. In particular suppose $p\ge \alpha R$, for otherwise the left-hand-side of \eqref{eq.new.label1} is zero.
It follows that the event $\{(z_p,x)_{z_0} \ge R \} \cap A_R^c \cap B_R^c \cap C_R^c $ is contained in the set 
	$$ \{ \exists i \in [\alpha R, \min(p,\beta R)] \ , \  4 R /5 \le d(z_0, z_i) \le R \} \ . $$

Note also that any path in $D_{p,R}^c$ must satisfy that for any $ 0 \le i \le \min(p, \beta R)$, 
	$$ (z_p, z_0)_{z_i} \le R/5 \ . $$ 
	
Therefore, the event 
$$ \{(z_p,x)_{z_0} \ge R \} \cap A_R^c \cap B_R^c \cap C_R^c \cap D_{p,R}^c $$

is contained in the event 
\begin{align*}
	 \{(z_p,x)_{z_0} \ge R \} \cap  
	 	 \{ \exists i \in [\alpha R, \min(p,\beta R)] \ , \  4 R /5 \le d(z_0, z_i) \le R \ , \  (z_p, z_0)_{z_i} \le R/5 \} \ .
\end{align*}

Using the union bound, 
\begin{align*}\PP\Big(  \{(z_p,x)_{z_0} \ge R \} \ \cap \ 
	 	 \{ \exists i \in [\alpha R, \min(p,\beta R)] \ , \ (4 R /5 \le d(z_0, z_i) \le R \ , \  (z_p, z_0)_{z_i} \le R/5) \}  \Big) 
		 \end{align*} 
$$ \le \somme{ \alpha R \le i \le \beta R } \PP \left( (z_p,x)_{z_0} \ge R \ , \  4 R /5 \le d(z_0, z_i) \le R \ , \  (z_p, z_0)_{z_i} \le R/5 \right) \ .$$

To conclude the proof, we show that there is a constant $c_7>0$ (independent of $p$, $x$ and $R$) such that for all $\alpha R \le i \le \beta R$ we have 
	$$ \PP \left( (z_p,x)_{z_0} \ge R \ , \  4 R /5 \le d(z_0, z_i) \le R \ , \  (z_p, z_0)_{z_i} \le R/5 \right) \le c_7^{-1} \ e^{-c_7 R} \ .$$

Using Lemma \ref{lemgeomdemhitting}, we get that, for any $\alpha R \le i \le \beta R$, we have  
\begin{align*}
	 \PP ( (z_p,x)_{z_0} \ge R \ , \  4 R /5 \le  d(z_0, z_i) \le R \ , \  (z_p, z_0)_{z_i} \le R/5 ) & \le \PP \left( d(z_0,x) - d(z_i,x) \ge \frac{R}{10} \right) \\
		& \le \PP \left( d(z_0,x) - d(z_i,x) \ge 0 \right) \ . 
\end{align*}

Using the walking-away property, we get a constant $c_7 > 0$ (independent of $x$ and $R$, and clearly, also of $p$) such that for every $\alpha R \le i \le \beta R $
\begin{align*}
	 \PP \left( d(z_0,x) - d(z_i,x) \ge 0 \right) & \le c_7^{-1} e^{-c_7 i}  
	  \le c_7^{-1} e^{-c_7 \alpha R} \ ,   
\end{align*}
which finishes the proof.
\hfill $\blacksquare$

\section{Large deviation principle for translation distance}
\label{sec.tau}
This section is devoted to the proof of Theorem \ref{theo.tau} that we recall here for reader's convenience.
\begin{theoreme}
\label{theo.tau2}
Let $\Gamma$ be a countable group acting by isometries on a geodesic Gromov-hyperbolic space $X$, $\mu$ a non-elementary probability measure on $\Gamma$ of bounded support. Then the sequence of random variables 
	$ ( \frac{1}{n} \tau(\gamma_n))_{n \in \NN} $
satisfies a large deviation principle with the same rate function as the one given by Theorem \ref{maintheo}.
\end{theoreme}

In order to prove this theorem, we will again make use of the criterion given by Lemma \ref{lemma.Psi.and.I} based on Theorem \ref{thm.LDP.criterion}.\\

As before, we shall distinguish the deviations from above and from below. Let us recall that, by definition, we have for any $g \in \Isom(X)$ and any $x \in X$
	\begin{equation}\label{eq.tau.smaller}
	     \tau(g) \le d(x, g \cdot x) \ .
	 \end{equation}

In particular for any $\alpha > l$  and any $n \in \NN$, we have
	$$ \PP(\tau(\gamma_n) \ge \alpha n) \le \PP( d_n \ge \alpha n) \ .$$
Recall that we denoted $d_n := d(z_0, z_n)$. Then, for any $n \in \NN^*$, we have 
	$$ \frac{-1}{n} \ln \PP(\tau(\gamma_n) \ge \alpha n) \ge \frac{-1}{n} \ln \PP( d_n \ge \alpha n) \ .$$

In view of Lemma \ref{lemma.Psi.and.I}, regarding deviations from above ($\alpha >l$), one is then left to show that for $l_{\max} \neq \alpha > l$, we have 
\begin{equation}
	\label{eq-from above translations}
 	\Psi(\alpha) := \limi{n \to \infty} \ \frac{-1}{n} \ln \left(  \PP (d_n \ge \alpha n ) \right) \ge   \limisup{n \to \infty}  \ \frac{-1}{n} \ln \left(  \PP (\tau(\gamma_n) \ge \alpha n ) \right) \ .
 \end{equation} 

The proof of the above inequality will be carried out in Subsection \ref{subsec comparison above}. It is very close in spirit to the proof of Proposition \ref{propupper} and relies on a Schottky-like argument with insertion trick. \\

For what concerns deviations from below ($\alpha < l$), in view of \eqref{eq.tau.smaller}, for all $n \in \NN$, we have 
	$$ \PP( d_n \le \alpha n)   \le \PP(\tau(\gamma_n) \le \alpha n)  \ .$$
We are then left to prove that for all $ l_{\min} \neq \alpha < l$ we have 
\begin{equation}
	\label{eq-from below translations}
 	\Psi(\alpha) := \limi{n \to \infty} \ \frac{-1}{n} \ln \left(  \PP (d_n \le \alpha n ) \right) \le  \limiinf{n \to \infty}  \ \frac{-1}{n} \ln \left(  \PP (\tau(\gamma_n) \le \alpha n ) \right) \ .
 \end{equation} 

The strategy to prove the above inequality is more involved. We shall detail it in Subsection \ref{subsec comparision below}. The proof is based on a geometric tool whose proof is postponed to Subsection \ref{subsec.ingredients.tau}.

\subsection{Comparison from above}
\label{subsec comparison above} The goal of this subsection is to show that \eqref{eq-from above translations} holds.  It is a consequence of the following

\begin{lemma}
	\label{lemma trans from above}
There exist a constant $c>0$ and an integer $p \in \NN$ such that for any $\alpha \in \RR_+$ and $\varepsilon > 0$, there is an integer $n_0 = n_0(\varepsilon, \alpha)$ such that for any $n \ge n_0$ we have 
	$$ c \ \PP( d_n \ge (\alpha + \varepsilon) n) \le  \PP (\tau(\gamma_{n +p}) \ge \alpha (n +p) ) \ . $$
\end{lemma}

To see that the inequality given by this lemma implies \eqref{eq-from above translations}, one observes that applying logarithm, dividing by $n$ and taking the limsup, we get that for every $ \alpha \in 
(l, l_{\max})$ and $\varepsilon > 0$  
	$$ 	\limisup{n \to \infty} \ \frac{-1}{n} \ln \PP( \tau(\gamma_n) \ge \alpha n)  \le \Psi(\alpha + \varepsilon) \ ,$$
which gives \eqref{eq-from above translations} by continuity of $\Psi$. \\

The proof of Lemma \ref{lemma trans from above} relies on the following geometric ingredient that we will use in combination with Proposition \ref{lemmapingpong}.

\begin{lemma}
\label{lemmaschottky2}
For any $x \in X$ and any Schottky set $S$, there is a constant $L > 0$ with the property that for every $g \in \Isom(X)$, there exists $s \in S$ such that $\tau(s g) \geq d(x,g \cdot x) -L$.
\end{lemma}

\textbf{Proof.}
To simplify the notation, let us denote the basepoint by $z_0=:o$. It clearly suffices to show the claim for $o \in X$. It is well-known (see e.g.\ \cite[Ch.9, Lemma 2.2]{CDP}) that for every $g \in \Isom(X)$, we have 
$$
d(g\cdot o,o) \geq \ell(g) \geq d(g\cdot o,o)- 2(g\cdot o,g^{-1}\cdot o)_o -2\delta
$$
and that $|\ell(.)-\tau(.)|$ is bounded (by $16 \delta$, see \cite[Ch.10, Prop 6.4]{CDP}). Let $S_{\sup}:=\max_{s \in S}d(s \cdot o, o)$. By triangle inequality, $|d(o, gs \cdot o)-d(o,g \cdot o)| \leq S_{\sup}$. Therefore,  we only need to show that there exists a constant $L'>0$ such that for every $g \in \Isom(X)$, there exists $s \in S$ satisfying $(sg\cdot o,g^{-1}s^{-1}\cdot o)_o \leq L'$. Again by triangle inequality and definition of Gromov product, we have 
$|(sg\cdot o,g^{-1}s^{-1}\cdot o)_o-(sg\cdot o,g^{-1}\cdot o)_o| \leq S_{\sup}$ and hence we only need to show that there exists a constant $L''$ such that for any $g \in \Isom(X)$, there exists $s\in S$ satisfying $(sg \cdot o, g^{-1} \cdot o)_o \leq L''$. This follows directly by the defining property of a Schottky set (see Definition \ref{defschottkyset}) applied to $y=g \cdot o$ and $z=g^{-1} \cdot o$. \hfill $\blacksquare$\\

\textbf{Proof of Lemma \ref{lemma trans from above}.}  We shall use an insertion trick similar to the one employed in Section \ref{secupper}. Using Proposition \ref{lemmapingpong}, let $S$ be a Schottky set in the subsemigroup generated by $\mu$ and $p \in \NN$ such that $S \subset \mu^{*p}$.\\

Let $s \in S$, we start by getting a lower bound to the following pivotal quantity. 
	$$ \PP (d(z_p, z_{n+p}) \ge (\alpha + \varepsilon) n \ , \ \gamma_p =s  ) \ .$$ 
On the one hand, since $\gamma_{p} =s$ is independent of $d(z_p, z_{n+p})$ we have
\begin{align*}
	\PP (d(z_p, z_{n+p}) \ge (\alpha + \varepsilon) n,  \ \gamma_p =s ) 
	& = \PP (d(z_p, z_{n+p}) \ge (\alpha + \varepsilon) n) \ \PP(\gamma_p =s ) \\ 
	& = \PP (d_n \ge (\alpha + \varepsilon) n) \ \mu^{*p}(s) \\ 
		& \ge \zeta \cdot \PP (d_n \ge (\alpha + \varepsilon) n) \ , 
\end{align*}
where $\zeta := \infi{s \in S} \ \mu^{*p}(s) > 0$ since we assumed that $S \subset \supp{\mu^{*p}}$ and because $d_n$ and $d(z_p, z_{n+p})$ follows the same law. \\

On the other hand, we have 
	$$ \PP (d(z_p, z_{n+p}) \ge (\alpha + \varepsilon) n, \ \gamma_p =s  ) = 
 \EE ( \mathds{1}_{\{d(z_p, z_{n+p}) \ge (\alpha + \varepsilon) n \ , \ \gamma_p =s \}} ) \ .$$

We get then for any $s \in S$ and any $n \in \NN$ that
 $$\zeta \cdot \PP (d_n \ge (\alpha + \varepsilon) n)  \le \EE ( \mathds{1}_{\{ d(z_p, z_{n+p}) \ge (\alpha + \varepsilon) n \ , \  \gamma_p =s \}} ) \ .$$ 
Averaging over the finite set $S$ yields
 $$  \ \zeta \cdot \PP (d_n \ge (\alpha + \varepsilon) n)  \le \EE \left( \frac{1}{\# S} \somme{s \in S} \ \mathds{1}_{\{d(z_p, z_{n+p}) \ge (\alpha + \varepsilon) n \ , \ \gamma_p =s \}} \right) \ .$$ 
 
Thanks to Lemma \ref{lemmaschottky2}, we have the deterministic upper bound
$$ \somme{s \in S} \ \mathds{1}_{\{d(z_p, z_{n+p}) \ge (\alpha + \varepsilon) n \ , \ \gamma_p =s \}} \le \# S \ \mathds{1}_{\{\tau(\gamma_{n+p}) + L  \ge (\alpha + \varepsilon) n \} } \ , $$ 
and then
\begin{align*}
 	 \zeta \cdot \PP (d_n \ge (\alpha + \varepsilon) n) & \le \PP  (\tau(\gamma_{n+p}) + L  \ge (\alpha + \varepsilon) n ) \\ 
 	 	& \le  \PP  (\tau(\gamma_{n+p}) \ge \alpha (n +p) - L + \varepsilon n - \alpha p ) \\
 	 	& \le   \PP  (\tau(\gamma_{n+p}) \ge \alpha (n +p) ) \ , 
\end{align*}
for $n \ge n_0$ with any $n_0$ such that $- L + \varepsilon n_0 - \alpha p > 0$, concluding the proof. \hfill $\blacksquare$ 

\subsection{Comparison from below} 
\label{subsec comparision below} 
The goal of this subsection is to show that \eqref{eq-from below translations} holds.  We shall actually prove the following stronger bound. 

\begin{lemma}
	\label{lemma trans from below}
For any $\alpha > 0$ and any $\varepsilon > 0$, there is an integer $N \in \NN$ such that for any $n \ge N$, we have
\begin{equation*}	
		\PP( \tau(\gamma_n) \le \alpha n)  \le (n+1) \ \PP(d_n \le (\alpha + \varepsilon) n ) \ . 
	\end{equation*}
\end{lemma}

One easily sees that this lemma implies \eqref{eq-from below translations}. Indeed, taking logarithm, dividing by $n$ and considering the liminf, we get that for $ \alpha \in (l_{\min}, l)$ and all $\varepsilon > 0$, 
	$$ 	\limiinf{n \to \infty} \ \frac{-1}{n} \ln \PP( \tau(\gamma_n) \le \alpha n)  \ge \Psi(\alpha + \varepsilon) \ .$$
This gives \eqref{eq-from below translations} in view of the continuity of $\Psi$. \\

The proof of Lemma \ref{lemma trans from below} relies on the following geometric result. 

\begin{proposition}\label{prop.moving.tau}
For any bounded subset $\mathcal{B}$ of $\Isom(X)$, for each $\beta >0$ there exists $N\geq 1$ so that the following holds. Let $b_1,\dots,b_n\in \mathcal{B}$, for some $n\geq N$, and let $g_i=b_1\dots b_i$ and $r_i=b_{i+1}\dots b_n$. Then for every $r\in [\tau(g_n), d(z_0,g_n \cdot z_0)]$ there exists $i$ so that $|\, d(r_ig_i \cdot z_0, z_0) - r\,|\leq \beta n$.
\end{proposition}

We postpone the proof of the geometric proposition to Subsection \ref{subsec.ingredients.tau}. \\ 

\textbf{Proof of Lemma \ref{lemma trans from below}.} Let $\alpha < l$ and $\varepsilon > 0$. We start with rewriting 
	\begin{equation}	
		\label{eqtranslenghtfrombelow100}
		\begin{split}
		\PP( \tau(\gamma_n) \le \alpha n) = \PP( \tau(\gamma_n) & \le \alpha n \ , \ d_n \ge (\alpha + \varepsilon)n) \\ 
			 	& + \PP( \tau(\gamma_n) \le \alpha n \ ,  \ d_n < (\alpha + \varepsilon)n)  \ .
		\end{split}
	\end{equation}

We shall deal with the above two probabilities separately; for the second one we use the rough upper bound
\begin{equation}
	\label{eqtranslenghtfrombelow1000} 
			\PP( \tau(\gamma_n) \le \alpha n \ ,  \ d_n < (\alpha + \varepsilon)n)  \le \PP( d_n \le (\alpha + \varepsilon)n) \ . 
\end{equation}

For the first one, we rely on the use of Proposition \ref{prop.moving.tau}. \\

We fix $\beta := \varepsilon/2$ and $\mathcal{B}  := \supp{\mu}$ which is bounded by assumption. Let $N$ be large enough so as to get the conclusions of Proposition \ref{prop.moving.tau}. We shall use it with
 \begin{itemize}
 	\item $ b_i := \omega_i$, the successive increments of the walk (and accordingly $g_i =  \gamma_i$) ; 
	 \item $ \tau(\gamma_n) \le r := (\alpha + \varepsilon/2) n \le d_n$ \ .
 \end{itemize}
 
Using that proposition, we deduce that for any $n \ge N$,  we have 
\begin{align*}
	 \{ \tau(\gamma_n)  \le \alpha n \ , \ d_n \ge (\alpha + \varepsilon)n \} & \subset \underset{1 \le i \le n}{\cup} \{ | d(r_i \gamma_i \cdot z_0, z_0) - (\alpha + \varepsilon/2)n | \le \varepsilon n /2 \} \\
	 		& \subset \underset{1 \le i \le n}{\cup} \{  d(r_i \gamma_i \cdot z_0, z_0) \le (\alpha + \varepsilon) n \} \ .
 \end{align*}
 
Note that the random variables $d(z_0, r_i \gamma_i \cdot z_0)$ follows the same law as $d_n$ for every $i \in [1, n]$ since we assumed the increments to be independent and identically distributed. In particular we get 
\begin{align*}
	 \PP (\tau(\gamma_n)  \le \alpha n \ , \ d_n \ge (\alpha + \varepsilon)n ) & \le \somme{1 \le i \le n} \ \PP( d(z_0, r_i \gamma_i \cdot z_0)  \le (\alpha + \varepsilon) n ) \\
	  	& \le n \ \PP(d_n \le (\alpha + \varepsilon) n ) \ . 
\end{align*}

Therefore, looking backward to \eqref{eqtranslenghtfrombelow1000} and \eqref{eqtranslenghtfrombelow100} we get
\begin{equation*}	
		\begin{split}
		\PP( \tau(\gamma_n) \le \alpha n) & \le n \ \PP(d_n \le (\alpha + \varepsilon) n ) + \PP(d_n \le (\alpha + \varepsilon) n ) \\ 
		& \le (n+1) \ \PP(d_n \le (\alpha + \varepsilon) n ) \ ,
		\end{split}
	\end{equation*}
concluding the proof. \hfill $\blacksquare$

\subsection{Proof of Proposition \ref{prop.moving.tau}}
\label{subsec.ingredients.tau} We start with some geometric preliminaries. Some of the results in this subsection might be known to experts, but we provide detailed proofs for completeness. \\ 

Let $(X,d)$ be a Gromov-hyperbolic space and let $\mathcal{B}$ be an arbitrary bounded subset of $\Isom(X)$. To simplify the notation, let $x \in X$ denote the choice of a basepoint.

\begin{remark}
For convenience, in the proofs below we will assume that $d(x,b \cdot x)\leq 1$ for each $b\in \mathcal{B}$. This can be achieved by rescaling $X$, and it is readily seen that all the statements hold for $X$ if and only if they hold for a rescaling of $X$, up to changing the constants.
\end{remark}

The following lemma has a more general version where there is no group action involved, and the sequence of the $g_i \cdot x$ is replaced by any discrete path with bounded jumps. We prefer to state the lemma in the form in which it will get used.

\begin{lemma}\label{lem:cant_avoid_geod}
Let $\mathcal{B}$ be a bounded subset of $\Isom(X)$. For every $\varepsilon>0$ there exist $D_0,N\geq 1$ so that the following holds. Let $b_1,\dots, b_n\in \mathcal{B}$, for some $n\geq N$, and let $g_i=b_1\dots b_i$. Let $\upsilon$ be a subpath of length $\geq \varepsilon n$ of a geodesic from $x$ to $g_n \cdot x$. Then there exists $i$ with $d(g_i \cdot x,\upsilon)\leq D_0$.
\end{lemma}

\textbf{Proof.}
This can be deduced from \cite[Claim 2 within Lemma 2.6]{HS:coarse}, which in our setting says the following. There exist $\varepsilon_0>0$ and $D'>0$ (independent of $x$ and $g_n$) so that, given disjoint balls $B_1,\dots,B_k$ of radius $D\geq D'$ centered on $\upsilon$, any path $\alpha$ from $x$ to $g_nx$ that avoids all $B_i$ satisfies $l(\alpha)\geq k(1+\varepsilon_0)^D$. Choose $D_0\geq D'+1$ so that $(1+\varepsilon_0)^{D_0-1}>3D_0/\varepsilon$. Also, we let $N\geq 6D_0/\varepsilon$, and check that these choices work. In the setting of the statement, suppose by contradiction that we have $d(g_i \cdot x,\upsilon)> D_0$ for all $i$. Then we can find at least $k\geq \varepsilon n/(2D_0)-1\geq \varepsilon n/(3D_0)$ disjoint balls $B_i$ of radius $D_0-1$ centered on $\upsilon$ so that the path $\alpha$ in $X$ obtained concatenating geodesics from $s_ix$ to $s_{i+1}x$ avoids all $B_i$. The length of $\alpha$ is at most $n$, so we obtain:
$$n\geq \frac{\varepsilon n}{3D_0} (1+\varepsilon_0)^{D_0-1}>n,$$
a contradiction. \hfill $\blacksquare$  \\

Let $\delta\geq 1$ be a hyperbolicity constant for $X$. For $g \in \Isom(X)$, define
$$\mathrm{Min}(g) :=\{z\in X: d(z,g \cdot z)\leq \tau(g)+4\delta\}.$$
Also, for $z\in Z$, denote by $\pi^g(z)$ a point in $\mathrm{Min}(g)$ so that $d(z,\pi^g(z))\leq d(z,\mathrm{Min}(g))+1$. (That is, $\pi^g$ is coarsely the closest-point projection to $\mathrm{Min}(g)$.) We can and will assume that $g \cdot \pi^g(z)=\pi^g( g \cdot z)$ holds for all $g$ and $z$. \\

It is known that $\mathrm{Min}(g)$ is quasiconvex (see e.g.\ \cite[Proposition 2.3.3]{delzant-gromov} and \cite[Proposition 2.28]{coulon}), but we will only need the following special case of quasiconvexity, which has a very short proof:

\begin{lemma}\label{lem:Min_qconv}
Let $g \in \Isom(X)$. If $y\in \mathrm{Min}(g)$, then any point on any geodesic from $y$ to $g \cdot y$ is also contained in $\mathrm{Min}(g)$.
\end{lemma}

\textbf{Proof.}
First, observe that given $y\in X$ and a geodesic $[y,g \cdot y]$, any $z\in [y,g \cdot y]$ has
$$d(z,g \cdot z)\leq d(z, g \cdot y)+d(g \cdot y, g \cdot z)=d(z,g \cdot y)+d(y,z)=d(y,g \cdot y).$$
The desired statement easily follows.
\hfill $\blacksquare$  \\

We now show that geodesics from $z$ to $g \cdot z$ pass close to the projection points of the endpoints onto $\mathrm{Min}(g)$. 

\begin{lemma}\label{lem:far_from_min}
There exists $D_1\geq 0$ so that the following holds. For every $g \in \Isom(X)$ and $z\in X$, we have
that any geodesic $\upsilon$ from $z$ to $g \cdot z$ passes $D_1$-close to $\pi^g(z)$ and $\pi^g(g \cdot z)$. Moreover, we have
$$d(z,g \cdot z)\geq 2d(z,\pi^g(z))+\tau(g)-D_1.$$
\end{lemma}

\textbf{Proof.}
Consider any geodesic $\upsilon$ from $z$ to $g \cdot z$. We will show that $\upsilon$ passes $(4\delta+2)$-close to $\pi^g(z)$, the argument for $\pi^g(g \cdot z)$ being similar. We will use $2\delta$-thinness of a quadrangle with vertices $z,\pi^g(z),\pi^g(g \cdot z),g \cdot z$. \\

Suppose by contradiction that $\upsilon$ does not pass $(4\delta+2)$-close to $\pi^g(z)$. Consider the point $z'$ on a geodesic from $z$ to $\pi^g(z)$ at distance $2\delta+2$ from $\pi^g(z)$. We observe that $z'$ cannot be $2\delta$-close to any geodesic $[\pi^g(z),\pi^g(g \cdot z)]$, for otherwise there would be a point $q$ on said geodesic, whence on $\mathrm{Min}(g)$ by Lemma \ref{lem:Min_qconv}, which satisfies $d(z,q)<d(z,\pi^g(z))-1$, contradicting the defining property of $\pi^g$. \\

Also, $z'$ cannot be $2\delta$-close to $\upsilon$ by hypothesis, so $z'$ is $2\delta$-close to $g([z,\pi^g(z)])$. But then it must be $4\delta$-close to the point on that geodesic at distance $4\delta+2$ from $\pi^g(g \cdot z)$, this point being $g \cdot z'$. We just showed $d(z',g \cdot z')\leq 4\delta$, which implies $z'\in \mathrm{Min}(g)$. But $d(z,z')<d(z,\pi^g(z))-1$, contradicting the defining property of $\pi^g(z)$. \\

Now, the fact that $\upsilon$ passes $(4\delta+2)$-close to $\pi^g(z)$ and $\pi^g(g \cdot z)$ implies the following inequality:
$$d(z,g \cdot z)\geq d(z,\pi^g(z))+d(\pi^g(z),\pi^g(g \cdot z))+d(\pi^g(g \cdot z),g \cdot z)-4(4\delta+2).$$
The first and third terms on the right-hand side are both equal to $d(z,\pi^g(z))$, while the second term is at least $\tau(g)$. Therefore, we can conclude by setting $D_1=4(4\delta+2)$.
\hfill $\blacksquare$  \\

We are now ready to prove Proposition \ref{prop.moving.tau}.\\[10pt] 
\textbf{Proof} (of Proposition \ref{prop.moving.tau}).
First, notice that for $r\geq d(g_n \cdot x, x)-  \varepsilon n$ we can just choose $i=n$, so in the arguments below we assume $r\leq d(g_n \cdot x, x) -\varepsilon n$. \\

Set $d :=(r-\tau(g_n))/2$, and assume that $n$ is larger than the $N$ from Lemma \ref{lem:cant_avoid_geod} with $\varepsilon/4$ replacing $\varepsilon$. We impose further constraints on $n$ later. 

\begin{lemma} 
	\label{lemme proof prop translation} If $n$ is sufficiently large, then we can find a subgeodesic $\upsilon$ of length $\varepsilon n/4$ of a geodesic $\upsilon'$ from $x$ to $g_nx$  so that any $p\in\upsilon$ has 
\begin{itemize}
    \item $d(p,q)\leq D_1+\delta$ for some $q$ on a geodesic from $x$ to $\pi^{g_n}(x)$,
    \item $|d(p,\pi^{g_n}(x))-d|\leq \varepsilon n/3$.
\end{itemize} 
\end{lemma}

\textbf{Proof.}
We let $p_0$ be the point along $\upsilon'$ so that
$$d(x,p_0)=d(x,\pi^{g_n}(x))-d-\varepsilon n/4=\hat d,$$
and we let $\upsilon$ be the subgeodesic of $\upsilon'$ of length $\varepsilon n/4$ with starting point $p_0$. We now check that, for $n$ large enough, this is all well-defined, and that $\upsilon$ has the required property. Let us make the preliminary observation that
\begin{align*}
	 d(x,g_n \cdot x) & \leq d(x,\pi^{g_n}(x))+d(\pi^{g_n}(x),\pi^{g_n}(g_nx))+ d(\pi^{g_n}(g_nx),g_nx) \\
	 & \leq 2d(x,\pi^{g_n}(x))+\tau(g_n)+4\delta.
\end{align*}
Observe now that we have
$$d\leq (d(x,g_n \cdot x)-\varepsilon n-\tau(g_n))/2\leq d(x,\pi^{g_n}(x))+4\delta-\varepsilon n/2,$$
implying $\hat d\geq d-4\delta +\varepsilon n/4$, which is a positive quantity if $n$ is sufficiently large.

Also, again by Lemma \ref{lem:far_from_min}, there exists $p'\in \upsilon'$ so that $d(p',\pi^{g_n}(x))\leq D_1$; denote by $\upsilon''$ the initial subgeodesic of $\upsilon'$ with terminal point $p'$.

Notice that for $n$ large enough we have
$$d(x,p')\geq d(x,\pi^{g_n}(x))-D_1= d+\hat d +\varepsilon n/4-D_1 \geq \hat d +\varepsilon n/4.$$
The inequalities we just showed imply that $p_0$ and $\upsilon$ are well-defined and, furthermore, that $\upsilon$ is a subgeodesic of $\upsilon''$. \\

Considering a triangle with vertices $x,p', \pi^{g_n}(x)$ and containing $\upsilon''$, we see that any point on $\upsilon''$, whence any point on $\upsilon$, is $(D_1+\delta)$-close to a point on a geodesic from $x$ to $\pi^{g_n}(x)$. In particular, for any $p\in \upsilon$ we have
$$|d(p,\pi^{g_n}(x))+d(x,p) -d(x,\pi^{g_n}(x)|\leq 2D_1+2\delta,$$
and hence
\begin{align*}
|d(p,\pi^{g_n}(x)) - d| & \leq |d(x,\pi^{g_n}(x))-d(x,p)- d| +2D_1+2\delta \\ 
  & \leq |d(x,\pi^{g_n}(x))-d(x,p_0)- d- d(p_0,p)| + 2D_1+2\delta \\
  & \le |\varepsilon n/4 - d(p_0,p)| + 2D_1+2\delta \\
  & \leq\varepsilon n/4 +2D_1+2\delta.
\end{align*}
Provided that $n$ is large enough, this concludes the proof of the claim.
\hfill $\blacksquare$  \\

By Lemma \ref{lem:cant_avoid_geod}, there exists $i$ so that we have $d(g_i,p)\leq D_0$ for some $p\in\upsilon$, and by the claim we have $d(p,q)\leq D_1+\delta$, whence $d(g_i \cdot x,q)\leq D_0+D_1+\delta$, for some $q$ on a geodesic from $x$ to $\pi^{g_n}(x)$. Notice that we can assume that $\pi^{g_n}(q)=\pi^{g_n}(x)$ since
$$d(q,\pi^{g_n}(x))=d(z,\pi^{g_n}(x))-d(q,x)\leq d(x,\mathrm{Min}(g_n))+1-d(q,x)\leq d(q,\mathrm{Min}(g_n))+1.$$
In particular, in view of Lemma \ref{lem:far_from_min}, we have
$$d(q,g_n \cdot q)\geq 2d(q,\pi^{g_n}(x))+\tau(g_n)-D_1.$$
We can now compute
\begin{align*}
d(x,r_ig_i \cdot x) & =d(g_i \cdot x,g_n g_i \cdot x)  \\
	& \geq d(q,g_n \cdot q)- 2(D_0+D_1+\delta) \\
	& \geq 2d(q,\pi^{g_n}(x))+\tau(g_n)- 3(D_0+D_1+\delta) \ .
\end{align*}
Hence,
$$d(x,r_ig_i \cdot x) \geq  2d +\tau(g_n) - 2\varepsilon n/3- 5(D_0+D_1+\delta)=  r - 2\varepsilon n/3 -5(D_0+D_1+\delta).$$
For $n$ sufficiently large, this last quantity is $\geq r-\varepsilon n$.

On the other hand, we also have
\begin{align*}
d(x,r_ig_i \cdot x) & = d(g_i \cdot x,g_n g_i \cdot x) \\ 
		& \leq d(p,\pi^{g_n}(x))+d(\pi^{g_n}(x),\pi^{g_n}(g_n \cdot x)) +d(g_n \cdot p,\pi^{g_n}(g_n \cdot x))+2D_0 \\  		& \leq  2d+2\varepsilon n/3 +2D_1+\tau(g_n)+4\delta+2D_0 \\
		& \leq r+2\varepsilon n/3+2D_0+2D_1+4\delta.
\end{align*}

For $n$ sufficiently large, this last quantity is $\leq r+\varepsilon n$, concluding the proof.
\hfill $\blacksquare$  

\section{Support of the rate function}\label{sec.support}
We start by recording a characterization of non-arithmeticity of a non-elementary set in \S \ref{subsec.arithmeticity} which will then be used to prove Theorem \ref{thm.support} in \S \ref{subsec.support}. In \S \ref{subsec.examples}, we discuss the examples mentioned in Remark \ref{rk.examples.intro} and finally in \S \ref{subsec.proof.of.joint.spectrum}, we prove Proposition \ref{prop.joint.spectrum}.

\subsection{A characterization of non-arithmeticity of a non-elementary set}\label{subsec.arithmeticity}

Recall from \eqref{eq.def.ell.gamma} that for an element $\gamma \in \Isom(X)$, 
$\ell(\gamma)$ denotes the asymptotic translation length given by $\lim_{n \to \infty} d(x, \gamma^n \cdot x)/n$ for any $x \in X$. Furthermore, recall that a subset $\mathcal{B} \subseteq \Isom(X)$ is called non-arithmetic if there exist $n \in \mathbb{N}$ and $g_1,g_2\in \mathcal{B}^n$ such that $\ell(g_1)\neq \ell(g_2)$. In the proof of Theorem \ref{thm.support}, we will use the following characterization of non-arithmeticity of a non-elementary set $\mathcal{B}$ in terms of the  asymptotic joint displacements $\ell(\mathcal{B})$ and $\ell_{\mathrm{sub}}(\mathcal{B})$ defined in \eqref{eq.def.l.lsub}:

\begin{proposition}\label{prop.arithmeticity}
Given a non-elementary subset $\mathcal{B}$ of the isometry group of a Gromov-hyperbolic space $(X,d)$, the set $\mathcal{B}$ is non-arithmetic if and only if $\ell_{\sub}(\mathcal{B}) \neq \ell(\mathcal{B})$.
\end{proposition}

For the proof, we will need the geometric Berger--Wang equality proved recently in \cite{oregon.reyes:properties,breuillard-fujiwara}. We will provide a brief proof of this equality for the non-elementary case using the tools we developed. To state it, for a subset $\mathcal{B}$ of $\Isom(X)$, we consider the following numerical invariant 
	$$\ell_{\infty}(\mathcal{B}):=\limsup_{k \in \mathbb{N}^\ast} \sup_{g \in \mathcal{B}^k} \frac{1}{k} \ell(g) \ . $$ 
Clearly, $\ell_\infty(\mathcal{B})$ is a conjugacy invariant and we have $\ell_{\infty}(\mathcal{B}) \leq \ell(\mathcal{B})$ (see also \cite[Lemma 1.1]{breuillard-fujiwara}). Using a Schottky-like argument, one gets
\begin{lemma}[Geometric Berger--Wang equality, \cite{oregon.reyes:properties,breuillard-fujiwara}]\label{prop.berger-wang}
For a non-elementary subset $\mathcal{B}$ of $\Isom(X,d)$, we have $\ell_\infty(\mathcal{B})=\ell(\mathcal{B})$.
\end{lemma}

\textbf{Proof.} In view of the definitions, we can suppose $\mathcal{B}$ to be countable. Assume for a contradiction that we have $\ell_\infty(\mathcal{B})<\ell(\mathcal{B})$. It follows that there exists $\delta>0$ such that for every $n \geq 1$, there exists $g_n \in \mathcal{B}^n$ with $\frac{1}{n}d(z_0,g_n \cdot z_0) \geq \ell_\infty(\mathcal{B})+\delta$. By Proposition \ref{lemmapingpong}, Lemma \ref{lemmaschottky2} and the fact that $|\tau(.)-\ell(.)|$ is uniformly bounded, we deduce that there exist $p \in \mathbb{N}$ and a Schottky set $S \subset \mathcal{B}^p$ such that for every $n$ large enough, there exists $s_n \in S$ such that $\frac{1}{n}\ell(s_ng_n) \geq \ell_\infty(\mathcal{B})+\delta/2$. Since $s_n g_n \in \mathcal{B}^{p+n}$, we get that for every $n \in \mathbb{N}$ large enough, we have
$$
\frac{\ell(s_ng_n)}{p+n} \geq \ell_\infty(\mathcal{B})+\frac{\delta}{4}.
$$
This clearly yields a contradiction in view of the definition of $\ell_{\infty}(\mathcal{B})$ using the fact that $\ell(g^n)=n\ell(g)$ for every $n \in \mathbb{N}$ and $g \in \Isom(X)$.
\hfill $\blacksquare$\\  

\textbf{Proof of Proposition \ref{prop.arithmeticity}}
Suppose $\mathcal{B}$ is non-arithmetic, i.e.\  there exist $n \in \mathbb{N}$ and $g_1,g_2\in \mathcal{B}^n$ such that $\ell(g_2)> \ell(g_1)$. From the facts that for every $g \in \Isom(M)$ and $m \in \mathbb{N}$, we have $\ell(g^m)=m \ell(g)$ and $d(z_0,g \cdot z_0) \geq \ell(g)$, we deduce
$$
\ell(\mathcal{B}) \geq \frac{\ell(g_2)}{n}>\frac{\ell(g_1)}{n} \geq \ell_{\mathrm{sub}}(\mathcal{B}),
$$
proving the first implication.
For the other implication, suppose that $\ell(\mathcal{B})>\ell_{\mathrm{sub}}(\mathcal{B})$. By Lemma \ref{prop.berger-wang}, we have $\ell_\infty(\mathcal{B})>\ell_{\mathrm{sub}}(\mathcal{B})$. This clearly implies that $\mathcal{B}$ is non-arithmetic. \hfill $\blacksquare$


\subsection{Proof of Theorem \ref{thm.support}}\label{subsec.support}
Let $\mu$ be a non-elementary probability measure on $\Isom(X)$ and $I:[0,\infty) \to [0,\infty]$ the rate function given by Theorem \ref{thm.ldp.general} (equivalently, by Theorem \ref{maintheo} if $\mu$ has a finite exponential moment). Denote by $\mathcal{B}$ the support of $\mu$.\\

We first study the lower bound of $D_I$; this part does not require any additional ingredients.  Now let $\alpha> \ell_{\mathrm{sub}}(\mathcal{B})$. Then, there exists $n_0 \in \mathbb{N}$ and $g \in \mathcal{B}^{n_0}$ such that $\frac{1}{n_0} d(g \cdot z_0,z_0)<\alpha$. By triangle inequality, we also have $\frac{1}{kn_0} d(g^k \cdot z_0,z_0) <\alpha$ for every positive $k \in \mathbb{N}$. We deduce 
$$\limsup_{n \to \infty}\frac{1}{n}\log \mathbb{P}(\frac{1}{n} d(z_n,z_0) \leq \alpha) \geq \limsup_{k \to \infty} \frac{1}{kn_0} \log \mathbb{P}(\gamma_{k n_0}=g^k) \geq \frac{1}{n_0} \log \mu^{\ast n_0}(g)>-\infty,$$
where in the second inequality we used the fact that $\mu^{n_0 k}(g^k) \geq \mu^{n_0}(g)$ which is an immediate consequence of the I.I.D. increments assumption. Since $\alpha>\ell_{\sub}(\mathcal{B})$ is arbitrary, it follows that for every $\varepsilon>0$, we have 
\begin{equation}\label{eq.lower.support}
D_I \cap (-\infty,\ell_{\sub}(\mathcal{B})+\varepsilon) \neq \emptyset.  
\end{equation}

To study the upper bound of $D_I$, we will make use the existence of a Schottky set (Proposition \ref{lemmapingpong}) and Lemma \ref{lemmaschottky2}. Let $\beta<\ell(\mathcal{B})$ be given. Then for every $\delta<\ell(\mathcal{B})-\beta$, for every $n \in \mathbb{N}$ large enough, there exists $g_n\in \mathcal{B}^{n}$ such that $\frac{1}{n}d(g_n\cdot z_0,z_0)>\beta+\delta$. By Proposition \ref{lemmapingpong}, there exists $p \in \mathbb{N}$ such that $\mathcal{B}^p$ contains a Schottky set $S$. It then follows by Lemma \ref{lemmaschottky2} and the fact that $|\tau(.)-\ell(.)|$ is uniformly bounded, that for every $n \in \mathbb{N}$ large enough, there exists $s_n \in \mathcal{B}^p$ such that $\frac{1}{n+p} \ell(s_n g_n) \geq \beta + \delta/2$. Fix a large enough $n_0 \in \mathbb{N}$ such that the latter inequality holds. Now since for every $g \in \Isom(X)$ and $k \in \mathbb{N}$, we have $\ell(g^k)=k \ell(g)$, we have for every $k \in \mathbb{N}$, 
$$
\frac{1}{k(n_0+p)}d((s_{n_0}g_{n_0})^k \cdot z_0,z_0) \geq \frac{1}{k(n_0+p)}\ell( (s_{n_0}g_{n_0})^k)=\frac{1}{n_0+p} \ell(s_{n_0}g_{n_0}) \geq \beta + \frac{\delta}{2}.
$$
Therefore we deduce
\begin{equation*}
\begin{aligned}
\limsup_{n \to \infty} \frac{1}{n} \log \mathbb{P}(\frac{1}{n} d(z_n,z_0) \geq \beta) & \geq \limsup_{k \to \infty} \frac{1}{k(n_0+p)} \mathbb{P}(\gamma_{k(n_0+p)}=(s_{n_0}g_{n_0})^k)\\  &\geq \frac{1}{n_0+p} \log \mu^{\ast (n_0+p)}(s_{n_0}g_{n_0})>-\infty.
\end{aligned}
\end{equation*}

Since $\beta<\ell(\mathcal{B})$ is arbitrary, it follows that for every $\varepsilon>0$, we have $D_I \cap (\ell(\mathcal{B}-\varepsilon),\infty) \neq \emptyset$. Since $D_I$ is an interval, together with \eqref{eq.lower.support}, this implies that $D_I \supseteq (\ell_{\sub}(\mathcal{B}),\ell(B))$. \\

On the other hand, unfolding the definitions, it is plain that we have $l_{\min} \geq \ell_{\mathrm{sub}}(\mathcal{B})$ and $l_{\max} \leq \ell(\mathcal{B})$. This proves that $l_{\min}=\ell_{\mathrm{sub}}(\mathcal{B})$ and $l_{\max}=\ell(\mathcal{B})$. The fact that $D_I$ has non-empty interior if and only if $\mu$ is non-arithmetic now follows from Proposition \ref{prop.arithmeticity}.\\

To prove the last statement, suppose that $\mu$ is finitely supported.
If $\ell_{\sub}(\mathcal{B})=\ell(\mathcal{B})$, then it is not hard to see that $D_I=\{\ell(\mathcal{B})\}$. Therefore we suppose that $\ell_{\sub}(\mathcal{B})<\ell(\mathcal{B})$, in other words, by Proposition \ref{prop.arithmeticity}, $\mu$ is non-arithmetic. We only need to show that on $(\ell_{\sub}(\mathcal{B}),\ell(\mathcal{B}))$, the rate function $I$ is bounded above by $-\min_{g\in \mathcal{B}}\ln \mu(g)<\infty$; lower semi-continuity of $I$ then entails that $I$ is bounded above by the same quantity on $[\ell_{\sub}(\mathcal{B}),\ell(\mathcal{B})]$ proving the claim. So let $\alpha \in (\ell_{\sub}(\mathcal{B}),\ell(\mathcal{B}))$. Then we have $\alpha \in D_{I}$ and by Theorem \ref{maintheo}, for any $\delta>0$, for every small enough $r>0$ and large enough $n \in \mathbb{N}$, we have
\begin{equation}\label{eq.supp2}
0<e^{-n(I(\alpha)+\delta)}\leq \mathbb{P}(\alpha-r \leq \frac{1}{n}d_n \leq \alpha +r) \leq e^{-n(I(\alpha)-\delta)}.
\end{equation}
It follows that for every such $n \in \mathbb{N}$, there exists $g_n \in \mathcal{B}^n$ with $\frac{1}{n}d(z_0,g_n \cdot z_0) \in [\alpha-r,\alpha+r]$. By the I.I.D.\ property of random walk increments, writing $g_n$ as a product $h_1\ldots h_n$ with $h_i \in \mathcal{B}$, it follows that $\mathbb{P}(\gamma_n=g_n) \geq (\min_{h \in \mathcal{B}}\mu(h))^n$. Plugging this in \eqref{eq.supp2}, since $\delta>0$ is arbitrary, we deduce that $I(\alpha)\leq -\min_{h \in \mathcal{B}} \ln \mu(h)$, as required. \hfill $\blacksquare$ \\

\subsection{Examples of rate function exploding on the boundary}\label{subsec.examples}
 
We now construct some examples illustrating in the setting of Theorem \ref{thm.support} that when the support of the probability measure is not finite, the rate function of LDP can explode on $l_{\min}$ or $l_{\max}$ or both. In fact, by considering the action of $\SL_2(\mathbb{R})$ on the Poincar\'{e} disc $\mathbb{D}$ and using the relation $\frac{1}{2}\ln\|g\|=d(g \cdot o, o)$ where $g \in G$, $o$ denotes the origin in $\mathbb{D}$ and $\|.\|$ is the operator norm induced by the Euclidean norm on $\mathbb{R}^2$, \cite[Example 5.5]{sert.LDP} already provides an example of a rate function that explodes on $l_{\max}$. Below, we shall give more examples where rate function explodes on any subset of $\{l_{\min},l_{\max}\}$.

\begin{exxample}
Consider $G=\SL_2(\mathbb{R})$ acting isometrically on the Poincar\'{e} disc $\mathbb{D}$ endowed with the 
usual hyperbolic metric $d$. Let $c_a$ and $c_b$ be two geodesics in $\mathbb{D}$ that are of distance $d=1$ to the origin and denoting their endpoints on $\partial \mathbb{D}$, respectively, by $\{x_a^+,x_a^-\}$ and $\{x_b^+,x_b^-\}$, suppose that these are ordered as $(x_b^-,x_b^+,x_a^+,x_a^-)$. For $k \geq 1$, $t \in \{a,b\}$, $t_k$ be hyperbolic elements of $G$ with 
translation axis $c_t$, attracting/repelling fixed points $x_t^+$ and $x_t^-$ and translation distance $10-\frac{1}{k}$ for $t=b$ and $\frac{1}{k}$ for $t=a$. Let $S=\{a_k,b_k \, |\, k \geq 1\}$ and let $Y_i$'s be 
the coordinate functions on $S^{\mathbb{N}}$. It is easy to see that the subsemigroup $\Gamma$ generated by $S$ 
consists of hyperbolic elements whose translation axes is contained in the connected region bounded by $c_a$ and $c_b$. It follows, e.g.\ by \cite[Lemma 6.3]{breuillard-sert}, that denoting by $D>0$ twice the 
distance between $c_a$ and $c_b$, for any $g,h \in \Gamma$, we have
\begin{equation}\label{eq.ex.1}
\tau(g)+\tau(h) \leq \tau(gh) \leq \tau(g)+\tau(h)+D.    
\end{equation}
Now for $t \in \{a,b\}$, $n \geq k \geq 1$ and $(s_1,\ldots,s_n) \in S^n$, denote by $\hat{t}_{n,k}$ the number of $t_{i}$'s in $(s_1,\ldots,s_n)$ with $i \geq k$. It is readily seen by \eqref{eq.ex.1} that we have the following inclusion of events for every $n \geq 1$ and $1 \leq  k \leq n$:
\begin{equation}\label{eq.ex.inc1}
    \left\{\frac{1}{n}\tau(\omega_1\ldots \omega_n)>10-\frac{1}{2k}\right\} \subset \left\{\hat{b}_{n,k} \geq \frac{n}{2}\right\}
\end{equation}
and
\begin{equation}\label{eq.ex.inc2}
    \left\{\frac{1}{n}\tau(\omega_1\ldots \omega_n)<\frac{1}{2k}\right\} \subset \left\{\hat{a}_{n,k} \geq \frac{n}{2}\right\}.
\end{equation}
Notice also that by elementary plane hyperbolic geometry, for every $g \in \Gamma$, we have 
\begin{equation}\label{eq.ex.2}
0 \leq  d(g \cdot o,o)-\tau(\gamma) \leq D +4.    
\end{equation}
For any probability measure $\mu$ on $S$, the random variables $d_n=d(z_n,z_0)$ satisfies a LDP with some rate function $I$; this follows from Theorem \ref{maintheo} if the support of $\mu$ contains $a_i$'s and $b_i$'s (so that $\mu$ is non-elementary) and from classical theorem of Cram\'{e}r if the support contains only $a_i$'s or $b_i$'s. 
The inequality \eqref{eq.ex.2} entails by Theorem \ref{thm.LDP.criterion} (or using \cite[Theorem 4.2.13]{dembo-zeitouni}) that $\frac{1}{n} \tau(\gamma_n)$ satisfies a LDP with rate function $I$ too. Now let $\alpha_k>0$ be such that $\sum_{k \geq 1} \alpha_k=\frac{1}{2}$ and consider $\mu_1=\frac{1}{2}\delta_{a_1}+\sum_{k \geq 1} \alpha_k \delta_{b_k}$ supported on $S_1$ and $\mu_2=\frac{1}{2}\delta_{b_1}+\sum_{k \geq 1} \alpha_k \delta_{a_k}$ supported on $S_2$ and $\mu=\frac{1}{2}(\mu_1+\mu_2)$ supported on $S$. Denote by $I_1$, $I_2$ and $I$ the rate functions of the LDP of $\frac{1}{n}d(z_n,z_0)$ when the driving probability measure is, respectively, $\mu_1, \mu_2$ and $\mu$. Using Theorem \ref{thm.support} and Stirling's formula, it is not hard to deduce from \eqref{eq.ex.inc1} that $\ell(S_1) \notin D_{I_1}$, and from \eqref{eq.ex.inc2} that $\ell_{\sub}(S_2) \notin D_{I_2}$, and finally, that we have $D_I=(\ell_{\sub}(S),\ell(S))$. Moreover, using \eqref{eq.ex.1}, one sees that  $D_{I_1}=[\ell_{\sub}(S_1),\ell(S_1))$ and $D_{I_2}=(\ell_{\sub}(S_2),\ell(S_2)]$.
\end{exxample}

\subsection{Proof of Proposition \ref{prop.joint.spectrum}}\label{subsec.proof.of.joint.spectrum}

One can use Lemma \ref{lemmaschottky2} and Proposition \ref{prop.moving.tau} to give a direct proof of Proposition \ref{prop.joint.spectrum}. Here, we give a short proof based on our large deviations results.\\

It follows from the definitions that for every $\varepsilon>0$, there exists $N \in \mathbb{N}$ such that for every $n \geq N$ and $g\in \mathcal{B}^n$, $\ell_{\sub}(\mathcal{B}) \leq \frac{1}{n} d(z_0,g \cdot z_0) \leq \ell(\mathcal{B})+\varepsilon$. This already implies the statement if $\ell_{\sub}(\mathcal{B})=\ell(\mathcal{B})$, so let $\ell_{\sub}(\mathcal{B})<\ell(\mathcal{B})$.

Then, for every $\varepsilon \in (0,(\ell(\mathcal{B})-\ell_{\sub}(\mathcal{B}))/2) $, the set $\mathcal{B}$ contains a finite subset $\mathcal{B}'$ such that $\ell(\mathcal{B}')\geq \ell(\mathcal{B})-\varepsilon$ and $\ell_{\sub}(\mathcal{B}')\leq \ell_{\sub}(\mathcal{B})-\varepsilon$. This follows from the definitions of $\ell(\mathcal{B})$ and $\ell_{\sub}(\mathcal{B})$. We can therefore suppose that $\mathcal{B}$ is finite. Now endow $\mathcal{B}$ with the uniform probability measure and consider the corresponding random walk on $\Isom(X)$. Given an interval $J$ of non-empty interior in $[\ell_{\sub}(\mathcal{B}),\ell(\mathcal{B})]$, by Theorem \ref{maintheo} and  Theorem \ref{thm.support}, we have $-\liminf_{n \to \infty} \frac{1}{n} \ln \mathbb{P}(\frac{1}{n}d(\gamma_n \cdot z_0,z_0) \in J)<\infty$. This says, in particular, that for every $n \in \mathbb{N}$ large enough, we have $J \cap \frac{1}{n}d(\mathcal{B}^n \cdot z_0,z_0) \neq \emptyset$. Together with the first paragraph above, this shows the Hausdorff convergence of $\frac{1}{n}d(\mathcal{B}^n \cdot z_0,z_0)$. The convergence of the sets $\frac{1}{n}\tau(\mathcal{B}^n)$ is deduced similarly using Theorem \ref{theo.tau}.
\hfill $\blacksquare$

\appendix

\section{Existence of Schottky sets}
\label{appendixB}

We prove Proposition \ref{lemmapingpong} that we recall here for the reader's convenience.

\begin{proposition}[Existence of Schottky sets]
\label{propschottkyappendix} Let $\Gamma$ be a countable group acting by isometries on a geodesic Gromov-hyperbolic space $X$, $z_0 \in X$ and $\mu$ a non-elementary probability measure on $\Gamma$. Then there is $p \in \NN$ such that $\supp{\mu^{*p}}$ contains a Schottky set. 
\end{proposition}

\textbf{Proof.} We first reduce the proof to a purely geometric statement. Since we assumed that $\supp{\mu}$ generates a non-elementary subsemigroup there are two independent loxodromic elements $\gamma_1, \gamma_2 \in \Gamma$ and $p_1, p_2 \in \NN$ such that 
	$$ \mu^{*p_i}(\gamma_i) > 0 $$
for $ i \in \{1, 2\}$. In particular we have
	$$ \left\{
	 	\begin{array}{l}	
			\mu^{*(p_1 p_2)}(\gamma_1^{p_2}) > 0 \\
			\mu^{*(p_1 p_2)}(\gamma_2^{p_1}) > 0 \ .
		\end{array} \right.
	$$

Because $\gamma_1^{p_2}$ (resp. $\gamma_2^{p_1}$) has the same fixed points as $\gamma_1$ (resp. $\gamma_2$), the pair $(\gamma_1^{p_2}, \gamma_2^{p_1})$  is still a pair of two independent loxodromic isometries. Therefore, up to taking some power of $\mu$ one can suppose that $\supp{\mu}$ contains two independent loxodromic elements. \\

For any pair $\gamma_1, \gamma_2 \in \Gamma$, let $S_k(\gamma_1, \gamma_2) \subset \Gamma$ be the set of all elements of $\Gamma$ which can be written as a product of exactly $k$ elements in $\{\gamma_1, \gamma_2\}$. Note that $S_k(\gamma_1, \gamma_2) $ is contained in the support of $\mu^{*k}$. Proposition \ref{propschottkyappendix} is an immediate consequence of the following 

\begin{proposition}
\label{propappendice}
Let $\gamma_1, \gamma_2$ two independent loxodromic isometries. Then there is $k \in \NN$ such that $S_k(\gamma_1, \gamma_2)$ contains a Schottky set as in Definition \ref{defschottkyset}.
\end{proposition}

\textbf{Proof.} For any points $x,y \in X$ and any $C >0$, we define 
	 $$ \mathcal{O}_C(x,y) := \{z \in X \ , \ (y,z)_x \ge d(x,y) - C\} \ ,$$
which we call the \textbf{$C$-shadow} of $y$ seen from $x$. Note that one can define it equivalently as
 $$ \mathcal{O}_C(x,y) := \{z \in X \ , \ (x,z)_y \le C \} \ ,$$
which is to say, when $X$ is geodesic and up to a constant depending on $\delta$, the set of all points $z$ such that any geodesic from $z$ to $x$ passes through the ball $B(y,C)$. \\

An easy consequence of the Morse lemma is the following. 

\begin{lemma}\label{lem:complement_inclusion}
For any $\lambda, C > 0$ there is a constant $K_0> 0$ such that for any $(\lambda,C)$-quasi-geodesic $(x_n)_{n \in \ZZ}$, any $m \le n \le p$, and any $K\geq K_0$ we have
\begin{equation}
	\label{eqlemballcomplement}
	\mathcal{O}_{K}(x_m,x_n)^c \subset  \mathcal{O}_{K}(x_p,x_n) \ .\end{equation}
\end{lemma}

\textbf{Proof.} 
The Morse Lemma gives some $L$ so that for all $m \le n \le p$ any geodesic from $x_m$ to $x_p$ passes $L$-close to $x_n$. In particular, this implies $(x_m,x_p)_{x_n}\leq L$. 

Fix any $K$ larger than $L+\delta$, and consider any $x\in \mathcal{O}_{K}(x_m,x_n)^c$, where $m \le n \le p$. By definition, we have $(x_m,x)_{x_n} > K$. Keeping $(x_m,x_p)_{x_n}\leq L$ into account, hyperbolicity yields
$$L\geq (x_m,x_p)_{x_n}\geq \min\{(x_m,x)_{x_n},(x,x_p)_{x_n}\}-\delta.$$
This forces $(x,x_p)_{x_n}\leq L+\delta\leq K$, which gives $x\in \mathcal{O}_{K}(x_p,x_n)$, as required.
\hfill $\blacksquare$ \\







\begin{lemma}\label{lem:inclusionshadow}
For any $\lambda, C > 0$ there is a constant $K_0> 0$ such that for any $K\geq K_0$ there exists $N$ with the following property. For any $(\lambda,C)$-quasi-geodesic $(x_n)_{n \in \ZZ}$, any $m \le n \le p$ with $n-m\geq N$ we have
\begin{equation}
	\label{eqinclusionshadow}
	\mathcal{O}_K(x_n,x_m) \subset \mathcal{O}_{K+\delta}(x_p,x_m) \ .\end{equation}
\end{lemma}

\textbf{Proof.}
As in the proof of Lemma \ref{lem:complement_inclusion}, let $L$ (depending on $\lambda$ and $C$) be so that for all $m \le n \le p$ we have $(x_m,x_p)_{x_n}\leq L$, so that we also have $(x_n,x_p)_{x_m}\geq d(x_m,x_n)-L$. Fix any $K\geq L+\delta$. If $n-m$ is sufficiently large (depending on $K$), then we have $(x_n,x_p)_{x_m}> K+\delta$. By definition, if $x\in \mathcal{O}_K(x_n,x_m)$, then $(x_n,x)_{x_m}\leq K$. Using hyperbolicity, we get
$$\min\{(x_p,x)_{x_m},(x_n,x_p)_{x_m}\}\leq (x,x_n)_{x_m}+\delta\leq K+\delta.$$
This forces $(x_p,x)_{x_m}\leq K+\delta$, that is, $x\in\mathcal{O}_{K+\delta}(x_p,x_m)$, as required. \hfill $\blacksquare$ \\


We will also need the next  lemma to set up the ping-pong table.  

\begin{lemma}
\label{lemloxo}
 Let $x_0 \in X$ and $\gamma$ be a loxodromic isometry of $X$. Then, there exists $K_1 > 0$ such that for any sufficiently large $n > 0$ we have
 	$$ \gamma^{2n} \left( \mathcal{O}_{K_1}( x_0, \gamma^{-n} \cdot x_0)^c \right) \subset \mathcal{O}_{K_1}(x_0, \gamma^{n} \cdot x_0) \ .$$ 
\end{lemma}

Note that, since $\gamma$ is an invertible isometry, we also have  
 	$$ \gamma^{-2n} \left( \mathcal{O}_{K_1}( x_0, \gamma^{n} \cdot x_0)^c \right) \subset \mathcal{O}_{K_1}(x_0, \gamma^{-n} \cdot x_0) \ .$$ 

\textbf{Proof.} Since, by definition, the sequence $(\gamma^n \cdot x_0)_{n \in \ZZ}$ is a quasi-geodesic one deduces from Inclusion \eqref{eqinclusionshadow} that
 $$	  \mathcal{O}_K(x_0, \gamma^{-n} \cdot x_0) \subset \mathcal{O}_{K+ \delta}(\gamma^n \cdot x_0,\gamma^{-n} \cdot x_0)  \ ,$$
where we choose $K$ satisfying both Lemma \ref{lem:complement_inclusion} and Lemma \ref{lem:inclusionshadow} and $n$ is sufficiently large. ($K$ only depends on the coefficients $(\lambda, C)$ of the quasi-geodesic $(\gamma^n \cdot x_0)_{n \in \ZZ}$.) Taking the complementary sets, we get that 
 $$	  \mathcal{O}_{K+ \delta}(x_0,\gamma^{-n} \cdot x_0) ^c \subset  \mathcal{O}_K(\gamma^{n} \cdot x_0, \gamma^{-n} \cdot x_0)^c \ .$$ 
 
We set $K_1 := K+ \delta$. Now we apply $\gamma^{2n}$ to get 
\begin{align*}	
\gamma^{2n} \left( \mathcal{O}_{K_1}( x_0, \gamma^{-n} \cdot x_0)^c \right) & \subset \gamma^{2n} \big( \mathcal{O}_K(\gamma^{n} \cdot x_0, \gamma^{-n} \cdot x_0)^c \big) \\
	& \subset  \mathcal{O}_K(\gamma^{3n} \cdot x_0, \gamma^{n} \cdot x_0)^c \ .
\end{align*}

We now use \eqref{eqlemballcomplement} (applied to  the reverse of $(\gamma^n \cdot x_0)_{n \in \ZZ}$) with $0\le n\le 3n$ to get that 
\begin{align*}
	\mathcal{O}_K(\gamma^{3n} \cdot x_0, \gamma^{n} \cdot x_0)^c & \subset  \mathcal{O}_K(x_0, \gamma^{n} \cdot x_0) \\
		& \subset  \mathcal{O}_{K_1}(x_0, \gamma^{n} \cdot x_0) \ ,
\end{align*}
concluding the proof.  \hfill $\blacksquare$ \\

Let $\gamma_1, \gamma_2$ be two loxodromic isometries as in the hypothesis. Fix $K_1$ satisfying Lemma \ref{lemloxo} for both $\gamma_1$ and $\gamma_2$ (notice that increasing $K_1$ does not affect the conclusion of the lemma).

\begin{lemma}\label{lem:shadow_position}
There exists $n>0$ so that the following hold.
\begin{enumerate}
    \item $\mathcal{O}_{K_1}(x_0, \gamma_1^{n} \cdot x_0)\cap\mathcal{O}_{K_1}(x_0,\gamma_2^{n} \cdot x_0)=\emptyset$,
    \item $\sup\{ (x,y)_{x_0} : x\in \mathcal{O}_{K_1}(x_0, \gamma_1^{n} \cdot x_0),y\in\mathcal{O}_{K_1}(x_0, \gamma_2^{n} \cdot x_0)\}<+\infty$,
    \item both items above also hold replacing "$n$" with "$-n$",
    \item the conclusion of Lemma \ref{lemloxo} holds for both $\gamma_1$ and $\gamma_2$, for the given $n$,
    \item $x_0\notin \mathcal{O}_{K_1}(x_0, \gamma_1^{-n} \cdot x_0)\cup\mathcal{O}_{K_1}(x_0,\gamma_2^{-n} \cdot x_0)$.
\end{enumerate}
\end{lemma}

\textbf{Proof.}
Since the quasigeodesic rays $(\gamma_1^n \cdot x_0)_{n\geq0}$ and $(\gamma_2^n \cdot x_0)_{n\geq0}$ have distinct endpoints at infinity, there exists $D$ so that for all $n,m\geq 0$ we have $(\gamma_1^n \cdot x_0,\gamma_2^m \cdot x_0)_{x_0}\leq D$. \\

For all sufficiently large $n$, we have $d(x_0,\gamma^n_i \cdot x_0)>D+K_1+2\delta$, for $i=1,2$. For $x\in \mathcal{O}_{K_1}(x_0, \gamma_1^{n}\cdot x_0)$ we claim that we have $(x,\gamma_2^n \cdot x_0)_{x_0}\leq D+\delta$.
Indeed, by definition of shadow we have
$$(\gamma_1^n \cdot x_0,x)_{x_0}\geq d(x_0,\gamma_1^n \cdot x_0)-K_1>D+2\delta,$$
and by hyperbolicity we have 
$$\min\{(x,\gamma_2^n \cdot x_0)_{x_0},(\gamma_1^n \cdot x_0,x)_{x_0}\}\leq (\gamma_1^n \cdot x_0,\gamma_2^n \cdot x_0)_{x_0}+\delta\leq D+\delta,$$
thereby showing the claim. In particular, $x\notin \mathcal{O}_{K_1}(x_0,\gamma_2^{n} \cdot x_0)$, since for any $z\in \mathcal{O}_{K_1}(x_0,\gamma_2^{n} \cdot x_0)$ we have $(z,\gamma_2^n \cdot x_0)_{x_0}> D+2\delta$ (we just did this computation for $\gamma_1$ above). This shows item 1. Now, if $y\in\mathcal{O}_{K_1}(x_0, \gamma_2^{n} \cdot x_0)$, then
$$\min\{(x,y)_{x_0},(\gamma_2^n \cdot x_0,y)_{x_0}\}\leq (x,\gamma_2^n \cdot x_0)_{x_0}+\delta\leq D+2\delta,$$
so that, in fact, we have $(x,y)_{x_0}\leq D+2\delta$. This shows item 2. Item 3 follows using the same arguments, again for any sufficiently large $n$. Up to increasing $n$, Lemma \ref{lemloxo} applies. For the last item, notice that we have $(x_0,x_0)_{\gamma_i^{-n} \cdot x_0}=d(x_0,\gamma_i^{-n} \cdot x_0)$, which is larger than $K_1$ for $n$ sufficiently large as above. \hfill $\blacksquare$ \\

Fix $n$ as in the previous lemma and denote $\mathcal O^\pm_i=\mathcal{O}_{K_1}(x_0, \gamma_i^{\pm n} \cdot x_0)$. We call a word in the alphabet $\{\gamma_1^{2n},\gamma_2^{2n}\}$ a positive word, while a negative word is a word in $\{\gamma_1^{-2n},\gamma_2^{-2n}\}$. In what follows we will conflate positive words and the corresponding group element. (A priori, different positive words might correspond to the same group element; we will deal with this later.) \\

For $w$ a positive word, denote $\mathcal O(w):=w \cdot \left( X- (\mathcal O^-_1\cup\mathcal O^-_2)\right)$.

\begin{lemma}\label{lem:positive_ping_pong}
For any integer $k$ there exists $D$ so that the following holds. If $w,w'$ are distinct positive words of the same length $k$ then $\mathcal O(w)\cap \mathcal O(w')=\emptyset$ and whenever $x\in \mathcal O(w)$ and $y\in \mathcal O(w')$, we have $(x,y)_{x_0}\leq D$.
\end{lemma}

\textbf{Proof.}
Consider distinct positive words $w,w'$ of length $k$. Up to swapping them, we can write them as $w=uv, w'=uv'$, where $v$ starts with $\gamma_1^{2n}$ and $v'$ starts with $\gamma_2^{2n}$ (and we allow $u$ to be empty).  By lemma \ref{lemloxo} and induction, we have $\mathcal O(v)\subseteq \mathcal O^+_1$ and $\mathcal O(v')\subseteq \mathcal O^+_2$, so that $\mathcal O(v)\cap \mathcal O(v')=\emptyset$ by Lemma \ref{lem:shadow_position}-(1). Since $\mathcal O(w)=u \cdot \mathcal O(v)$ and similarly for $w'$, we also have $\mathcal O(w)\cap \mathcal O(w')=\emptyset$, as required. \\

Consider now $x\in \mathcal O(w)$ and $y\in \mathcal O(w')$, so that $x=u\hat{x}$ for $\hat{x}\in \mathcal O(v)$, and similarly for $y$. Since $\mathcal O(v)\subseteq \mathcal O^+_1$ and $\mathcal O(v')\subseteq \mathcal O^+_2$, Lemma \ref{lem:shadow_position}-(2) implies $(\hat{x},\hat{y})_{x_0}\leq B$, where $B$ is the supremum in the statement. But then $(x,y)_{x_0}\leq B+d(x_0,u \cdot x_0)$, and the second term is bounded depending on $k$ only. This completes the proof of the lemma. \hfill $\blacksquare$ \\

Notice that the lemma implies that distinct words of the same length correspond to distinct group elements (since the $\mathcal O(w)$ are non-empty by Lemma \ref{lem:shadow_position}-(5)). Similar arguments as in the previous lemma also give

\begin{lemma}\label{lem:negative_ping_pong}
If $w,w'$ are distinct negative words of the same length then, for $i=1,2$, we have $w \cdot \mathcal O^-_i\cap w' \cdot \mathcal O^-_i=\emptyset$. 
\end{lemma}

We claim that the set $S$ of all (group elements corresponding to) positive words of length 7 is a Schottky set, where the constant $C$ is any constant larger than $D+\delta$ for $D$ as in Lemma \ref{lem:positive_ping_pong} with $k=4$. Let $y,z\in X$. For $i=1,2$, let $v_i$ be the positive word constructed as follows. If there is a positive word $v$ of length 3 so that $z\in v^{-1} \cdot \mathcal O^-_i$, then set $v=v_i$; note that there is at most one such word by Lemma \ref{lem:negative_ping_pong}. If there is no such $v$, choose any positive word of length 3 as $v_i$. We might have $v_1=v_2$. \\

Suppose by contradiction that at least one third of all $s\in S$ are so that $(y,s \cdot z)_{x_0}> C$. Then we have a subset $S'$ of $S$ with $\# S'\geq \# S/3$, that is, $\# S'\geq 43$, so that for any $s_1,s_2\in S'$ we have \begin{equation}
\label{eq.schottkyset.1}
    (s_1 \cdot z,s_2 \cdot z)_{x_0}\geq C-\delta \ , 
\end{equation}
by hyperbolicity. \\

From now and until the end of the proof, we refer to positive words of length 7 simply as words. Since there are at most $32 = 2 \times 2^4$ words ending with either $v_1$ or $v_2$, there must be at least 11 words which belong to $S'$ and not ending with $v_1$ or $v_2$. We are then left with a set of 11 words which do not end with $v_1$ or $v_2$ and such that Inequality \eqref{eq.schottkyset.1} holds for any pair of such words. Moreover, since there are at most 8 words that start with 4 given letters, out of these 11 words there must be 2 which have different initial subword of length 4. To sum up, we have shown so far that there are 2 words $\omega = u v$ and $\omega' = u' v'$ with $u \neq u'$ and $v,v' \notin \{v_1, v_2\}$. Let us come to the desired contradiction by showing that 
    $$ (\omega \cdot z,\omega' \cdot z)_{x_0} < C - \delta \ ,$$
contradicting \eqref{eq.schottkyset.1}. Indeed, since $v,v'\notin \{ v_1, v_2\}$, we have $v \cdot z,v' \cdot z\notin \mathcal O^-_1\cup\mathcal O^-_2$ by construction. Therefore, we have $\omega \cdot z = u v \cdot z \in \mathcal O(u)$ and $\omega' \cdot z \in \mathcal O(u')$ which, by Lemma \ref{lem:positive_ping_pong}, implies that $(s_1 \cdot z,s_2 \cdot z)\leq D$, a contradiction since we assumed $C>D+\delta$. \hfill $\blacksquare$

\section{Hamana's argument} 
\label{appendixA}

 We mainly repeat arguments from \cite{arthamadaupperdev} requiring only sub-additivity. Let $X$ be a metric space, $\mu$ a probability measure on $\Isom(X)$ with a finite exponential moment and $z_0 \in X$. Recall that we denoted by $z_n$ the position in $X$ at time $n$ of the random walk driven by $\mu$. By the triangle inequality we have for every $n,m \in \NN$. 
 	$$ d(z_0, z_{n + m}) \le d(z_0, z_n) + d(z_n, z_{n+m}) \ .$$
Recall that by sub-additivity and because $d(z_0, z_m)$ follows the same law as $d(z_n, z_{n+m})$ the following limit is well defined
	$$ l := \limi{n \to \infty} \ \frac{\EE(d_n)}{n} \ , $$
where we denoted $d_n := d(z_0,z_n)$. \\

Moreover, since we assumed that $\mu$ has a finite exponential moment and that the increments are I.I.D, there exists $\lambda_0 > 0$ such that for all $\lambda < \lambda_0$, one has 
 	$$ \EE(e^{\lambda d_{n + m}}) \le \EE(e^{ \lambda d_n}) \cdot \EE(e^{ \lambda d_m}) \ .$$
We conclude using the following purely analytical lemma.

\begin{lemma}\label{lemma.hamana}
Let $(d_n)_{n \in \NN}$ be a sequence of non-negative real valued random variables such that
\begin{itemize}
\item $d_1$ has a finite exponential moment;
\item there is $\lambda_0 > 0$ such that for any $0 \le \lambda < \lambda_0$ and for for any $m,n \in \NN$ one has $$\EE(e^{\lambda d_{m+n}}) \le \EE(e^{\lambda d_n}) \cdot \EE(e^{\lambda d_m}) \ . $$
\end{itemize}
Then for any 
	$$ a > l:=\limi{ n \to \infty} \left( \frac{\EE(d_n)}{n}\right)$$
one has 
 		$$ \limiinf{n \to \infty} \ \frac{ - \ln \left( \PP( d_n \ge a n ) \right) }{n} > 0 \ . $$
 \end{lemma}
 
 The range of validity of the above proposition is much wider than for random walks. It could be used in the setting of a sub-additive defective adapted cocycle as defined in \cite{artmathieusisto} for example. \\

\textbf{Proof.} First observe that the condition $\EE(e^{\lambda d_{m+n}}) \le \EE(e^{\lambda d_n}) \cdot \EE(e^{\lambda d_m})$ implies that 
$\EE(d_{n+m})\leq \EE(d_n)+\EE(d_m)$ and therefore that the limit defining $l$ does exist. \\ 

Let us introduce the notation 	$$ \Lambda_n(\lambda) := \ln \EE (e^{\lambda d_n}) \ . $$ 
Our two assumptions imply that there is $\lambda_0 > 0$ such that for all $ n \in \NN$ and all $\lambda < \lambda_0$ we have $\EE( e^{\lambda d_n}) < \infty$. Since $\EE(e^{\lambda d_{m+n}}) \le \EE(e^{\lambda d_n}) \cdot \EE(e^{\lambda d_m})$, we have 
	$ \Lambda_{n+m}(\lambda) \le \Lambda_n(\lambda) + \Lambda_m(\lambda)$, 		
which is to say that the sequence $( \Lambda_n(\lambda))_{n \in \NN}$ is sub-additive. Fekete's lemma implies that 
	$$ \frac{\Lambda_n(\lambda)}{n} \tends{n \to \infty } \Lambda(\lambda) := \infi{p \in \NN} \left( \frac{\Lambda_p(\lambda)}{p} \right) \ .$$

Using Markov's inequality we get that, for all $n$ and $\lambda>0$, 
$$ \PP ( d_n \ge  a n )  = \PP( e^{\lambda d_n} \ge e^{ \lambda a n} )  
	 		 \le e^{-\lambda a n} \ \EE ( e^{\lambda d_n}) \ . $$

Applying the logarithm and dividing by $\lambda n$ we get
$$
	 \frac{1}{\lambda} \frac{ \ln \Big( \PP ( d_n \ge  a n ) \Big) }{n} \le - a + \frac{ \Lambda_n(\lambda)}{\lambda n} \ . 
$$

Therefore, for all $\lambda > 0$, 
	$$	\limisup{n \to \infty} \ \frac{1}{\lambda} \frac{ \ln \Big( \PP ( d_n \ge  a n ) \Big) }{n} \le - a + \frac{\Lambda(\lambda)}\lambda \ .$$
	
It remains then to show that 
\begin{equation}\label{eq.it.remains.then.to}
\limisup{\lambda \to 0}  \left( \frac{\Lambda(\lambda)}{ \lambda } \right) \le l.
\end{equation}

At the cost of slightly reducing the value of $\lambda_0$,  one can suppose that $ \EE (d_n e^{\lambda_0 d_n}) < \infty$ for all $n > 0$. Because of the upper bound $e^x \le 1 + x + x^2e^x$, we have for all $\lambda < \lambda_0$.
\begin{align*}
	 \EE ( e^{\lambda d_n} ) & \le 1 + \lambda \ \EE ( d_n) + \lambda^2 \ \EE (d_n^2 \ e^{\lambda d_n}) \\
	 & \le 1 + \lambda \ \EE ( d_n) + \lambda^2 \ \EE (d_n^2 \ e^{\lambda_0 d_n}) \\
	 & \le 1 + \lambda \ \EE ( d_n) + \lambda^2 \ C_n  \ ,
\end{align*}
where $C_n := \EE (d_n^2 \ e^{\lambda_0 d_n})$. \\

Applying the logarithm, dividing both sides by $n$ and using the inequality $\ln(1+x) \le x$, we get
	$$ \frac{1}{n} \ln \left( \EE ( e^{\lambda d_n} ) \right) \le \lambda \  \frac{\EE ( d_n)}{n} + \frac{ C_n \lambda^2}{n} \ . $$
	
Therefore, for all $n > 0$ and for all $\lambda > 0$
	$$ \frac{\Lambda_n(\lambda)}{n} \le \lambda \  \frac{\EE ( d_n)}{n} + \frac{C_n \lambda^2}{n} \ . $$

In particular for all $ \lambda < \lambda_0$ and all $n \in \NN$ 
$$\Lambda(\lambda)= \infi{k \in \NN} \left( \frac{\Lambda_k(\lambda)}{  k} \right) \le \lambda \frac{\EE ( d_n)}{n} + \frac{C_n \lambda^2}{n} \ . $$

Letting $\lambda \to 0$ we deduce that for all $n \in \NN$
 $$ \limisup{\lambda \to 0} \left(  \frac{\Lambda(\lambda)}{ \lambda }  \right) \le \frac{\EE ( d_n)}{n} \ . $$
 
Finally taking $n$ to $\infty$ gives \eqref{eq.it.remains.then.to}. \hfill $\blacksquare$


\section{Properness and identification of the rate function} 
\label{appendix.rate}

Here, we show that the rate function appearing in Theorem \ref{maintheo} is proper and indicate a way of identifying the rate function as a Legendre transform of a generating function, under a stronger moment condition. These admit simple proofs and should be well-known to experts; however, we did not find an explicit reference and hence we indicate the argument for the convenience of the reader who may not be well-versed in large deviation theory. Finally, we give an explicit example of a rate function and suggest some open questions.

\subsubsection{Exponential tightness}

We show that the finite exponential moment assumption implies exponential tightness of the sequence $\frac{1}{n}d_n$ of random variables where, as before, we denote $d_n=d(z_n,z_0)$. We provide the proof for reader's convenience.

\begin{lemma}\label{lemma.proper.rate}
Let $\mu$ be a non-elementary probability measure on $\Isom(X)$ with a finite exponential moment. Then the sequence $\frac{1}{n} d_n$ is exponentially tight. 
\end{lemma}

\textbf{Proof.}
By Definition \ref{def.exp.tight} of exponential tightness, it suffices to show that 
\begin{equation*}
\lim_{t \rightarrow \infty} \limsup_{n \rightarrow \infty} \frac{1}{n} \ln \mathbb{P}\left(\frac{1}{n}d_n\geq t\right)=-\infty.
\end{equation*}
To see this, note that by Chebyshev inequality, for every $\lambda \geq 0$, we have
\begin{equation}\label{eq.proper}
\mathbb{P}( d_n \geq tn) \leq \mathbb{E}[e^{\lambda d_n}] e^{-\lambda tn}.
\end{equation} Using finite exponential moment assumption, let $\lambda_0>0$ be such that $\mathbb{E}[e^{\lambda_0 d_1}]<\infty$. In \eqref{eq.proper}, taking logarithm, dividing by $n$ and specializing to some $\lambda_{1}>0$ such that $\lambda_0 \geq \lambda_{1}$, we get
\begin{equation*}
\frac{1}{n} \ln \mathbb{P}(d_n \geq tn) \leq -(\lambda_{1}t - \frac{1}{n}\ln \mathbb{E}[e^{\lambda_{1} d_n}])
\end{equation*} 
On the other hand, it follows by the independence of random walk increments and the subadditivity that for all $n  \geq 1$, we have $\frac{1}{n} \ln \mathbb{E}[e^{\lambda_{1} d_n}] \leq \ln \mathbb{E}[e^{\lambda_{1} d_1}]$. Therefore, we obtain $$
\limsup_{n \rightarrow \infty} \frac{1}{n} \ln \mathbb{P}\left(\frac{1}{n}d_n\geq t\right) \leq -(\lambda_{1}t-\mathbb{E}[e^{\lambda_{1}d_1}])$$ Since $\mathbb{E}[e^{\lambda_{1}d_1}]<\infty$ and $\lambda_1>0$, the result follows by taking the limit as $t \to \infty$. \hfill $\blacksquare$

\subsubsection{Identification of the rate function}
In this part, let $\mu$ be a non-elementary probability measure on $\Isom(X)$ which has \textit{strong exponential moment} in the sense that $\mathbb{E}[\exp(\alpha d_n)]<\infty$ for every $\alpha \geq 0$. Note that clearly, a probability measure $\mu$ of bounded support has strong exponential moment. The limit Laplace transform of the sequence $(d_n)$ is the function $\Lambda: \mathbb{R} \to [0,\infty)$ defined by 
\begin{equation*}
\Lambda(\lambda) = \lim_{n \rightarrow \infty}\frac{1}{n}\ln \mathbb{E}[e^{\lambda d_n}].
\end{equation*} 
This function already appeared in the proof of Lemma \ref{lemma.hamana}. As mentioned there, for the random variables $(d_n)$, this limit exists by subadditivity arguments without appeal to LDP. More generally, provided that the sequence $\frac{1}{n}d_n$ satisfies a LDP with convex rate function, this convergence is also a consequence of Varadhan's integral lemma (see \cite{dembo-zeitouni} section 4.3), which, moreover, identifies the limit. In the other direction, we note that nice analytic properties (e.g.\ differentiability, steepness) of this function have direct implications for the LDP (see e.g.\ G\"{a}rtner-Ellis theorem \cite[\S 4]{dembo-zeitouni}). 

The following lemma gives an expression of the rate function appearing in Theorem \ref{maintheo} under strong exponential moment assumption.

\begin{lemma}\label{lemma.identify}
Let $\mu$ be a non-elementary probability measure on $\Isom(X)$ with a strong exponential moment. Let $I:[0,\infty) \to [0,\infty]$ be the rate function given by Theorem \ref{maintheo} and let $\lambda:\mathbb{R} \to [0,\infty)$ be the limit Laplace transform of the sequence $(d_n)$. Then, for every $\lambda \in [0,\infty)$, we have 
$$
I(\lambda)=\sup_{\alpha \in \mathbb{R}}(\lambda \alpha- \Lambda(\alpha)).
$$
\end{lemma}

\textbf{Proof.}
Thanks to the strong exponential moment assumption, for every $\lambda \in \mathbb{R}$, the functional $x \mapsto \lambda x$ composed with $d_n$ satisfies the moment assumption of Varadhan's integral lemma (see \cite[(4.3.3)]{dembo-zeitouni}). Therefore, for every $\lambda \in \mathbb{R}$, we have 
\begin{equation*}
\Lambda(\lambda)=\lim_{n \to \infty} \frac{1}{n} \ln \mathbb{E}[e^{\lambda d_n}] =\underset{ \alpha \in \mathbb{R}}{\sup}(\lambda \alpha-I(\alpha))
\end{equation*} where $I$ is the proper rate function of the LDP of the sequence $(\frac{1}{n}d_n)$.\\

For a function $f$ on $\mathbb{R}$, denote its convex conjugate (Legendre transform) by $f^{\ast}(.)$, where $f^{\ast}(\lambda):= \sup_{\alpha \in \mathbb{R}}(\lambda \alpha -f( \alpha)) $. The above conclusion of Varadhan's integral lemma reads as $\Lambda(\lambda)=I^{\ast}(\lambda)$. Now, since $I$ is a convex rate function, Fenchel--Moreau duality tells us that $I(\alpha)=I^{\ast \ast}(\alpha)=\Lambda^{\ast}(\alpha)$, identifying $I(\alpha)$ with $\Lambda^{\ast}(\alpha)$ and completing the proof. \hfill $\blacksquare$\\

Let us finish with an example of a rate function that one can obtain using the previous lemma, and some questions.\\

It is not difficult to pinpoint the explicit expression of the rate function for the standard random walk on the free group $F_q$ of rank $q \geq 1$. It is given by the following 

\begin{equation*}
I(\alpha)= \begin{cases} \frac{1+\alpha}{2}\ln(1+\alpha) + \frac{1-\alpha}{2} \ln(1-\alpha) +\ln (q) - \frac{1+\alpha}{2} \ln(2q-1) \qquad &\alpha \in [0,1]\\ \infty &\text{otherwise}
\end{cases}
\end{equation*}
We remark that, among others, this function satisfies the following properties:\\[2pt]
\indent 1) it is analytic and strictly convex on its effective support,\\[2pt]
\indent 2) $I(0)=-\ln \frac{\sqrt{2q-1}}{q}$ where $\frac{\sqrt{2q-1}}{q}$ is the spectral radius of the standard random walk on $F_q$ calculated by Kesten \cite{kesten.symmetric},\\[2pt]
\indent 3) the drift $\frac{q-1}{q}$ is the unique zero of $I$,\\[2pt]
\indent 4) if $\Lambda(\lambda)$ denotes the Legendre transform of $I$ given by $\Lambda(\lambda)=\sup_{\alpha \in \mathbb{R}} (\lambda \alpha-I(\alpha))$, then $\Lambda''(0)-(\frac{q-1}{q})^2$ is the variance appearing in the central limit theorem for the standard random walk on the free group (this fact can be deduced either directly or as in \cite[Lemma 5.2]{bougerol-lacroix}).\\[4pt]
Whereas finding an explicit expression for the rate function $I$ in Theorem \ref{maintheo} does not seem to be feasible in general, pinning down some of its general properties, paralleling the above ones, is a more tractable challenge. As we showed, the property 3)  holds under very general assumptions, and it is not hard to see that the same is true of 2). In turn, the properties 1) and 4) naturally suggest the corresponding open problems. We mention only a few of them:\\

\textbf{Question C.1.} Is the rate function appearing in Theorem \ref{maintheo} strictly convex? Analytic? Do these properties depend on generating set or probability measure?

\bibliographystyle{alpha}

\bibliography{bibliography} 

\end{document}

%% file: dessinssistohyp.pdf_tex
\begingroup%
  \makeatletter%
  \providecommand\color[2][]{%
    \errmessage{(Inkscape) Color is used for the text in Inkscape, but the package 'color.sty' is not loaded}%
    \renewcommand\color[2][]{}%
  }%
  \providecommand\transparent[1]{%
    \errmessage{(Inkscape) Transparency is used (non-zero) for the text in Inkscape, but the package 'transparent.sty' is not loaded}%
    \renewcommand\transparent[1]{}%
  }%
  \providecommand\rotatebox[2]{#2}%
  \ifx\svgwidth\undefined%
    \setlength{\unitlength}{1190.5511811bp}%
    \ifx\svgscale\undefined%
      \relax%
    \else%
      \setlength{\unitlength}{\unitlength * \real{\svgscale}}%
    \fi%
  \else%
    \setlength{\unitlength}{\svgwidth}%
  \fi%
  \global\let\svgwidth\undefined%
  \global\let\svgscale\undefined%
  \makeatother%
  \begin{picture}(1,0.5)%
    \put(0,0){\includegraphics[width=\unitlength,page=1]{dessinssistohyp.pdf}}%
    \put(0.11519275,0.33657881){\color[rgb]{0,0,0}\makebox(0,0)[lb]{\smash{$x_{k_1}$}}}%
    \put(0.11339287,0.40857428){\color[rgb]{0,0,0}\makebox(0,0)[lb]{\smash{$x_{k_1 +1}$}}}%
    \put(0.42657314,0.46977041){\color[rgb]{0,0,0}\makebox(0,0)[lb]{\smash{$x_k$}}}%
    \put(0.34107848,0.35457767){\color[rgb]{0,0,0}\makebox(0,0)[lb]{\smash{$x_{i+1}$}}}%
    \put(0.27808246,0.43737248){\color[rgb]{0,0,0}\makebox(0,0)[lb]{\smash{$x_i$}}}%
    \put(0.44187219,0.35637752){\color[rgb]{0,0,0}\makebox(0,0)[lb]{\smash{}}}%
    \put(0,0){\includegraphics[width=\unitlength,page=2]{dessinssistohyp.pdf}}%
    \put(0.39768183,0.07434445){\color[rgb]{0,0,0}\makebox(0,0)[lb]{\smash{$\pi_{\nu}(x_k)$}}}%
    \put(0.0827263,0.07127184){\color[rgb]{0,0,0}\makebox(0,0)[lb]{\smash{$\pi_{\nu}(x_{k_1})$}}}%
    \put(0,0){\includegraphics[width=\unitlength,page=3]{dessinssistohyp.pdf}}%
    \put(0.80904785,0.37977606){\color[rgb]{0,0,0}\makebox(0,0)[lb]{\smash{$x_{k_2 +1}$}}}%
    \put(0.70195581,0.44637188){\color[rgb]{0,0,0}\makebox(0,0)[lb]{\smash{$x_{k_2}$}}}%
    \put(0.76617333,0.07127187){\color[rgb]{0,0,0}\makebox(0,0)[lb]{\smash{$\pi_{\nu}(x_{k_2 +1})$}}}%
    \put(0,0){\includegraphics[width=\unitlength,page=4]{dessinssistohyp.pdf}}%
    \put(0.85353757,0.13496047){\color[rgb]{0,0,0}\makebox(0,0)[lb]{\smash{$\rho$}}}%
    \put(0.01781798,0.11581682){\color[rgb]{0,0,0}\makebox(0,0)[lb]{\smash{$\nu$}}}%
  \end{picture}%
\endgroup%

%% file: lemmegeomsimple.pdf_tex
\begingroup%
  \makeatletter%
  \providecommand\color[2][]{%
    \errmessage{(Inkscape) Color is used for the text in Inkscape, but the package 'color.sty' is not loaded}%
    \renewcommand\color[2][]{}%
  }%
  \providecommand\transparent[1]{%
    \errmessage{(Inkscape) Transparency is used (non-zero) for the text in Inkscape, but the package 'transparent.sty' is not loaded}%
    \renewcommand\transparent[1]{}%
  }%
  \providecommand\rotatebox[2]{#2}%
  \ifx\svgwidth\undefined%
    \setlength{\unitlength}{841.88976378bp}%
    \ifx\svgscale\undefined%
      \relax%
    \else%
      \setlength{\unitlength}{\unitlength * \real{\svgscale}}%
    \fi%
  \else%
    \setlength{\unitlength}{\svgwidth}%
  \fi%
  \global\let\svgwidth\undefined%
  \global\let\svgscale\undefined%
  \makeatother%
  \begin{picture}(1,0.70707071)%
    \put(0,0){\includegraphics[width=\unitlength,page=1]{lemmegeomsimple.pdf}}%
    \put(0.07595958,0.29919994){\color[rgb]{0,0,0}\makebox(0,0)[lb]{\smash{$z$}}}%
    \put(0.66698859,0.62413975){\color[rgb]{0,0,0}\makebox(0,0)[lb]{\smash{$q$}}}%
    \put(0.90454338,0.17522613){\color[rgb]{0,0,0}\makebox(0,0)[lb]{\smash{$x$}}}%
    \put(0.45051708,0.43906323){\color[rgb]{0,0,0}\makebox(0,0)[lb]{\smash{$y$}}}%
    \put(0,0){\includegraphics[width=\unitlength,page=2]{lemmegeomsimple.pdf}}%
    \put(0.30798057,0.24562088){\color[rgb]{0,0,0}\makebox(0,0)[lb]{\smash{$\ge R$}}}%
    \put(0.25198413,0.3856121){\color[rgb]{0,0,0}\makebox(0,0)[lb]{\smash{$\le R$}}}%
    \put(0.17689794,0.35252319){\color[rgb]{0,0,0}\makebox(0,0)[lb]{\smash{$4R/5 \le$}}}%
    \put(0.46451623,0.36524974){\color[rgb]{0,0,0}\makebox(0,0)[lb]{\smash{$\le R/5$}}}%
  \end{picture}%
\endgroup%